%% file: lekspde.tex
\documentclass[12pt,a4paper]{article}

\title{Lyapunov exponents of the Kuramoto--Sivashinsky PDE}

\author{Russell A. Edson
\thanks{School of Mathematical Sciences, University of Adelaide, South Australia.  \protect\url{mailto:russell.edson@adelaide.edu.au} ,
\protect\url{https://orcid.org/0000-0002-4607-5396}}
\and  J.~E. Bunder
\thanks{School of Mathematical Sciences, University of Adelaide, South Australia.  \protect\url{mailto:judith.bunder@adelaide.edu.au} ,
\protect\url{http://orcid.org/0000-0001-5355-2288}}
\and  Trent W. Mattner
\thanks{School of Mathematical Sciences, University of Adelaide, South Australia.  \protect\url{mailto:trent.mattner@adelaide.edu.au} ,
\protect\url{https://orcid.org/0000-0002-5313-5887}}
\and A.~J. Roberts
\thanks{School of Mathematical Sciences, University of Adelaide, South Australia.  \protect\url{mailto:anthony.roberts@adelaide.edu.au} ,
\protect\url{http://orcid.org/0000-0001-8930-1552} }
}

\date{\today}

\usepackage{algorithm}
\usepackage{algorithmicx}
\usepackage{algpseudocode}
\algnewcommand\algorithmicinput{\textsc{Input:}}
\algnewcommand\Input{\item[\algorithmicinput]}
\algnewcommand\algorithmicoutput{\textsc{Output:}}
\algnewcommand\Output{\item[\algorithmicoutput]}
\algnewcommand\algorithmicbreak{\textbf{break}}
\algnewcommand\Break{\State \algorithmicbreak}

\usepackage{amsmath,defns,amssymb,eulervm,reducecode,doi}
\allowdisplaybreaks

\IfFileExists{jeb.sty}{\usepackage{jeb}}{\usepackage{ajr}}

\usepackage[capitalise,nameinlink,noabbrev]{cleveref}
\crefname{equation}{}{}
\let\eqref\cref
\let\theorem\relax

\makeatletter\let\c@theorem\relax\makeatother
\theoremstyle{theorem} %
\newtheorem{theorem}{Theorem}
\newtheorem{lemma}[theorem]{Lemma}
\newtheorem{proposition}[theorem]{Proposition}
\newtheorem{assumption}[theorem]{Assumption}
\newtheorem{definition}[theorem]{Definition}
\atBegin{definition}{\bigskip}
\atBegin{theorem}{\bigskip}
\atBegin{proposition}{\bigskip}
\atBegin{lemma}{\bigskip}
\atBegin{assumption}{\bigskip}

\usepackage[backend=bibtex8 %
    , style=authoryear %
    ,backref=true %
    ,maxbibnames=12 %
    ,giveninits=true %
    ,isbn=false
    ,eprint=true
    ]{biblatex}
\bibliography{bibexport}
\AtEndDocument{{\raggedright\printbibliography}}
\let\citet\textcite
\let\citep\parencite
\makeatletter%
\def\cite{\@ifnextchar[{\parencite}{\textcite}}
\makeatother
\DeclareFieldFormat{url}{%
  \iffieldundef{doi}{%
  \mkbibacro{URL}\addcolon\space\url{#1}}{}}
\DeclareFieldFormat{urldate}{%
  \iffieldundef{doi}{%
  \mkbibparens{\bibstring{urlseen}\space#1}}{}}

\usepackage{pgfplots}
\pgfplotsset{compat=newest}
\usetikzlibrary{decorations.pathreplacing}
\usetikzlibrary{shapes.geometric}
\usepgfplotslibrary{external} %
\newcommand{\myinput}[1]{}%

\newcommand{\kspde}{Kuramoto--Sivashinsky \textsc{pde}~\eqref{eqn:kspde}}
\Bb R
\Vec u \Vec w\Vec q

\iftrue
\usepackage{pdfcomment}

\else
\fi

\begin{document}

\maketitle

\begin{abstract}
The Kuramoto--Sivashinsky equation is a prototypical chaotic nonlinear partial differential equation (\textsc{pde}) in which the size of the spatial domain plays the role of a bifurcation parameter.
We investigate the changing dynamics of the Kuramoto--Sivashinsky \pde\ by calculating the Lyapunov spectra over a large range of domain sizes.
Our comprehensive computation and analysis of the Lyapunov exponents and the associated Kaplan--Yorke dimension provides new insights into the chaotic dynamics of the Kuramoto--Sivashinsky \pde, and the transition to its 1D turbulence.
\end{abstract}

\tableofcontents

\section{Introduction}

The Kuramoto--Sivashinsky \pde~\cref{eqn:kspde} models a wide variety of nonlinear systems with intrinsic instabilities, such as wave propagation in chemical reaction-diffusion systems \citep{Kuramoto76}, the velocity of laminar flame front instabilities \citep{Sivashinsky77}, thin fluid film flow down inclined planes \citep{Sivashinsky80}, and hydrodynamic turbulence \citep[e.g.,][]{Pomeau85, Hohenberg89, Dankowicz96}.
In the Kuramoto--Sivashinsky \pde\ the large scale dynamics are dominated by a destabilising `diffusion', whereas small scale dynamics are dominated by stabilising hyperdiffusion, and a nonlinear advective term stabilises the system by transferring energy from the large unstable modes to the small stable modes \citep[e.g.,][p.~199]{Sprott10}.
The interplay between these contrasting features leads to significant spatio-temporal complexity \citep[e.g.,][]{Hyman86,Cross93,Cvitanovic10}, from intermittent disorder, through to chaos, hyperchaos and~turbulence. Lyapunov exponents characterise this chaos and turbulence \citep[e.g.,][]{Eckmann85,Ruelle71,Takens81}, and are increasingly
used to analyse such spatio-temporal complexity in various applications
such as turbulent Poiseulle flow \cite[]{Keefe1992}, turbulence in 
flames \cite[]{Hassanaly2018}, and Rayleigh--B{\'e}rnard fluid convection \cite[]{Chertovskih2015, MuXu2017}.
Here we show new details of how the dynamics of the Kuramoto--Sivashinsky \pde\ become increasingly chaotic as the size of the domain increases, for both periodic and odd-periodic boundary conditions. 
To measure the degree of chaos, \cref{sec:ele} computes the Lyapunov exponents using the classic algorithm introduced by \citet{Benettin80b} and \citet{Shimada79}, but now in new detail over a comprehensive range of domain sizes.
By comparison, \citet{Tajima02} explored the Kuramoto--Sivashinsky \pde~\eqref{eqn:kspde} with rigid boundary conditions over a range of domain lengths, whereas we explore periodic~\eqref{eq:pbc} and odd-periodic~\eqref{eq:opbc} cases, we use an order of magnitude increased resolution in the domain lengths, and we also cover the transition to chaos regime.
The Lyapunov exponents describe the rate at which neighbouring trajectories diverge under a chaotic flow, and thus provide a quantitative measure of the degree of chaos in a system.
\cref{sec:ele} analyses the growth of the Lyapunov exponents with increasing domain size, and then uses the Lyapunov spectra to identify the onset of chaos and to characterise new details of the increasingly complex spatio-temporal dynamics of the Kuramoto--Sivashinsky~\pde.

A further use of the Lyapunov exponents is in defining the Kaplan--Yorke dimension of the attractor of a dynamical system~\citep{Kaplan79}.
The Kaplan--Yorke dimension bounds above the fractal dimension of the chaotic attractor, and approximates the number of effective modes necessary to describe the dynamics on the attractor~\citep{Grassberger83}.
For a Kuramoto--Sivashinsky \pde\ defined on either a periodic or odd-periodic domain, \cref{sec:ckyd} confirms more accurately how the Kaplan--Yorke dimension scales roughly linearly with the domain size. 
This linear scaling corresponds well with the scaling observed by \citet{Manneville85} and \citet{Tajima02} for the Kuramoto--Sivashinsky \pde\ with rigid boundary conditions.

Our detailed analysis of the chaotic dynamics of the Kuramoto--Sivashinsky \pde\ and its dependence on domain size provides new insights into the onset of chaos. 
In many chaotic systems, we can identify that discrete point at which a bifurcation parameter permits chaos, but our new visualisation of the comprehensive computation of Lyapunov exponents highlights the gradual changes which drive a system into the chaotic regime, and thence into 1D~turbulence.

\section{The Kuramoto--Sivashinsky equation}
\label{sub:ckspde-pde}

On the spatial domain~$0 \leq x \leq L$ for some domain size~$L$, the one-dimensional Kuramoto--Sivashinsky~\textsc{pde} for field~$u(x,t)$ is
\begin{equation}
  \partial_t u + \partial^4_x u  + \partial^2_x u  + u \partial_x u= 0\,.
  \label{eqn:kspde}
\end{equation}
We apply either periodic boundary conditions,
\begin{equation}
u(x+L,t)=u(x,t)\quad\text{for all }0\leq x\leq L\,,  
\label{eq:pbc}
\end{equation}
or odd-periodic boundary conditions, 
\begin{equation}
u(x,t)=\partial_x^2u(x,t)=0\quad \text{at }x=0,L\,.
\label{eq:opbc}
\end{equation}
In the Kuramoto--Sivashinsky~\textsc{pde}, the second order diffusive term $\partial^2_x u$ is destabilising whereas the fourth order hyperdiffusion term $\partial^4_x u$ is stabilising, resulting in large scale instabilities and small scale dissipation, with the transfer of energy from large to small scales, mediated by the nonlinear term~$u \partial_x u$, producing a stabilising influence on the system \citep[e.g.,][p.~199]{Sprott10}.

The  \kspde\ supports several symmetries, although in the turbulent regime they only hold in a time-averaged sense~\citep[e.g.,][]{Wittenberg99, Cvitanovic10}.
Of particular interest for periodic domains~\eqref{eq:pbc} are the Galilean invariance, $u(x,t) \to u(x-ct,t)+c$ for all speeds~$c$, and the spatial translation invariance, $u(x,t) \to u(x+d,t)$ for all~$d$.
These two symmetries do not hold for odd-periodic domains~\eqref{eq:opbc}.
The \kspde\ with odd-periodic domains~\eqref{eq:opbc} is particularly well studied~\citep[e.g.,][]{Rempel04, Lan08, Foias86, Eguiluz99} compared to that with periodic domains, as the removal of periodic symmetries simplifies somewhat the analysis of the dynamics. 
The relatively simpler dynamics of the odd-periodic case is observed in \cref{sec:ele} when comparing the periodic and odd-periodic Lyapunov spectra.
The most obvious point of difference is that the Galilean and translation invariances support two zero Lyapunov exponents which are absent from the Lyapunov spectra for the odd-periodic case.
Commonly for the periodic case, a zero mean condition is imposed to remove the Galilean invariance and consequently one of the zero Lyapunov exponents~\citep[e.g.,][]{Cvitanovic10, Wittenberg99, Dankowicz96}; we do not impose the zero mean condition.

\begin{figure}
  \centering
  \tikzsetfigurename{ksoddperiodic_plots}
  \input{ksoddperiodic_plots.tex}
  \caption{Simulations of the \kspde+\eqref{eq:opbc} depend upon the size of the odd-periodic spatial domain~$L$: (top left)~for $L=17.5$, an oscillating cell emerges; (top right)~for $L=18.2$\,, intermittent bursts disrupt the cell structure; (bottom left)~for $L=41$, chaotic cellular structures criss-cross and interact; and (bottom right)~for $L=100$\,, spatio-temporally complex patterns of `turbulence'.}
  \label{fig:cksodd-plots}
\end{figure}
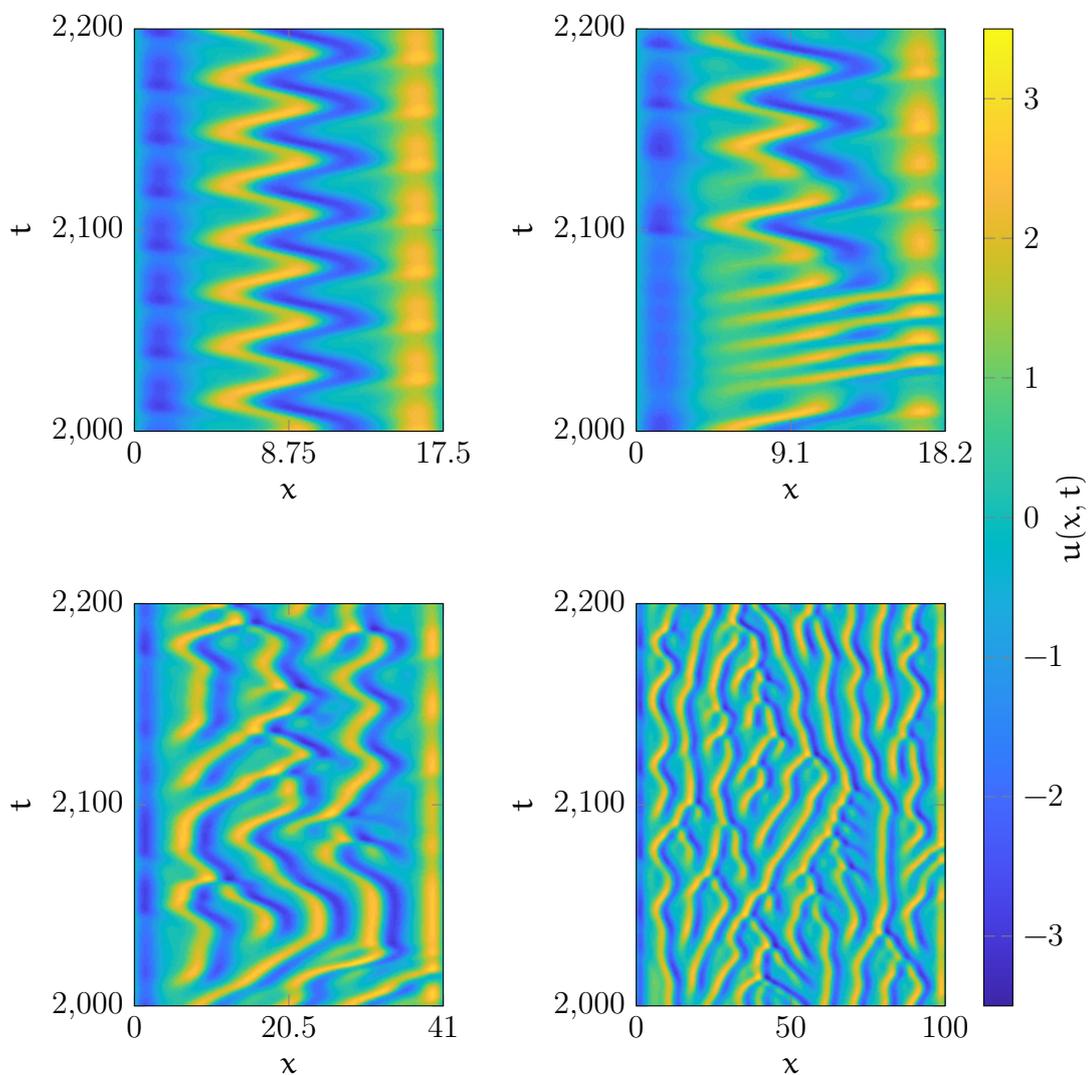

As the size~\(L\) of the spatial domain varies, the \kspde\ produces distinctly different dynamics \citep[e.g.,][]{Hyman86, Cross93, Cvitanovic10}.
For periodic boundary conditions~\eqref{eq:pbc}, \cref{fig:ckspde-plots} shows the increasing complexity of the Kuramoto--Sivashinsky dynamics as domain size~$L$ increases, from stable travelling wave when $L\lesssim 13$ through to a spatio-temporal chaotic turbulence when $L\approx100$\,.
\cref{fig:cksodd-plots} for odd-periodic boundary conditions~\eqref{eq:opbc} also shows an increase in the complexity of the dynamics as~$L$ increases, progressing from an oscillating cell when $L\lesssim 17$ through to a spatio-temporal turbulence when $L\approx100$\,.  
For both types of boundary conditions, the spatial domain size~$L$ plays the role of a bifurcation parameter.

\begin{figure}
  \centering
  \tikzsetfigurename{ksperiodic_plots}
  \input{ksperiodic_plots.tex}
  \caption{Simulations of the \kspde+\eqref{eq:pbc} depend upon the size of the periodic spatial domain~$L$: (top left)~for $L=12$\,, a travelling wave emerges; (top right)~for $L=13.5$\,,  intermittent bursts disrupt the travelling wave structure; (bottom left)~for $L=36$\,, chaotic cellular structures criss-cross and interact; and (bottom right)~for $L=100$\,, spatio-temporally complex patterns of `turbulence'.}
  \label{fig:ckspde-plots}
\end{figure}

Our aim is to provide new details of the character of the trend to spatio-temporal complexity and `turbulence' with increasingly long domains~$L$. 
To do this,  \cref{sec:ele} comprehensively computes the $24$~most positive Lyapunov exponents across a significant range of domain lengths~$L$.
It is these `most positive' Lyapunov exponents which determine the nature of the chaotic dynamics.
Then, \cref{sec:ckyd} evaluates the Kaplan--Yorke dimension over the same range of~$L$ to show more information about how the dimension of the chaotic attractor grows linearly with~$L$.

\section{Evaluate the Lyapunov exponents}
\label{sec:ele}

In a dynamical system, the Lyapunov exponents measure the exponential divergence of initially close trajectories, with positive Lyapunov exponents indicating divergent trajectories and negative Lyapunov exponents indicating convergence~\citep{Chicone06,Eckmann85}.
A chaotic system, due to its sensitivity to initial conditions, must have at least one positive Lyapunov exponent. Furthermore, an increasingly chaotic system has an increasing number of positive Lyapunov exponents.
This section evaluates Lyapunov exponents of the \kspde\ for increasing domain size~$L$ and interprets each increase in the number of positive Lyapunov exponents as a transition to an increasingly chaotic system.
 
Formally,  Lyapunov exponents measure trajectory divergences in the infinite time limit, with different Lyapunov exponents corresponding to divergences in different orthogonal directions.
The existence of these time limits are assured (almost everywhere) by Oseledec's multiplicative ergodic theorem~\citep[e.g.,][]{Eckmann85,Ruelle79}.
In numerical calculations of Lyapunov exponents, complications due to the infinite time limit are avoided by computing $N$~iterations of the divergence of  trajectories over finite time intervals~$T$, with $N$~large but finite~\citep[e.g.,][]{Shimada79, Benettin1980, Geist90, Dieci97, Skokos2010}.
At the end of each iteration, the divergent trajectories are reorthonormalised.
This reorthonormalisation ensures the tracked directions remain orthogonal, rather than all converging to that of the largest positive Lyapunov exponent~\citep[e.g.,][]{Geist90}. 
This rescaling is valid in the ergodic case because the Lyapunov exponents are (almost everywhere) independent of a trajectory's initial condition. However, the finite time numerical approximation generally results in some numerical error~\citep[e.g.,][]{Dieci97}.

\cref{alg:clyap-lyaps} assumes a vector function of time, $\uv (t)\in\RR^n$\,, satisfies the dynamical system
\begin{equation}
\dot{\uv }=\vec{f}(t,\uv )\,,\label{eq:genu}
\end{equation}
with initial condition~$\uv (0)$.
Trajectories are first evolved for time~$\tau$ to ensure initial transients have decayed and thus that the system is close to an attractor.
Then \cref{alg:clyap-lyaps} numerically solves the \ode~\eqref{eq:genu} for a time~\(N\cdot T\) to compute the $m$~most positive Lyapunov exponents~$\lambda_i$ for $i = 1,\ldots,m\leq n$ using reduced \textsc{qr}~decomposition to reorthonormalise after each of \(N\)~time intervals of length~\(T\) \citep[e.g.,][]{Shimada79, Benettin1980, Geist90, Dieci97, Skokos2010}.
As is standard, the resulting Lyapunov exponents are ordered such that $\lambda_1\geq\lambda_2\geq\ldots\geq \lambda_m$\,.

\begin{algorithm}
  \caption{The classic algorithm for computing the spectrum of 
  Lyapunov exponents for a dynamical system, introduced by 
  \protect\citet{Benettin80b}, and \protect\citet{Shimada79}.}
  \label{alg:clyap-lyaps}
  \begin{algorithmic}[1]
      \Statex $d\uv /dt = \vec{f}(t,\uv )$: the dynamical 
      system \ode
      \Statex $\uv(0)$: the initial value of~\uv
      \Statex $m$: the number of the most positive exponents 
      to compute
      \Statex $\tau$: time to simulate system before computing
      exponents
      \Statex $T$: time between reorthonormalisation steps
      \Statex $N$: the total number of reorthonormalisation steps
      \Statex $\epsilon$: perturbation magnitude (typically take 
        $\epsilon = 10^{-6}$).

    \Output
      \Statex $\lambda_i$: the $m$ most positive Lyapunov exponents,
      $i=1,\ldots,m$.

    \phantom{text}

    \State compute \(\uv ^{(0)} := \uv(\tau)\) via solving \ode\ on \([0,\tau]\)
    \State set~$t_j := \tau+jT$, for~$j=0,1,2,\ldots,N$
    \State choose initial orthogonal directions 
    $Q^{(0)} := \begin{bmatrix}
    \qv _1^{(0)} & \cdots & \qv _m^{(0)}
    \end{bmatrix}$

    \For{$j=1:N$}
      \State compute $\uv ^{(j)} := \uv(t_j)$ via solving \ode\ with \(\uv(t_{j-1})=\uv^{(j-1)}\)  
      \For{$i=1:m$}
      \State compute $\wv_i ^{(j)} := \uv(t_j)$ via \ode\ with \(\uv(t_{j-1})=\uv ^{(j-1)} + 
        \epsilon \qv _i^{(j-1)}\)  
        \State approximate $\Psi(t_j,t_{j-1})\qv _i^{(j-1)}
        := (\wv _i^{(j)}-\uv ^{(j)})/\epsilon$
      \EndFor  
      \State construct $\Psi(t_j,t_{j-1}) Q^{(j-1)} := 
      \begin{bmatrix} \Psi(t_j,t_{j-1})\qv _1^{(j-1)} & 
      \cdots & \Psi(t_j,t_{j-1})\qv _m^{(j-1)} \end{bmatrix}$
      \State compute $Q^{(j)}R^{(j)} := \,\,$\textsc{qr}$\left(
      \Psi(t_j,t_{j-1})Q^{(j-1)}\right)$
    \EndFor  
      \For{$i=1:m$}
      \State compute~$\lambda_i := \sum^N_{j=1}\log R_{i,i}^{(j)}/(NT)$
      \EndFor  
  \State \Return $\{\lambda_i\}$.
\end{algorithmic}
\end{algorithm}

In implementing \cref{alg:clyap-lyaps} for the \kspde\ a \(n\)-D approximate system is used, either spectral in space for the periodic case~\eqref{eq:pbc} or finite differences for the odd-periodic case~\eqref{eq:opbc}. 
In either case we choose truncations so that the maximum wavenumber resolved was $k_{\max}\approx9$ (which decays extremely rapidly, on a time scale of \(1/k_{\max}^4\approx 10^{-4}\)). 
Initial conditions were random and normally distributed$\sim\mathcal{N}(0,1)$.
After some testing of different transient times~$\tau$, we selected $\tau=2000$, which is  smaller than some other studies (e.g., \citet{Wittenberg99} used $\tau=100\,000$), but repeatedly provided consistent and expected dynamics. A total of~$N=1000$ reorthonormalisation steps are performed in the exponent computations.
This choice of~$N$ computes fairly accurate Lyapunov exponents within a reasonable time frame.

\begin{figure}
  \centering
\begin{itemize}
\item  domain size \(L=20\)\\
\tikzsetfigurename{ksoddperiodic_lyaps_L20_100varyT_scaled}
  \input{ksoddperiodic_lyaps_L20varyT_scaled.tex}
\item domain size \(L=100\)\\
    \tikzsetfigurename{ksoddperiodic_lyaps_L100varyT}
  \input{ksoddperiodic_lyaps_L100varyT.tex}

\end{itemize}
  \caption{The $m=24$ most positive Lyapunov exponents~$\lambda_i$ for (top)~$L=20$ and (bottom)~$L=100$ calculated using \cref{alg:clyap-lyaps} with odd-periodic boundary conditions~\eqref{eq:opbc} and different interval times~$T$.
    An accurate Lyapunov exponent is approximately constant as~$T$ varies. When $L=20$ there is no range of~$T$ for which \emph{all} $24$~exponents are constant, but the most negative ones are least important. 
    For $L=100$\,, all $m=24$ Lyapunov exponents are reasonably constant for $1<T<10$\,.
  }
  \label{fig:ksoddperiodic_lyaps_L20_100varyT}
\end{figure}
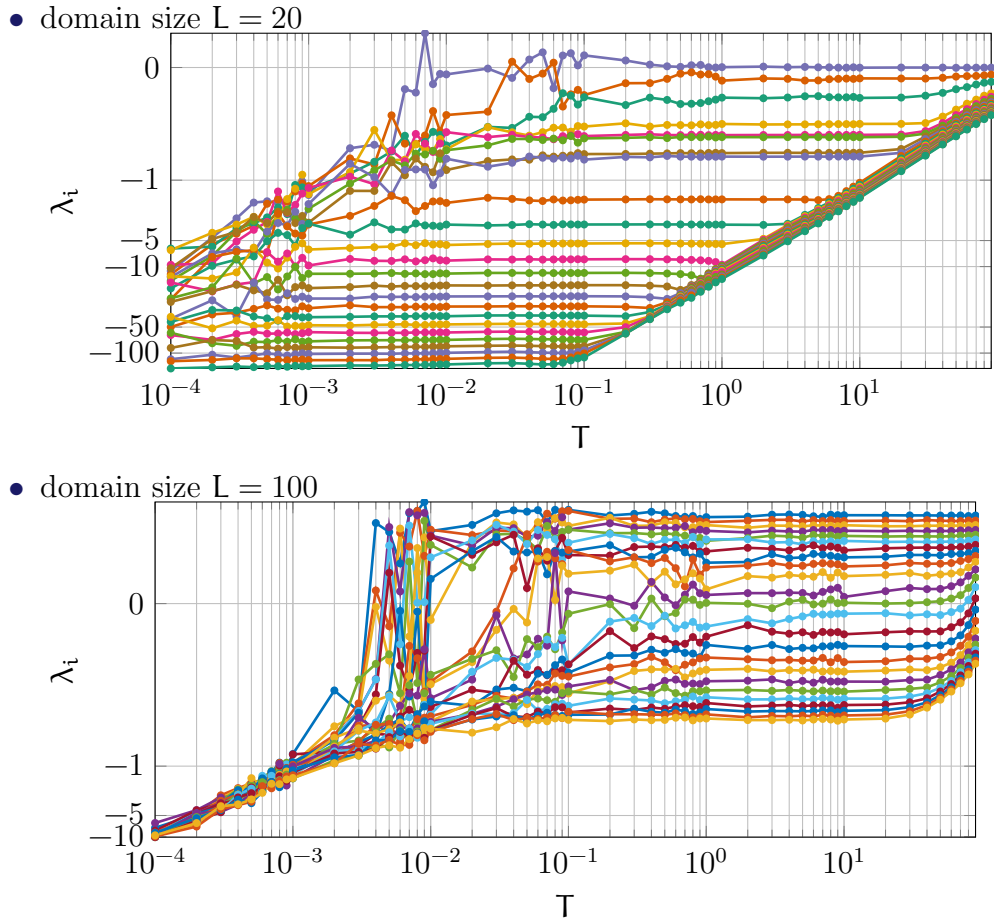

An important decision in the numerical calculation of the Lyapunov spectrum is the size of each time interval~$T$ between reorthonormalisations.
\cref{fig:ksoddperiodic_lyaps_L20_100varyT} demonstrates how the choice of interval~$T$ in \cref{alg:clyap-lyaps} affects the calculation of the $m=24$ most positive Lyapunov exponents, for small domain $L=20$ or large domain $L=100$, in the case of odd-periodic boundary conditions~\eqref{eq:opbc} (the periodic case~\eqref{eq:pbc} provides similar plots). 
For small interval~$T$ ($T\lesssim 0.1$ for $L=20$, and $T\lesssim 1$ for $L=100$) the most positive Lyapunov exponents are inaccurate because $T$~is too small to sufficiently capture the trajectory divergence, leading to an unstable \textsc{qr}~decomposition.
For larger interval~$T$ ($T\gtrsim 0.1$ for $L=20$, and $T\gtrsim 10$ for $L=100$) the most negative of the $m=24$ Lyapunov exponents evolve for too long and are corrupted towards the more positive exponents.
Based upon \cref{fig:ksoddperiodic_lyaps_L20_100varyT} we generally chose reorthonormalisation interval \(T=2\)\,.
With this choice of~$T$ we accurately resolve the most positive Lyapunov exponents, while also computing a sufficient number of negative Lyapunov exponents for the evaluation of the Kaplan--Yorke dimension for the chosen range of domains.

\begin{figure}
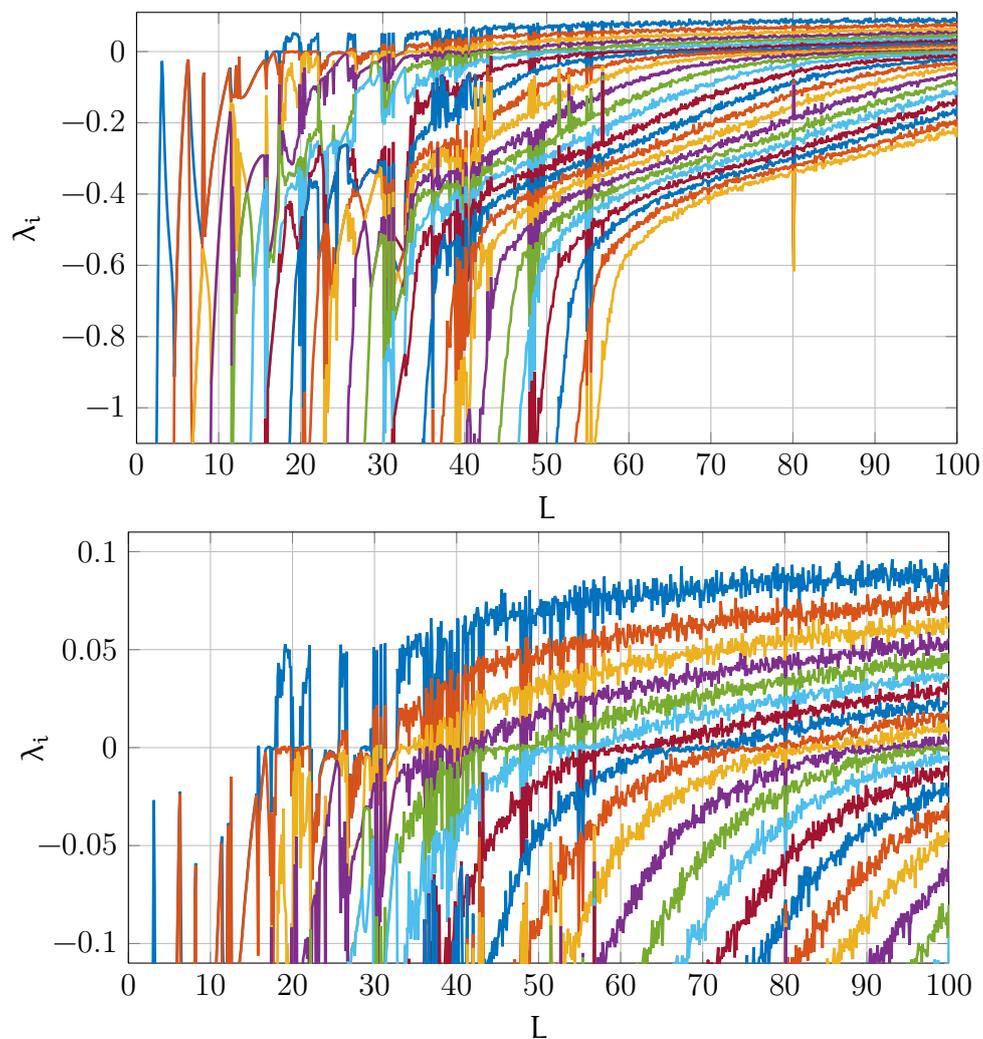

\centering
    \tikzsetfigurename{ksoddperiodic_lyaps_full}
    \input{ksoddperiodic_lyaps_full.tex}
    \tikzsetfigurename{ksoddperiodic_lyaps_zoom}
    \input{ksoddperiodic_lyaps_zoom.tex}
  \caption{The $24$~most positive Lyapunov exponents  $\lambda_1,\lambda_2,\ldots, \lambda_{24}$, computed for the  \kspde\ on the odd-periodic~\eqref{eq:opbc} spatial domain for domain sizes $0 < L \leq 100$\,.  The bottom plot provides a more detailed look at those Lyapunov exponents near zero. }
  \label{fig:cksodd-lyap}
\end{figure}

\Cref{fig:cks-lyap,fig:cksodd-lyap} plot the $24$~most positive Lyapunov exponents for the \kspde\ over different domain sizes, $0<L\leq 100$ with periodic~\eqref{eq:pbc} and odd-periodic~\eqref{eq:opbc} boundary conditions, respectively.
Although these calculations of the Lyapunov exponents contain noise, the exponents generally increase as $L$~increases, and larger values of~$L$ generally have more positive exponents.
However, the increase in the Lyapunov exponents is limited as they appear to be bounded above by about~$0.1$, for both the periodic and odd-periodic cases. This upper bound~$0.1$ matches with the upper bound observed for the rigid boundary condition case \citep{Yang09,Manneville85}.

\begin{figure}
\centering
    \tikzsetfigurename{ksperiodic_lyaps_full}
    \input{ksperiodic_lyaps_full.tex}
    \tikzsetfigurename{ksperiodic_lyaps_zoom}
    \input{ksperiodic_lyaps_zoom.tex}
  \caption{The $24$~most positive Lyapunov exponents  $\lambda_1,\lambda_2,\ldots, \lambda_{24}$, computed for the  \kspde\ on the periodic~\eqref{eq:pbc} spatial domain for domain sizes $0 < L \leq 100$\,. 
  The bottom plot zooms into those Lyapunov exponents near zero. }
  \label{fig:cks-lyap}
\end{figure}

We now further explore the Lyapunov exponents in the periodic case.
The structure of the positive Lyapunov exponents is noisy but there are some reasonably clear trends: here we show  that the \(i\)th~positive Lyapunov exponent on a domain length~\(L\) is approximately  
\begin{equation}
\lambda_i(L)\approx 0.093-0.94(i-0.39)/L\,.
\label{eq:lamil}
\end{equation}
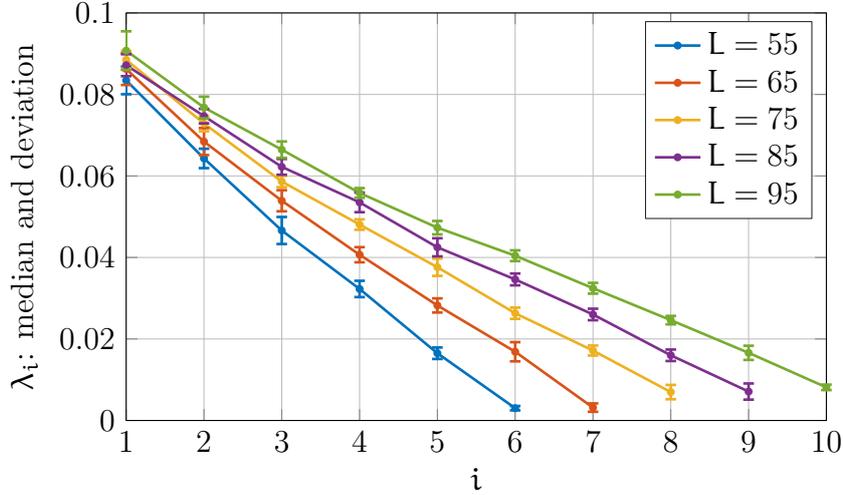
\begin{figure}
\centering
\input{formLEindex_2}
\caption{\label{fig:formLEindex}Positive Lyapunov exponents~\(\lambda_i\) for the periodic case~\eqref{eq:pbc}, \cref{fig:cks-lyap}, as a function of index~\(i\).
The joining lines are purely to aid visualisation.  
The data points plotted for each~\(L\) are the medians over the window~\([L-1,L+1]\), \(\pm\)-the mean absolute deviation.}
\end{figure}%
To derive this approximation, first look at the \(i\)-dependence for various fixed~\(L\): \cref{fig:formLEindex} plots Lyapunov exponents~\(\lambda_i(L)\) for fixed~\(L\) as a function of index~\(i\).
Since the data from \cref{fig:cks-lyap} is noisy, the median of the Lyapunov exponent is plotted where the median is over the window~\([L-1,L+1]\) (21~data points).
Then the vertical bars for each point represent plus-and-minus the mean absolute deviation (\textsc{mad}): these statistics are more robust to outliers than the usual mean and standard deviation, and so appear more suitable here.
\cref{fig:formLEindex} indicates the Lyapunov exponents are, for fixed~\(L\), approximately equi-spaced in~\(i\), especially for \(i\leq5\).
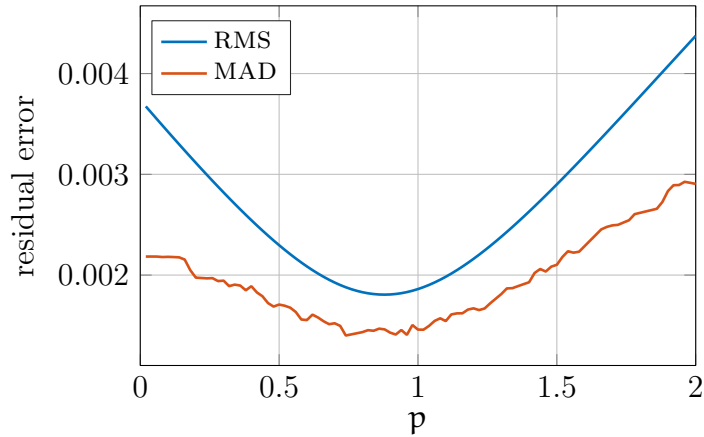
\begin{figure}
\centering
\input{formLEerror_2}
\caption{\label{fig:formLEerror}Fit the Lyapunov exponents \(\lambda_i(L)\approx a+(b+ci)/L^p\) for various exponents~\(p\) and then plot the residual error as a function of exponent~\(p\): here the root-mean-square error (\textsc{rms}), and the mean absolute deviation (\textsc{mad}).  The minimum suggests the optimum exponent \(p\approx1\).}
\end{figure}%
The magnitude of the slope of the \(i\)-dependence decreases as the domain length~\(L\) increases, so we try to fit a function of the power-law form \(\lambda_i(L)\approx a+(b+ci)/L^p\) for various exponents~\(p\).
\cref{fig:formLEerror} plots the residual error in the fit as a function of exponent~\(p\) showing that there is a minimum error at \(p\approx1\): this minimum occurs both in the \textsc{rms} error and the \textsc{mad} error.
In view of the fluctuations in the \textsc{mad}, it seems reasonable to choose the case of the exponent \(p=1\) reported by equation~\eqref{eq:lamil}.
Moreover, this is the exponent which best fits our preconception that the chaos in the Kuramoto--Sivashinsky \pde\ is `extensive'---that the number of positive Lyapunov exponents increases linearly with domain length~\(L\).
The next \cref{sec:ckyd} explores this issue further via the Kaplan--Yorke dimension, and finds results consistent with the approximate formula~\eqref{eq:lamil}.

\section{Compute the Kaplan--Yorke dimension}
\label{sec:ckyd}

The Kaplan--Yorke dimension is a measure of the dimension of an attractor~\citep{Kaplan79} and is defined in terms of a sum of the most positive Lyapunov exponents,
\begin{equation}
  D_{\text{KY}} = j + \frac{\sum_{i=1}^j \lambda_i}{\left| \lambda_{j+1} \right|},
  \label{eqn:dky}
\end{equation}
where~$j$ is the largest index such that $\sum^j_{i=1} \lambda_i \geq 0$\,.
The Kaplan--Yorke dimension is an upper bound of the Hausdorff dimension of the attractor, and as each of the $j$~Lyapunov exponents correspond to an orthogonal direction, the Kaplan--Yorke dimension approximates the minimum number of modes required to describe the emergent dynamics of the system on the attractor~\citep{Grassberger83}.

\begin{figure}
  \centering
\begin{itemize}
\item   periodic boundary conditions \cref{eq:pbc}\\
\tikzsetfigurename{ksperiodic_dky}
  \input{ksperiodic_dky.tex}
\item odd-periodic~\cref{eq:opbc}\\
  \tikzsetfigurename{ksoddperiodic_dky}
  \input{ksoddperiodic_dky.tex}
\end{itemize}
  \caption{The Kaplan--Yorke dimension~$D_{\text{KY}}$ dependence on domain size~$L$, computed for  the \kspde\ on the (top) periodic spatial domain~\eqref{eq:pbc}, and (bottom) odd-periodic spatial domain~\eqref{eq:opbc}.
In both cases, the Kaplan--Yorke dimension increases linearly with the domain size~$L$, with $D_{\text{KY}} \approx 0.226L-c$ for \(c=0.160\) and~\(2.106\), respectively.
This linearity appears most accurate for larger domains, $L\gtrsim 80$\,.
At lower domain sizes, $L\lesssim 80$, there are several non-chaotic regions where the Kaplan--Yorke dimension is not of interest.} 
  \label{fig:ckspde-dky}
\end{figure}
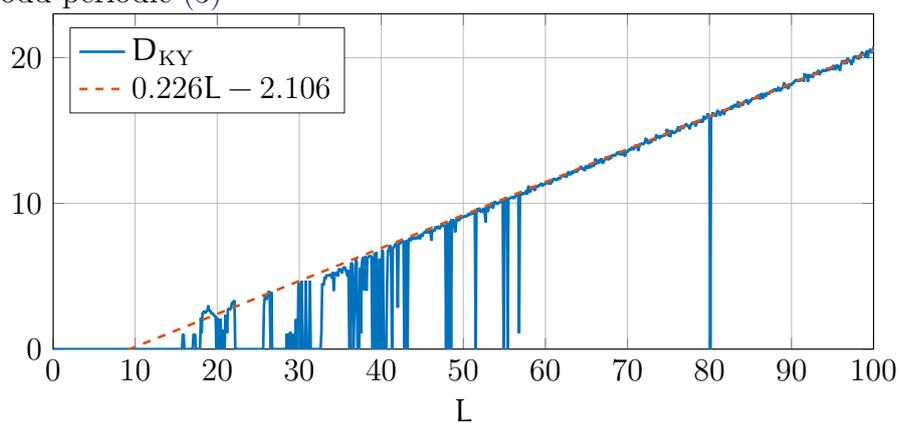

In the formula~\eqref{eqn:dky} the term \({\sum_{i=1}^j \lambda_i}\big/{\left| \lambda_{j+1} \right|}\) is usually a fraction in~\((0,1)\) and so the index~\(j\) is roughly the Kaplan--Yorke dimension.
Using the approximation~\eqref{eq:lamil} to the Lyapunov exponents for the periodic case, one may straightforwardly estimate the~\(j\) for which \(\sum_{i=1}^j \lambda_i\approx 0\), namely \(j\approx 0.2L-0.2\)\,.
This is acceptably close to the Kaplan--Yorke dimension \(D_{\text{KY}}\approx0.226\,L-0.160\) shown in \cref{fig:ckspde-dky} and obtained from extensive computation.

\cref{tab:cks-lyap} presents Lyapunov exponents and Kaplan--Yorke dimensions for six different periodic domain sizes~$L=12,13.5,22,36,60,100$\,,%
\footnote{The case \(L=22\) is chosen for comparison with the Lyapunov exponents of \cite{Cvitanovic10}: our Lyapunov exponents agree with theirs to a difference of about~\(0.002\).}
whereas \cref{tab:cksodd-lyap} presents Lyapunov exponents and Kaplan--Yorke dimensions for six different odd-periodic domain sizes~$L=17.5,18.1,18.2,41,60,100$\,.
These tables demonstrate how the increasingly positive Lyapunov exponents reveal the onset of chaotic dynamics and the increasing dimension of the chaotic attractor. 
These calculations of the Lyapunov exponents and Kaplan--Yorke dimensions with odd periodic domains are compatible with
previously calculated Kaplan--Yorke dimensions~\citep{Rempel04}.

\begin{table}
  \centering
  \caption{the $12$~most positive Lyapunov exponents and the Kaplan--Yorke dimension for several domain sizes~$L$ and periodic boundary conditions~\eqref{eq:pbc}.}  
\label{tab:cks-lyap}
$\begin{array}{r | r r r r r r}
& L=12 & L=13.5 & L=22 & L=36 & L=60 & L=100 \\
\hline
\lambda_1 & 0.003 & 0.059 & 0.043 & 0.080 & 0.089 & 0.088 \\
\lambda_2 & -0.005 & 0.004 & 0.003 & 0.056 & 0.067 & 0.082 \\
\lambda_3 & -0.088 & -0.004 & 0.002 & 0.014 & 0.055 & 0.070 \\
\lambda_4 & -0.089 & -0.227 & -0.004 & 0.003 & 0.041 & 0.061 \\
\lambda_5 & -0.186 & -0.730 & -0.008 & -0.003 & 0.030 & 0.048 \\
\lambda_6 & -3.524 & -1.467 & -0.185 & -0.004 & 0.005 & 0.041 \\
\lambda_7 & -3.525 & -1.529 & -0.253 & -0.021 & 0.003 & 0.033 \\
\lambda_8 & -9.835 & -6.956 & -0.296 & -0.088 & 0.000 & 0.028 \\
\lambda_9 & -9.849 & -6.963 & -0.309 & -0.160 & -0.004 & 0.018 \\
\lambda_{10} & -9.959 & -7.977 & -1.965 & -0.224 & -0.009 & 0.012 \\
\lambda_{11} & -10.01 & -7.993 & -1.967 & -0.309 & -0.029 & 0.005 \\
\lambda_{12} & -10.12 & -9.199 & -5.599 & -0.373 & -0.066 & 0.003 \\[8pt]
D_{\text{KY}} & 1.663 & 3.259 & 5.198 & 8.229 & 13.56 & 22.44
\end{array}$
\end{table}

\begin{table}
  \centering
  \caption{the $12$~most positive Lyapunov exponents and the Kaplan--Yorke dimension for several domain sizes~$L$ and odd-periodic boundary conditions~\eqref{eq:opbc}.}
  \label{tab:cksodd-lyap}
$\begin{array}{r | r r r r r r}
& L=17.5 & L=18.1 & L=18.2 & L=41 & L=60 & L=100 \\
\hline
\lambda_1 & -0.001 & 0.000 & 0.036 & 0.067 & 0.076 & 0.094 \\
\lambda_2 & -0.166 & -0.003 & -0.001 & 0.038 & 0.056 & 0.077 \\
\lambda_3 & -0.272 & -0.194 & -0.073 & 0.017 & 0.042 & 0.063 \\
\lambda_4 & -0.299 & -0.280 & -0.268 & 0.001 & 0.027 & 0.056 \\
\lambda_5 & -0.300 & -0.377 & -0.359 & -0.008 & 0.021 & 0.044 \\
\lambda_6 & -0.526 & -4.813 & -4.044 & -0.029 & 0.006 & 0.036 \\
\lambda_7 & -0.619 & -4.923 & -4.348 & -0.076 & 0.000 & 0.031 \\
\lambda_8 & -1.794 & -1.391 & -1.395 & -0.162 & -0.007 & 0.022 \\
\lambda_9 & -3.780 & -3.145 & -3.070 & -0.237 & -0.029 & 0.017 \\
\lambda_{10} & -6.513 & -5.525 & -5.383 & -0.283 & -0.050 & 0.008 \\
\lambda_{11} & -9.692 & -8.700 & -8.538 & -0.318 & -0.094 & 0.001 \\
\lambda_{12} & -9.854 & -9.540 & -10.10 & -0.355 & -0.146 & 0.000 \\[8pt]
D_{\text{KY}} & 0.000 & 1.081 & 2.482 & 7.056 & 11.35 & 20.75
\end{array}$
\end{table}

\cref{fig:ckspde-dky} shows that as the domain size~$L$ increases, the Kaplan--Yorke dimension scales linearly with~$D_{\text{KY}} \propto 0.226L$ for sufficiently large~$L$ ($L\gtrsim 80$). 
Similarly, in the case of rigid boundary conditions $u,\partial_x u=0$ at $x=0$ and at $x=L$ both \citet{Manneville85} and \citet{Tajima02} found the Kaplan--Yorke dimension to scale as~$0.230L$ when $50<L<400$\,.  
The small~\(2\%\) difference in the coefficient suggests that the  scaling of the Kaplan--Yorke dimension for the \kspde\ on domain sizes $L\gtrsim 80$ only depends on the nature of the given boundary condition through an additive constant. 
In contrast, for smaller domain sizes $L\lesssim 50$\,, all points in the domain are somewhat close to a boundary and boundary effects play a more dominant role in the dynamics. 
The linear scaling of the attractor, here measured with the Kaplan--Yorke dimension, is a defining feature of an extensively chaotic system~\citep[e.g.,][]{Cross93,Greenside96}.

\section{Conclusion}
\label{sec:cconc}

Through an exhaustive computation and analysis of the positive and least negative Lyapunov exponents, we investigated the development of spatio-temporal chaos in the Kuramoto--Sivashinsky \textsc{pde}~\cref{eqn:kspde} as the domain size~$L$ increases.
We found new details about how the Lyapunov exponents and the Kaplan--Yorke dimension increase with domain size, and are able to identify successive transitions into more chaotic regimes as individual Lyapunov exponents change sign from negative to positive, indicating additional directions in which trajectories of the chaotic system diverge. 

The spatial extensivity of the \kspde\ that we have confirmed here in new detail indicates that the system in a large domain may be viewed as composed of interacting subsystems, approximately uncorrelated for short enough times \cite[e.g.,][]{Wittenberg99, Yang09, Greenside96}.
This interpretation suggests that we should be able to successfully simulate the `turbulence' in the \kspde\ on very large domains by appropriately coupling relatively small patches of simulations across space using the equation-free paradigm \cite[e.g.,][]{Kevrekidis09a}.
Exactly what may be appropriate coupling is the subject of ongoing research.

\paragraph{Acknowledgements}
This research was partly supported by the research grant DP150102385 from the Australian Research Council, and RAE was financially supported by the Australian Government Research Training Program. 
\end{document}

%% file: ksperiodic_plots.tex
\definecolor{mycolor1}{rgb}{0.09640,0.75000,0.71204}%
\definecolor{mycolor2}{rgb}{0.22064,0.46026,0.99729}%
\definecolor{mycolor3}{rgb}{0.19633,0.48472,0.98915}%
\definecolor{mycolor4}{rgb}{0.18340,0.50737,0.97980}%
\definecolor{mycolor5}{rgb}{0.17644,0.54990,0.95202}%
\definecolor{mycolor6}{rgb}{0.16874,0.57026,0.93587}%
\definecolor{mycolor7}{rgb}{0.14603,0.60912,0.90786}%
\definecolor{mycolor8}{rgb}{0.15400,0.59020,0.92180}%
\definecolor{mycolor9}{rgb}{0.13802,0.62763,0.89729}%
\definecolor{mycolor10}{rgb}{0.12481,0.64593,0.88834}%
\definecolor{mycolor11}{rgb}{0.11125,0.66350,0.87631}%
\definecolor{mycolor12}{rgb}{0.09521,0.67983,0.85978}%
\definecolor{mycolor13}{rgb}{0.02967,0.70817,0.81633}%
\definecolor{mycolor14}{rgb}{0.00357,0.72027,0.79170}%
\definecolor{mycolor15}{rgb}{0.00666,0.73121,0.76601}%
\definecolor{mycolor16}{rgb}{0.14077,0.75840,0.68416}%
\definecolor{mycolor17}{rgb}{0.17170,0.76696,0.65544}%
\definecolor{mycolor18}{rgb}{0.19377,0.77577,0.62510}%
\definecolor{mycolor19}{rgb}{0.24696,0.79180,0.55674}%
\definecolor{mycolor20}{rgb}{0.29061,0.79729,0.51883}%
\definecolor{mycolor21}{rgb}{0.34064,0.80080,0.47886}%
\definecolor{mycolor22}{rgb}{0.39090,0.80287,0.43545}%
\definecolor{mycolor23}{rgb}{0.44563,0.80242,0.39092}%
\definecolor{mycolor24}{rgb}{0.50440,0.79930,0.34800}%
\definecolor{mycolor25}{rgb}{0.56156,0.79423,0.30448}%
\definecolor{mycolor26}{rgb}{0.61740,0.78762,0.26124}%
\definecolor{mycolor27}{rgb}{0.72420,0.76984,0.19103}%
\definecolor{mycolor28}{rgb}{0.77383,0.75980,0.16461}%
\definecolor{mycolor29}{rgb}{0.82031,0.74981,0.15353}%
\definecolor{mycolor30}{rgb}{0.86343,0.74060,0.15963}%
\begin{tikzpicture}

\begin{axis}[%
width=1.6in,
height=2.1in,
at={(0in,0in)},
scale only axis,
point meta min=-3.5,
point meta max=3.5,
xmin=0,
xmax=12,
xtick={0,6,12},
xlabel={$x$},
ymin=2000,
ymax=2200,
ytick={2000,2100,2200},
ylabel={$t$},
axis background/.style={fill=white},
colormap={mymap}{[1pt] rgb(0pt)=(0.2422,0.1504,0.6603); rgb(1pt)=(0.25039,0.164995,0.707614); rgb(2pt)=(0.257771,0.181781,0.751138); rgb(3pt)=(0.264729,0.197757,0.795214); rgb(4pt)=(0.270648,0.214676,0.836371); rgb(5pt)=(0.275114,0.234238,0.870986); rgb(6pt)=(0.2783,0.255871,0.899071); rgb(7pt)=(0.280333,0.278233,0.9221); rgb(8pt)=(0.281338,0.300595,0.941376); rgb(9pt)=(0.281014,0.322757,0.957886); rgb(10pt)=(0.279467,0.344671,0.971676); rgb(11pt)=(0.275971,0.366681,0.982905); rgb(12pt)=(0.269914,0.3892,0.9906); rgb(13pt)=(0.260243,0.412329,0.995157); rgb(14pt)=(0.244033,0.435833,0.998833); rgb(15pt)=(0.220643,0.460257,0.997286); rgb(16pt)=(0.196333,0.484719,0.989152); rgb(17pt)=(0.183405,0.507371,0.979795); rgb(18pt)=(0.178643,0.528857,0.968157); rgb(19pt)=(0.176438,0.549905,0.952019); rgb(20pt)=(0.168743,0.570262,0.935871); rgb(21pt)=(0.154,0.5902,0.9218); rgb(22pt)=(0.146029,0.609119,0.907857); rgb(23pt)=(0.138024,0.627629,0.89729); rgb(24pt)=(0.124814,0.645929,0.888343); rgb(25pt)=(0.111252,0.6635,0.876314); rgb(26pt)=(0.0952095,0.679829,0.859781); rgb(27pt)=(0.0688714,0.694771,0.839357); rgb(28pt)=(0.0296667,0.708167,0.816333); rgb(29pt)=(0.00357143,0.720267,0.7917); rgb(30pt)=(0.00665714,0.731214,0.766014); rgb(31pt)=(0.0433286,0.741095,0.73941); rgb(32pt)=(0.0963952,0.75,0.712038); rgb(33pt)=(0.140771,0.7584,0.684157); rgb(34pt)=(0.1717,0.766962,0.655443); rgb(35pt)=(0.193767,0.775767,0.6251); rgb(36pt)=(0.216086,0.7843,0.5923); rgb(37pt)=(0.246957,0.791795,0.556743); rgb(38pt)=(0.290614,0.79729,0.518829); rgb(39pt)=(0.340643,0.8008,0.478857); rgb(40pt)=(0.3909,0.802871,0.435448); rgb(41pt)=(0.445629,0.802419,0.390919); rgb(42pt)=(0.5044,0.7993,0.348); rgb(43pt)=(0.561562,0.794233,0.304481); rgb(44pt)=(0.617395,0.787619,0.261238); rgb(45pt)=(0.671986,0.779271,0.2227); rgb(46pt)=(0.7242,0.769843,0.191029); rgb(47pt)=(0.773833,0.759805,0.16461); rgb(48pt)=(0.820314,0.749814,0.153529); rgb(49pt)=(0.863433,0.7406,0.159633); rgb(50pt)=(0.903543,0.733029,0.177414); rgb(51pt)=(0.939257,0.728786,0.209957); rgb(52pt)=(0.972757,0.729771,0.239443); rgb(53pt)=(0.995648,0.743371,0.237148); rgb(54pt)=(0.996986,0.765857,0.219943); rgb(55pt)=(0.995205,0.789252,0.202762); rgb(56pt)=(0.9892,0.813567,0.188533); rgb(57pt)=(0.978629,0.838629,0.176557); rgb(58pt)=(0.967648,0.8639,0.16429); rgb(59pt)=(0.96101,0.889019,0.153676); rgb(60pt)=(0.959671,0.913457,0.142257); rgb(61pt)=(0.962795,0.937338,0.12651); rgb(62pt)=(0.969114,0.960629,0.106362); rgb(63pt)=(0.9769,0.9839,0.0805)},
colorbar,
colorbar style={ylabel={$u(x,t)$}, y label style={at={(2.10,0.5)}}, width={0.15in}, ytick={-3,-2, -1,0,1,2, 3}, at={(2.75,1.0)},
height=2*\pgfkeysvalueof{/pgfplots/parent axis height}+0.9in},
axis on top
]
\addplot[thick] graphics[xmin=0,ymin=2000,xmax=12,ymax=2200] {ksperiodic_L12plot.jpg};
\end{axis}

\begin{axis}[%
width=1.6in,
height=2.1in,
at={(2.6in,0in)},
scale only axis,
point meta min=-3.5,
point meta max=3.5,
xmin=0,
xmax=13.5,
xtick={0,6.75,13.5},
xlabel={$x$},
ymin=2000,
ymax=2200,
ytick={2000,2100,2200},
ylabel={$t$},
axis background/.style={fill=white},
colormap={mymap}{[1pt] rgb(0pt)=(0.2422,0.1504,0.6603); rgb(1pt)=(0.25039,0.164995,0.707614); rgb(2pt)=(0.257771,0.181781,0.751138); rgb(3pt)=(0.264729,0.197757,0.795214); rgb(4pt)=(0.270648,0.214676,0.836371); rgb(5pt)=(0.275114,0.234238,0.870986); rgb(6pt)=(0.2783,0.255871,0.899071); rgb(7pt)=(0.280333,0.278233,0.9221); rgb(8pt)=(0.281338,0.300595,0.941376); rgb(9pt)=(0.281014,0.322757,0.957886); rgb(10pt)=(0.279467,0.344671,0.971676); rgb(11pt)=(0.275971,0.366681,0.982905); rgb(12pt)=(0.269914,0.3892,0.9906); rgb(13pt)=(0.260243,0.412329,0.995157); rgb(14pt)=(0.244033,0.435833,0.998833); rgb(15pt)=(0.220643,0.460257,0.997286); rgb(16pt)=(0.196333,0.484719,0.989152); rgb(17pt)=(0.183405,0.507371,0.979795); rgb(18pt)=(0.178643,0.528857,0.968157); rgb(19pt)=(0.176438,0.549905,0.952019); rgb(20pt)=(0.168743,0.570262,0.935871); rgb(21pt)=(0.154,0.5902,0.9218); rgb(22pt)=(0.146029,0.609119,0.907857); rgb(23pt)=(0.138024,0.627629,0.89729); rgb(24pt)=(0.124814,0.645929,0.888343); rgb(25pt)=(0.111252,0.6635,0.876314); rgb(26pt)=(0.0952095,0.679829,0.859781); rgb(27pt)=(0.0688714,0.694771,0.839357); rgb(28pt)=(0.0296667,0.708167,0.816333); rgb(29pt)=(0.00357143,0.720267,0.7917); rgb(30pt)=(0.00665714,0.731214,0.766014); rgb(31pt)=(0.0433286,0.741095,0.73941); rgb(32pt)=(0.0963952,0.75,0.712038); rgb(33pt)=(0.140771,0.7584,0.684157); rgb(34pt)=(0.1717,0.766962,0.655443); rgb(35pt)=(0.193767,0.775767,0.6251); rgb(36pt)=(0.216086,0.7843,0.5923); rgb(37pt)=(0.246957,0.791795,0.556743); rgb(38pt)=(0.290614,0.79729,0.518829); rgb(39pt)=(0.340643,0.8008,0.478857); rgb(40pt)=(0.3909,0.802871,0.435448); rgb(41pt)=(0.445629,0.802419,0.390919); rgb(42pt)=(0.5044,0.7993,0.348); rgb(43pt)=(0.561562,0.794233,0.304481); rgb(44pt)=(0.617395,0.787619,0.261238); rgb(45pt)=(0.671986,0.779271,0.2227); rgb(46pt)=(0.7242,0.769843,0.191029); rgb(47pt)=(0.773833,0.759805,0.16461); rgb(48pt)=(0.820314,0.749814,0.153529); rgb(49pt)=(0.863433,0.7406,0.159633); rgb(50pt)=(0.903543,0.733029,0.177414); rgb(51pt)=(0.939257,0.728786,0.209957); rgb(52pt)=(0.972757,0.729771,0.239443); rgb(53pt)=(0.995648,0.743371,0.237148); rgb(54pt)=(0.996986,0.765857,0.219943); rgb(55pt)=(0.995205,0.789252,0.202762); rgb(56pt)=(0.9892,0.813567,0.188533); rgb(57pt)=(0.978629,0.838629,0.176557); rgb(58pt)=(0.967648,0.8639,0.16429); rgb(59pt)=(0.96101,0.889019,0.153676); rgb(60pt)=(0.959671,0.913457,0.142257); rgb(61pt)=(0.962795,0.937338,0.12651); rgb(62pt)=(0.969114,0.960629,0.106362); rgb(63pt)=(0.9769,0.9839,0.0805)},
axis on top
]
\addplot[thick] graphics[xmin=0,ymin=2000,xmax=13.5,ymax=2200] {ksperiodic_L135plot.jpg};
\end{axis}

\begin{axis}[%
width=1.6in,
height=2.1in,
at={(0in,-3in)},
scale only axis,
point meta min=-3.5,
point meta max=3.5,
xmin=0,
xmax=36,
xtick={0,18,36},
xlabel={$x$},
ymin=2000,
ymax=2200,
ytick={2000,2100,2200},
ylabel={$t$},
axis background/.style={fill=white},
colormap={mymap}{[1pt] rgb(0pt)=(0.2422,0.1504,0.6603); rgb(1pt)=(0.25039,0.164995,0.707614); rgb(2pt)=(0.257771,0.181781,0.751138); rgb(3pt)=(0.264729,0.197757,0.795214); rgb(4pt)=(0.270648,0.214676,0.836371); rgb(5pt)=(0.275114,0.234238,0.870986); rgb(6pt)=(0.2783,0.255871,0.899071); rgb(7pt)=(0.280333,0.278233,0.9221); rgb(8pt)=(0.281338,0.300595,0.941376); rgb(9pt)=(0.281014,0.322757,0.957886); rgb(10pt)=(0.279467,0.344671,0.971676); rgb(11pt)=(0.275971,0.366681,0.982905); rgb(12pt)=(0.269914,0.3892,0.9906); rgb(13pt)=(0.260243,0.412329,0.995157); rgb(14pt)=(0.244033,0.435833,0.998833); rgb(15pt)=(0.220643,0.460257,0.997286); rgb(16pt)=(0.196333,0.484719,0.989152); rgb(17pt)=(0.183405,0.507371,0.979795); rgb(18pt)=(0.178643,0.528857,0.968157); rgb(19pt)=(0.176438,0.549905,0.952019); rgb(20pt)=(0.168743,0.570262,0.935871); rgb(21pt)=(0.154,0.5902,0.9218); rgb(22pt)=(0.146029,0.609119,0.907857); rgb(23pt)=(0.138024,0.627629,0.89729); rgb(24pt)=(0.124814,0.645929,0.888343); rgb(25pt)=(0.111252,0.6635,0.876314); rgb(26pt)=(0.0952095,0.679829,0.859781); rgb(27pt)=(0.0688714,0.694771,0.839357); rgb(28pt)=(0.0296667,0.708167,0.816333); rgb(29pt)=(0.00357143,0.720267,0.7917); rgb(30pt)=(0.00665714,0.731214,0.766014); rgb(31pt)=(0.0433286,0.741095,0.73941); rgb(32pt)=(0.0963952,0.75,0.712038); rgb(33pt)=(0.140771,0.7584,0.684157); rgb(34pt)=(0.1717,0.766962,0.655443); rgb(35pt)=(0.193767,0.775767,0.6251); rgb(36pt)=(0.216086,0.7843,0.5923); rgb(37pt)=(0.246957,0.791795,0.556743); rgb(38pt)=(0.290614,0.79729,0.518829); rgb(39pt)=(0.340643,0.8008,0.478857); rgb(40pt)=(0.3909,0.802871,0.435448); rgb(41pt)=(0.445629,0.802419,0.390919); rgb(42pt)=(0.5044,0.7993,0.348); rgb(43pt)=(0.561562,0.794233,0.304481); rgb(44pt)=(0.617395,0.787619,0.261238); rgb(45pt)=(0.671986,0.779271,0.2227); rgb(46pt)=(0.7242,0.769843,0.191029); rgb(47pt)=(0.773833,0.759805,0.16461); rgb(48pt)=(0.820314,0.749814,0.153529); rgb(49pt)=(0.863433,0.7406,0.159633); rgb(50pt)=(0.903543,0.733029,0.177414); rgb(51pt)=(0.939257,0.728786,0.209957); rgb(52pt)=(0.972757,0.729771,0.239443); rgb(53pt)=(0.995648,0.743371,0.237148); rgb(54pt)=(0.996986,0.765857,0.219943); rgb(55pt)=(0.995205,0.789252,0.202762); rgb(56pt)=(0.9892,0.813567,0.188533); rgb(57pt)=(0.978629,0.838629,0.176557); rgb(58pt)=(0.967648,0.8639,0.16429); rgb(59pt)=(0.96101,0.889019,0.153676); rgb(60pt)=(0.959671,0.913457,0.142257); rgb(61pt)=(0.962795,0.937338,0.12651); rgb(62pt)=(0.969114,0.960629,0.106362); rgb(63pt)=(0.9769,0.9839,0.0805)},
axis on top
]
\addplot[thick] graphics[xmin=0,ymin=2000,xmax=36,ymax=2200] {ksperiodic_L36plot.jpg};
\end{axis}

\begin{axis}[%
width=1.6in,
height=2.1in,
at={(2.6in,-3in)},
scale only axis,
point meta min=-3.5,
point meta max=3.5,
xmin=0,
xmax=100,
xtick={0,50,100},
xlabel={$x$},
ymin=2000,
ymax=2200,
ytick={2000,2100,2200},
ylabel={$t$},
axis background/.style={fill=white},
colormap={mymap}{[1pt] rgb(0pt)=(0.2422,0.1504,0.6603); rgb(1pt)=(0.25039,0.164995,0.707614); rgb(2pt)=(0.257771,0.181781,0.751138); rgb(3pt)=(0.264729,0.197757,0.795214); rgb(4pt)=(0.270648,0.214676,0.836371); rgb(5pt)=(0.275114,0.234238,0.870986); rgb(6pt)=(0.2783,0.255871,0.899071); rgb(7pt)=(0.280333,0.278233,0.9221); rgb(8pt)=(0.281338,0.300595,0.941376); rgb(9pt)=(0.281014,0.322757,0.957886); rgb(10pt)=(0.279467,0.344671,0.971676); rgb(11pt)=(0.275971,0.366681,0.982905); rgb(12pt)=(0.269914,0.3892,0.9906); rgb(13pt)=(0.260243,0.412329,0.995157); rgb(14pt)=(0.244033,0.435833,0.998833); rgb(15pt)=(0.220643,0.460257,0.997286); rgb(16pt)=(0.196333,0.484719,0.989152); rgb(17pt)=(0.183405,0.507371,0.979795); rgb(18pt)=(0.178643,0.528857,0.968157); rgb(19pt)=(0.176438,0.549905,0.952019); rgb(20pt)=(0.168743,0.570262,0.935871); rgb(21pt)=(0.154,0.5902,0.9218); rgb(22pt)=(0.146029,0.609119,0.907857); rgb(23pt)=(0.138024,0.627629,0.89729); rgb(24pt)=(0.124814,0.645929,0.888343); rgb(25pt)=(0.111252,0.6635,0.876314); rgb(26pt)=(0.0952095,0.679829,0.859781); rgb(27pt)=(0.0688714,0.694771,0.839357); rgb(28pt)=(0.0296667,0.708167,0.816333); rgb(29pt)=(0.00357143,0.720267,0.7917); rgb(30pt)=(0.00665714,0.731214,0.766014); rgb(31pt)=(0.0433286,0.741095,0.73941); rgb(32pt)=(0.0963952,0.75,0.712038); rgb(33pt)=(0.140771,0.7584,0.684157); rgb(34pt)=(0.1717,0.766962,0.655443); rgb(35pt)=(0.193767,0.775767,0.6251); rgb(36pt)=(0.216086,0.7843,0.5923); rgb(37pt)=(0.246957,0.791795,0.556743); rgb(38pt)=(0.290614,0.79729,0.518829); rgb(39pt)=(0.340643,0.8008,0.478857); rgb(40pt)=(0.3909,0.802871,0.435448); rgb(41pt)=(0.445629,0.802419,0.390919); rgb(42pt)=(0.5044,0.7993,0.348); rgb(43pt)=(0.561562,0.794233,0.304481); rgb(44pt)=(0.617395,0.787619,0.261238); rgb(45pt)=(0.671986,0.779271,0.2227); rgb(46pt)=(0.7242,0.769843,0.191029); rgb(47pt)=(0.773833,0.759805,0.16461); rgb(48pt)=(0.820314,0.749814,0.153529); rgb(49pt)=(0.863433,0.7406,0.159633); rgb(50pt)=(0.903543,0.733029,0.177414); rgb(51pt)=(0.939257,0.728786,0.209957); rgb(52pt)=(0.972757,0.729771,0.239443); rgb(53pt)=(0.995648,0.743371,0.237148); rgb(54pt)=(0.996986,0.765857,0.219943); rgb(55pt)=(0.995205,0.789252,0.202762); rgb(56pt)=(0.9892,0.813567,0.188533); rgb(57pt)=(0.978629,0.838629,0.176557); rgb(58pt)=(0.967648,0.8639,0.16429); rgb(59pt)=(0.96101,0.889019,0.153676); rgb(60pt)=(0.959671,0.913457,0.142257); rgb(61pt)=(0.962795,0.937338,0.12651); rgb(62pt)=(0.969114,0.960629,0.106362); rgb(63pt)=(0.9769,0.9839,0.0805)},
axis on top
]
\addplot[thick] graphics[xmin=0,ymin=2000,xmax=100,ymax=2200] {ksperiodic_L100plot.jpg};
\end{axis}
\end{tikzpicture}%

%% file: ksoddperiodic_plots.tex
\definecolor{mycolor1}{rgb}{0.09640,0.75000,0.71204}%
\definecolor{mycolor2}{rgb}{0.06887,0.69477,0.83936}%
\definecolor{mycolor3}{rgb}{0.09521,0.67983,0.85978}%
\definecolor{mycolor4}{rgb}{0.12481,0.64593,0.88834}%
\definecolor{mycolor5}{rgb}{0.14603,0.60912,0.90786}%
\definecolor{mycolor6}{rgb}{0.15400,0.59020,0.92180}%
\definecolor{mycolor7}{rgb}{0.17644,0.54990,0.95202}%
\definecolor{mycolor8}{rgb}{0.18340,0.50737,0.97980}%
\definecolor{mycolor9}{rgb}{0.19633,0.48472,0.98915}%
\definecolor{mycolor10}{rgb}{0.24403,0.43583,0.99883}%
\definecolor{mycolor11}{rgb}{0.26991,0.38920,0.99060}%
\definecolor{mycolor12}{rgb}{0.27597,0.36668,0.98290}%
\definecolor{mycolor13}{rgb}{0.28101,0.32276,0.95789}%
\definecolor{mycolor14}{rgb}{0.28033,0.27823,0.92210}%
\definecolor{mycolor15}{rgb}{0.27830,0.25587,0.89907}%
\definecolor{mycolor16}{rgb}{0.00357,0.72027,0.79170}%
\definecolor{mycolor17}{rgb}{0.98920,0.81357,0.18853}%
\definecolor{mycolor18}{rgb}{0.93926,0.72879,0.20996}%
\definecolor{mycolor19}{rgb}{0.97276,0.72977,0.23944}%
\definecolor{mycolor20}{rgb}{0.99699,0.76586,0.21994}%
\definecolor{mycolor21}{rgb}{0.17170,0.76696,0.65544}%
\definecolor{mycolor22}{rgb}{0.86343,0.74060,0.15963}%
\definecolor{mycolor23}{rgb}{0.21609,0.78430,0.59230}%
\definecolor{mycolor24}{rgb}{0.77383,0.75980,0.16461}%
\definecolor{mycolor25}{rgb}{0.24696,0.79180,0.55674}%
\definecolor{mycolor26}{rgb}{0.34064,0.80080,0.47886}%
\definecolor{mycolor27}{rgb}{0.44563,0.80242,0.39092}%
\definecolor{mycolor28}{rgb}{0.50440,0.79930,0.34800}%
\definecolor{mycolor29}{rgb}{0.61740,0.78762,0.26124}%
\definecolor{mycolor30}{rgb}{0.72420,0.76984,0.19103}%
\begin{tikzpicture}

\begin{axis}[%
width=1.6in,
height=2.1in,
at={(0in,0in)},
scale only axis,
point meta min=-3.5,
point meta max=3.5,
xmin=0,
xmax=17.5,
xtick={0, 8.75, 17.5},
xlabel={$x$},
ymin=2000,
ymax=2200,
ytick={2000,2100,2200},
ylabel={$t$},
axis background/.style={fill=white},
colormap={mymap}{[1pt] rgb(0pt)=(0.2422,0.1504,0.6603); rgb(1pt)=(0.25039,0.164995,0.707614); rgb(2pt)=(0.257771,0.181781,0.751138); rgb(3pt)=(0.264729,0.197757,0.795214); rgb(4pt)=(0.270648,0.214676,0.836371); rgb(5pt)=(0.275114,0.234238,0.870986); rgb(6pt)=(0.2783,0.255871,0.899071); rgb(7pt)=(0.280333,0.278233,0.9221); rgb(8pt)=(0.281338,0.300595,0.941376); rgb(9pt)=(0.281014,0.322757,0.957886); rgb(10pt)=(0.279467,0.344671,0.971676); rgb(11pt)=(0.275971,0.366681,0.982905); rgb(12pt)=(0.269914,0.3892,0.9906); rgb(13pt)=(0.260243,0.412329,0.995157); rgb(14pt)=(0.244033,0.435833,0.998833); rgb(15pt)=(0.220643,0.460257,0.997286); rgb(16pt)=(0.196333,0.484719,0.989152); rgb(17pt)=(0.183405,0.507371,0.979795); rgb(18pt)=(0.178643,0.528857,0.968157); rgb(19pt)=(0.176438,0.549905,0.952019); rgb(20pt)=(0.168743,0.570262,0.935871); rgb(21pt)=(0.154,0.5902,0.9218); rgb(22pt)=(0.146029,0.609119,0.907857); rgb(23pt)=(0.138024,0.627629,0.89729); rgb(24pt)=(0.124814,0.645929,0.888343); rgb(25pt)=(0.111252,0.6635,0.876314); rgb(26pt)=(0.0952095,0.679829,0.859781); rgb(27pt)=(0.0688714,0.694771,0.839357); rgb(28pt)=(0.0296667,0.708167,0.816333); rgb(29pt)=(0.00357143,0.720267,0.7917); rgb(30pt)=(0.00665714,0.731214,0.766014); rgb(31pt)=(0.0433286,0.741095,0.73941); rgb(32pt)=(0.0963952,0.75,0.712038); rgb(33pt)=(0.140771,0.7584,0.684157); rgb(34pt)=(0.1717,0.766962,0.655443); rgb(35pt)=(0.193767,0.775767,0.6251); rgb(36pt)=(0.216086,0.7843,0.5923); rgb(37pt)=(0.246957,0.791795,0.556743); rgb(38pt)=(0.290614,0.79729,0.518829); rgb(39pt)=(0.340643,0.8008,0.478857); rgb(40pt)=(0.3909,0.802871,0.435448); rgb(41pt)=(0.445629,0.802419,0.390919); rgb(42pt)=(0.5044,0.7993,0.348); rgb(43pt)=(0.561562,0.794233,0.304481); rgb(44pt)=(0.617395,0.787619,0.261238); rgb(45pt)=(0.671986,0.779271,0.2227); rgb(46pt)=(0.7242,0.769843,0.191029); rgb(47pt)=(0.773833,0.759805,0.16461); rgb(48pt)=(0.820314,0.749814,0.153529); rgb(49pt)=(0.863433,0.7406,0.159633); rgb(50pt)=(0.903543,0.733029,0.177414); rgb(51pt)=(0.939257,0.728786,0.209957); rgb(52pt)=(0.972757,0.729771,0.239443); rgb(53pt)=(0.995648,0.743371,0.237148); rgb(54pt)=(0.996986,0.765857,0.219943); rgb(55pt)=(0.995205,0.789252,0.202762); rgb(56pt)=(0.9892,0.813567,0.188533); rgb(57pt)=(0.978629,0.838629,0.176557); rgb(58pt)=(0.967648,0.8639,0.16429); rgb(59pt)=(0.96101,0.889019,0.153676); rgb(60pt)=(0.959671,0.913457,0.142257); rgb(61pt)=(0.962795,0.937338,0.12651); rgb(62pt)=(0.969114,0.960629,0.106362); rgb(63pt)=(0.9769,0.9839,0.0805)},
colorbar,
colorbar style={ylabel={$u(x,t)$}, y label style={at={(2.10,0.5)}}, width={0.15in}, ytick={-3,-2, -1,0,1,2, 3}, at={(2.75,1.0)},
height=2*\pgfkeysvalueof{/pgfplots/parent axis height}+0.9in},
axis on top
]
\addplot[thick] graphics[xmin=0,ymin=2000,xmax=17.5,ymax=2200] {ksoddperiodic_L175plot.jpg};
\end{axis}

\begin{axis}[%
width=1.6in,
height=2.1in,
at={(2.6in,0in)},
scale only axis,
point meta min=-3.5,
point meta max=3.5,
xmin=0,
xmax=18.2,
xtick={0, 9.1, 18.2},
xlabel={$x$},
ymin=2000,
ymax=2200,
ytick={2000,2100,2200},
ylabel={$t$},
axis background/.style={fill=white},
colormap={mymap}{[1pt] rgb(0pt)=(0.2422,0.1504,0.6603); rgb(1pt)=(0.25039,0.164995,0.707614); rgb(2pt)=(0.257771,0.181781,0.751138); rgb(3pt)=(0.264729,0.197757,0.795214); rgb(4pt)=(0.270648,0.214676,0.836371); rgb(5pt)=(0.275114,0.234238,0.870986); rgb(6pt)=(0.2783,0.255871,0.899071); rgb(7pt)=(0.280333,0.278233,0.9221); rgb(8pt)=(0.281338,0.300595,0.941376); rgb(9pt)=(0.281014,0.322757,0.957886); rgb(10pt)=(0.279467,0.344671,0.971676); rgb(11pt)=(0.275971,0.366681,0.982905); rgb(12pt)=(0.269914,0.3892,0.9906); rgb(13pt)=(0.260243,0.412329,0.995157); rgb(14pt)=(0.244033,0.435833,0.998833); rgb(15pt)=(0.220643,0.460257,0.997286); rgb(16pt)=(0.196333,0.484719,0.989152); rgb(17pt)=(0.183405,0.507371,0.979795); rgb(18pt)=(0.178643,0.528857,0.968157); rgb(19pt)=(0.176438,0.549905,0.952019); rgb(20pt)=(0.168743,0.570262,0.935871); rgb(21pt)=(0.154,0.5902,0.9218); rgb(22pt)=(0.146029,0.609119,0.907857); rgb(23pt)=(0.138024,0.627629,0.89729); rgb(24pt)=(0.124814,0.645929,0.888343); rgb(25pt)=(0.111252,0.6635,0.876314); rgb(26pt)=(0.0952095,0.679829,0.859781); rgb(27pt)=(0.0688714,0.694771,0.839357); rgb(28pt)=(0.0296667,0.708167,0.816333); rgb(29pt)=(0.00357143,0.720267,0.7917); rgb(30pt)=(0.00665714,0.731214,0.766014); rgb(31pt)=(0.0433286,0.741095,0.73941); rgb(32pt)=(0.0963952,0.75,0.712038); rgb(33pt)=(0.140771,0.7584,0.684157); rgb(34pt)=(0.1717,0.766962,0.655443); rgb(35pt)=(0.193767,0.775767,0.6251); rgb(36pt)=(0.216086,0.7843,0.5923); rgb(37pt)=(0.246957,0.791795,0.556743); rgb(38pt)=(0.290614,0.79729,0.518829); rgb(39pt)=(0.340643,0.8008,0.478857); rgb(40pt)=(0.3909,0.802871,0.435448); rgb(41pt)=(0.445629,0.802419,0.390919); rgb(42pt)=(0.5044,0.7993,0.348); rgb(43pt)=(0.561562,0.794233,0.304481); rgb(44pt)=(0.617395,0.787619,0.261238); rgb(45pt)=(0.671986,0.779271,0.2227); rgb(46pt)=(0.7242,0.769843,0.191029); rgb(47pt)=(0.773833,0.759805,0.16461); rgb(48pt)=(0.820314,0.749814,0.153529); rgb(49pt)=(0.863433,0.7406,0.159633); rgb(50pt)=(0.903543,0.733029,0.177414); rgb(51pt)=(0.939257,0.728786,0.209957); rgb(52pt)=(0.972757,0.729771,0.239443); rgb(53pt)=(0.995648,0.743371,0.237148); rgb(54pt)=(0.996986,0.765857,0.219943); rgb(55pt)=(0.995205,0.789252,0.202762); rgb(56pt)=(0.9892,0.813567,0.188533); rgb(57pt)=(0.978629,0.838629,0.176557); rgb(58pt)=(0.967648,0.8639,0.16429); rgb(59pt)=(0.96101,0.889019,0.153676); rgb(60pt)=(0.959671,0.913457,0.142257); rgb(61pt)=(0.962795,0.937338,0.12651); rgb(62pt)=(0.969114,0.960629,0.106362); rgb(63pt)=(0.9769,0.9839,0.0805)},
axis on top
]
\addplot[thick] graphics[xmin=0,ymin=2000,xmax=18.2,ymax=2200] {ksoddperiodic_L182plot.jpg};
\end{axis}

\begin{axis}[%
width=1.6in,
height=2.1in,
at={(0in,-3in)},
scale only axis,
point meta min=-3.5,
point meta max=3.5,
xmin=0,
xmax=41,
xtick={0, 20.5, 41},
xlabel={$x$},
ymin=2000,
ymax=2200,
ytick={2000,2100,2200},
ylabel={$t$},
axis background/.style={fill=white},
colormap={mymap}{[1pt] rgb(0pt)=(0.2422,0.1504,0.6603); rgb(1pt)=(0.25039,0.164995,0.707614); rgb(2pt)=(0.257771,0.181781,0.751138); rgb(3pt)=(0.264729,0.197757,0.795214); rgb(4pt)=(0.270648,0.214676,0.836371); rgb(5pt)=(0.275114,0.234238,0.870986); rgb(6pt)=(0.2783,0.255871,0.899071); rgb(7pt)=(0.280333,0.278233,0.9221); rgb(8pt)=(0.281338,0.300595,0.941376); rgb(9pt)=(0.281014,0.322757,0.957886); rgb(10pt)=(0.279467,0.344671,0.971676); rgb(11pt)=(0.275971,0.366681,0.982905); rgb(12pt)=(0.269914,0.3892,0.9906); rgb(13pt)=(0.260243,0.412329,0.995157); rgb(14pt)=(0.244033,0.435833,0.998833); rgb(15pt)=(0.220643,0.460257,0.997286); rgb(16pt)=(0.196333,0.484719,0.989152); rgb(17pt)=(0.183405,0.507371,0.979795); rgb(18pt)=(0.178643,0.528857,0.968157); rgb(19pt)=(0.176438,0.549905,0.952019); rgb(20pt)=(0.168743,0.570262,0.935871); rgb(21pt)=(0.154,0.5902,0.9218); rgb(22pt)=(0.146029,0.609119,0.907857); rgb(23pt)=(0.138024,0.627629,0.89729); rgb(24pt)=(0.124814,0.645929,0.888343); rgb(25pt)=(0.111252,0.6635,0.876314); rgb(26pt)=(0.0952095,0.679829,0.859781); rgb(27pt)=(0.0688714,0.694771,0.839357); rgb(28pt)=(0.0296667,0.708167,0.816333); rgb(29pt)=(0.00357143,0.720267,0.7917); rgb(30pt)=(0.00665714,0.731214,0.766014); rgb(31pt)=(0.0433286,0.741095,0.73941); rgb(32pt)=(0.0963952,0.75,0.712038); rgb(33pt)=(0.140771,0.7584,0.684157); rgb(34pt)=(0.1717,0.766962,0.655443); rgb(35pt)=(0.193767,0.775767,0.6251); rgb(36pt)=(0.216086,0.7843,0.5923); rgb(37pt)=(0.246957,0.791795,0.556743); rgb(38pt)=(0.290614,0.79729,0.518829); rgb(39pt)=(0.340643,0.8008,0.478857); rgb(40pt)=(0.3909,0.802871,0.435448); rgb(41pt)=(0.445629,0.802419,0.390919); rgb(42pt)=(0.5044,0.7993,0.348); rgb(43pt)=(0.561562,0.794233,0.304481); rgb(44pt)=(0.617395,0.787619,0.261238); rgb(45pt)=(0.671986,0.779271,0.2227); rgb(46pt)=(0.7242,0.769843,0.191029); rgb(47pt)=(0.773833,0.759805,0.16461); rgb(48pt)=(0.820314,0.749814,0.153529); rgb(49pt)=(0.863433,0.7406,0.159633); rgb(50pt)=(0.903543,0.733029,0.177414); rgb(51pt)=(0.939257,0.728786,0.209957); rgb(52pt)=(0.972757,0.729771,0.239443); rgb(53pt)=(0.995648,0.743371,0.237148); rgb(54pt)=(0.996986,0.765857,0.219943); rgb(55pt)=(0.995205,0.789252,0.202762); rgb(56pt)=(0.9892,0.813567,0.188533); rgb(57pt)=(0.978629,0.838629,0.176557); rgb(58pt)=(0.967648,0.8639,0.16429); rgb(59pt)=(0.96101,0.889019,0.153676); rgb(60pt)=(0.959671,0.913457,0.142257); rgb(61pt)=(0.962795,0.937338,0.12651); rgb(62pt)=(0.969114,0.960629,0.106362); rgb(63pt)=(0.9769,0.9839,0.0805)},
axis on top
]
\addplot[thick] graphics[xmin=0,ymin=2000,xmax=41,ymax=2200] {ksoddperiodic_L41plot.jpg};
\end{axis}

\begin{axis}[%
width=1.6in,
height=2.1in,
at={(2.6in,-3in)},
scale only axis,
point meta min=-3.5,
point meta max=3.5,
xmin=0,
xmax=100,
xtick={0, 50, 100},
xlabel={$x$},
ymin=2000,
ymax=2200,
ytick={2000,2100,2200},
ylabel={$t$},
axis background/.style={fill=white},
colormap={mymap}{[1pt] rgb(0pt)=(0.2422,0.1504,0.6603); rgb(1pt)=(0.25039,0.164995,0.707614); rgb(2pt)=(0.257771,0.181781,0.751138); rgb(3pt)=(0.264729,0.197757,0.795214); rgb(4pt)=(0.270648,0.214676,0.836371); rgb(5pt)=(0.275114,0.234238,0.870986); rgb(6pt)=(0.2783,0.255871,0.899071); rgb(7pt)=(0.280333,0.278233,0.9221); rgb(8pt)=(0.281338,0.300595,0.941376); rgb(9pt)=(0.281014,0.322757,0.957886); rgb(10pt)=(0.279467,0.344671,0.971676); rgb(11pt)=(0.275971,0.366681,0.982905); rgb(12pt)=(0.269914,0.3892,0.9906); rgb(13pt)=(0.260243,0.412329,0.995157); rgb(14pt)=(0.244033,0.435833,0.998833); rgb(15pt)=(0.220643,0.460257,0.997286); rgb(16pt)=(0.196333,0.484719,0.989152); rgb(17pt)=(0.183405,0.507371,0.979795); rgb(18pt)=(0.178643,0.528857,0.968157); rgb(19pt)=(0.176438,0.549905,0.952019); rgb(20pt)=(0.168743,0.570262,0.935871); rgb(21pt)=(0.154,0.5902,0.9218); rgb(22pt)=(0.146029,0.609119,0.907857); rgb(23pt)=(0.138024,0.627629,0.89729); rgb(24pt)=(0.124814,0.645929,0.888343); rgb(25pt)=(0.111252,0.6635,0.876314); rgb(26pt)=(0.0952095,0.679829,0.859781); rgb(27pt)=(0.0688714,0.694771,0.839357); rgb(28pt)=(0.0296667,0.708167,0.816333); rgb(29pt)=(0.00357143,0.720267,0.7917); rgb(30pt)=(0.00665714,0.731214,0.766014); rgb(31pt)=(0.0433286,0.741095,0.73941); rgb(32pt)=(0.0963952,0.75,0.712038); rgb(33pt)=(0.140771,0.7584,0.684157); rgb(34pt)=(0.1717,0.766962,0.655443); rgb(35pt)=(0.193767,0.775767,0.6251); rgb(36pt)=(0.216086,0.7843,0.5923); rgb(37pt)=(0.246957,0.791795,0.556743); rgb(38pt)=(0.290614,0.79729,0.518829); rgb(39pt)=(0.340643,0.8008,0.478857); rgb(40pt)=(0.3909,0.802871,0.435448); rgb(41pt)=(0.445629,0.802419,0.390919); rgb(42pt)=(0.5044,0.7993,0.348); rgb(43pt)=(0.561562,0.794233,0.304481); rgb(44pt)=(0.617395,0.787619,0.261238); rgb(45pt)=(0.671986,0.779271,0.2227); rgb(46pt)=(0.7242,0.769843,0.191029); rgb(47pt)=(0.773833,0.759805,0.16461); rgb(48pt)=(0.820314,0.749814,0.153529); rgb(49pt)=(0.863433,0.7406,0.159633); rgb(50pt)=(0.903543,0.733029,0.177414); rgb(51pt)=(0.939257,0.728786,0.209957); rgb(52pt)=(0.972757,0.729771,0.239443); rgb(53pt)=(0.995648,0.743371,0.237148); rgb(54pt)=(0.996986,0.765857,0.219943); rgb(55pt)=(0.995205,0.789252,0.202762); rgb(56pt)=(0.9892,0.813567,0.188533); rgb(57pt)=(0.978629,0.838629,0.176557); rgb(58pt)=(0.967648,0.8639,0.16429); rgb(59pt)=(0.96101,0.889019,0.153676); rgb(60pt)=(0.959671,0.913457,0.142257); rgb(61pt)=(0.962795,0.937338,0.12651); rgb(62pt)=(0.969114,0.960629,0.106362); rgb(63pt)=(0.9769,0.9839,0.0805)},
axis on top
]
\addplot[thick] graphics[xmin=0,ymin=2000,xmax=100,ymax=2200] {ksoddperiodic_L100plot.jpg};
\end{axis}

\end{tikzpicture}

%% file: ksoddperiodic_lyaps_L20varyT_scaled.tex
\definecolor{mycolor1}{rgb}{0.00000,0.44700,0.74100}%
\definecolor{mycolor2}{rgb}{0.85000,0.32500,0.09800}%
\definecolor{mycolor3}{rgb}{0.92900,0.69400,0.12500}%
\definecolor{mycolor4}{rgb}{0.49400,0.18400,0.55600}%
\definecolor{mycolor5}{rgb}{0.46600,0.67400,0.18800}%
\definecolor{mycolor6}{rgb}{0.30100,0.74500,0.93300}%
\definecolor{mycolor7}{rgb}{0.63500,0.07800,0.18400}%
\begin{tikzpicture}

\begin{axis}[%
width=4.25in,
height=1.75in,
at={(0in,0in)},
scale only axis,
xmode=log,
xmin=0.0001,
xmax=90,
xminorticks=true,
xlabel={$T$},
xmajorgrids,
xminorgrids,
ymin=-8.00998678774122,
ymax=0.912399535314118,
ytick={-7.60090270954199,-6.90775627898064,-5.29834236561059,-4.60527017099142,-2.99822295029797,0,2.99822295029797},
yticklabels={{$-100$},{$-50$},{$-10$},{$-5$},{$-1$},{$0$},{$1$}},
ylabel={$\lambda_i$},
ylabel shift=-4pt,
ymajorgrids,
axis background/.style={fill=white}
]
\addplot [color=mycolor1,solid,line width=1pt,mark size=1pt,mark=*,mark options={solid},forget plot]
  table[row sep=crcr]{%
0.0001	-5.41577255247554\\
0.0002	-4.52332258359843\\
0.0003	-4.0145171309548\\
0.0004	-3.59479749147296\\
0.0005	-3.56062358132666\\
0.0006	-3.59988983258269\\
0.0007	-3.52142404594203\\
0.0008	-3.28021751932703\\
0.0009	-3.03589675137414\\
0.001	-3.39327295047854\\
0.002	-2.14831669307802\\
0.003	-2.05582238469554\\
0.004	-2.21608427933452\\
0.005	-0.575880152186128\\
0.006	-0.666373493429845\\
0.007	0.912399535314118\\
0.008	-0.441871765248287\\
0.009	-0.170771713552398\\
0.01	-0.182918059103128\\
0.02	-0.0281080466691806\\
0.03	-0.276303275824777\\
0.04	0.212465424834524\\
0.05	0.403024940833778\\
0.06	-0.54860837062566\\
0.07	0.319078112965522\\
0.08	0.379943535312701\\
0.09	0.056410924836234\\
0.1	0.32362509521812\\
0.2	0.187186639570424\\
0.3	0.0767269725274962\\
0.4	0.0318006888086905\\
0.5	0.016234149911999\\
0.6	0.0594392747191682\\
0.7	0.0680658620697823\\
0.8	0.00718548536730174\\
0.9	-0.000998222084220704\\
1	0.00927565409001293\\
2	0.0238865534536047\\
3	0.00847772614800461\\
4	0.0154950439418792\\
5	0.00427058181893615\\
6	0.00429016563953211\\
7	-0.00298329757474126\\
8	-0.00399401658112065\\
9	0.0056006998196462\\
10	-0.00434421863580087\\
20	-0.000671358379567325\\
30	-0.00128191374890455\\
40	-0.00138299455913035\\
50	-0.00712174949813689\\
60	-0.00295894058224821\\
70	-0.00407386593141701\\
80	-0.00359082268331256\\
90	-0.00280307022928339\\
};
\addplot [color=mycolor2,solid,line width=1pt,mark size=1pt,mark=*,mark options={solid},forget plot]
  table[row sep=crcr]{%
0.0001	-5.50426791913721\\
0.0002	-4.7100162150045\\
0.0003	-4.42366236489114\\
0.0004	-4.20870928124406\\
0.0005	-3.62667193618215\\
0.0006	-3.51734714536431\\
0.0007	-3.57031995964598\\
0.0008	-3.29820131177282\\
0.0009	-3.08793799753385\\
0.001	-3.15635751487056\\
0.002	-2.41884172630875\\
0.003	-2.58022006674761\\
0.004	-1.27280370614677\\
0.005	-2.02638225724152\\
0.006	-2.0399708858526\\
0.007	-1.80862317670486\\
0.008	-1.16006676287511\\
0.009	-1.69709243704807\\
0.01	-1.26505546054948\\
0.02	-1.17577665116439\\
0.03	0.154724113003285\\
0.04	-0.304245599534674\\
0.05	-0.162129260382481\\
0.06	0.12837093652203\\
0.07	-1.03942006314254\\
0.08	-0.714068396362627\\
0.09	-0.595892820622738\\
0.1	-0.73405079588349\\
0.2	-0.417754082488482\\
0.3	-0.41497210649199\\
0.4	-0.302189458055484\\
0.5	-0.164370084048705\\
0.6	-0.129968549861651\\
0.7	-0.186706307606206\\
0.8	-0.178871687262304\\
0.9	-0.243920208153784\\
1	-0.349381220922747\\
2	-0.291461008205544\\
3	-0.30156177590571\\
4	-0.346247171478024\\
5	-0.344168153502491\\
6	-0.310715055977911\\
7	-0.302008393405227\\
8	-0.295149536990177\\
9	-0.302672842746833\\
10	-0.289558208966653\\
20	-0.292374594857417\\
30	-0.294893931576478\\
40	-0.273788162132313\\
50	-0.241159993598719\\
60	-0.229216938668092\\
70	-0.224531329730174\\
80	-0.209642129553581\\
90	-0.189635720517815\\
};
\addplot [color=mycolor3,solid,line width=1pt,mark size=1pt,mark=*,mark options={solid},forget plot]
  table[row sep=crcr]{%
0.0001	-4.81838995332363\\
0.0002	-4.76709511403128\\
0.0003	-4.27947899079305\\
0.0004	-4.15063721655202\\
0.0005	-3.9894187385612\\
0.0006	-3.72633069485147\\
0.0007	-3.90180294707624\\
0.0008	-3.12574809716188\\
0.0009	-3.18986996180597\\
0.001	-3.24252918416049\\
0.002	-2.82034030334339\\
0.003	-2.51037135977017\\
0.004	-2.23175918191116\\
0.005	-1.81659575700486\\
0.006	-2.33681210935301\\
0.007	-1.97843442026397\\
0.008	-1.71295587343403\\
0.009	-1.79333397172545\\
0.01	-2.0744388371037\\
0.02	-1.58053994200035\\
0.03	-1.60767671250182\\
0.04	-1.31171707387931\\
0.05	-1.32168345596777\\
0.06	-1.10023841024268\\
0.07	-0.68405538743464\\
0.08	-0.748851638980104\\
0.09	-0.987368233037956\\
0.1	-0.790301127983739\\
0.2	-0.998412803037934\\
0.3	-0.8020495207188\\
0.4	-0.862585427914166\\
0.5	-0.970761320174154\\
0.6	-0.964704916158336\\
0.7	-0.945374081412997\\
0.8	-0.89806400525818\\
0.9	-0.855157241821511\\
1	-0.801765821939627\\
2	-0.817699048434392\\
3	-0.7986876247746\\
4	-0.765482915610174\\
5	-0.761911029363781\\
6	-0.788044742073341\\
7	-0.794239341796202\\
8	-0.793320799765063\\
9	-0.783096405113816\\
10	-0.809499218849058\\
20	-0.79998438630908\\
30	-0.749982727776937\\
40	-0.674391464053043\\
50	-0.57752713404534\\
60	-0.521933340611038\\
70	-0.4536503197757\\
80	-0.404337485359979\\
90	-0.381411173189036\\
};
\addplot [color=mycolor4,solid,line width=1pt,mark size=1pt,mark=*,mark options={solid},forget plot]
  table[row sep=crcr]{%
0.0001	-4.85943911634356\\
0.0002	-4.39944198076882\\
0.0003	-4.18054397363785\\
0.0004	-3.97553342237642\\
0.0005	-4.17517567426732\\
0.0006	-4.01633723309225\\
0.0007	-3.50138418512935\\
0.0008	-3.22640207466534\\
0.0009	-2.85341682607233\\
0.001	-3.37085872082624\\
0.002	-2.78037032545332\\
0.003	-1.66230625914128\\
0.004	-2.46762681239492\\
0.005	-2.45322116392006\\
0.006	-2.04356039149462\\
0.007	-2.16978071117954\\
0.008	-2.21805587621368\\
0.009	-1.90483103670784\\
0.01	-2.17704379707688\\
0.02	-1.58809195637561\\
0.03	-1.72462903379678\\
0.04	-1.81381945834022\\
0.05	-1.72093547829064\\
0.06	-1.60266316547081\\
0.07	-1.59892334907139\\
0.08	-1.70871440900834\\
0.09	-1.55897452031926\\
0.1	-1.57152181998528\\
0.2	-1.48513794947558\\
0.3	-1.57107825421943\\
0.4	-1.52860561887247\\
0.5	-1.5324280365932\\
0.6	-1.53793730505149\\
0.7	-1.51307479166427\\
0.8	-1.53277416502719\\
0.9	-1.5272119033757\\
1	-1.50289024577141\\
2	-1.51628190500886\\
3	-1.5139823128843\\
4	-1.51318375720677\\
5	-1.51091998554067\\
6	-1.50816973864209\\
7	-1.50604066145902\\
8	-1.50649632064055\\
9	-1.50938882900159\\
10	-1.50092714696141\\
20	-1.49973179344798\\
30	-1.47915116938042\\
40	-1.30675842360215\\
50	-1.12153652333607\\
60	-0.963178711743733\\
70	-0.855338694992098\\
80	-0.751320485602871\\
90	-0.68863929951963\\
};
\addplot [color=mycolor5,solid,line width=1pt,mark size=1pt,mark=*,mark options={solid},forget plot]
  table[row sep=crcr]{%
0.0001	-5.25119006276512\\
0.0002	-5.21545671945019\\
0.0003	-4.62929084654924\\
0.0004	-4.31958718062108\\
0.0005	-4.10906109401785\\
0.0006	-3.29972107035975\\
0.0007	-3.87018846593019\\
0.0008	-3.65778297211068\\
0.0009	-3.35516585427977\\
0.001	-3.20462038930246\\
0.002	-2.90574438305048\\
0.003	-3.10525269845975\\
0.004	-2.20994877669621\\
0.005	-2.49420979652158\\
0.006	-1.76871179270122\\
0.007	-2.02655287436096\\
0.008	-2.24820840003275\\
0.009	-2.02210485668469\\
0.01	-1.720130897014\\
0.02	-1.85336931171737\\
0.03	-1.79329206579235\\
0.04	-1.85003833619467\\
0.05	-1.86461137239127\\
0.06	-1.90913770164971\\
0.07	-1.82749885515082\\
0.08	-1.8101508163243\\
0.09	-1.77505099138946\\
0.1	-1.82516419571021\\
0.2	-1.79502659861794\\
0.3	-1.77719772322429\\
0.4	-1.79804200273027\\
0.5	-1.7998740097102\\
0.6	-1.79007898838957\\
0.7	-1.79109456828698\\
0.8	-1.80143013353577\\
0.9	-1.79498430707659\\
1	-1.799849903375\\
2	-1.79696017740637\\
3	-1.79705299112252\\
4	-1.79375688899384\\
5	-1.79637076850506\\
6	-1.79379790251354\\
7	-1.7945940327774\\
8	-1.79471881463007\\
9	-1.79359427592722\\
10	-1.79315205424731\\
20	-1.78552331041747\\
30	-1.69093991889078\\
40	-1.47681405877306\\
50	-1.2659258182016\\
60	-1.1075729511501\\
70	-0.972699524022171\\
80	-0.855339724765111\\
90	-0.779630993788391\\
};
\addplot [color=mycolor6,solid,line width=1pt,mark size=1pt,mark=*,mark options={solid},forget plot]
  table[row sep=crcr]{%
0.0001	-5.67536965779939\\
0.0002	-5.01471648479097\\
0.0003	-4.36506252768181\\
0.0004	-4.08171420919261\\
0.0005	-4.43418372849763\\
0.0006	-3.89253252672504\\
0.0007	-3.78429655644729\\
0.0008	-4.02296240759221\\
0.0009	-4.0427625757311\\
0.001	-3.76123451077836\\
0.002	-3.08844611825098\\
0.003	-2.7369897445218\\
0.004	-2.47518327785463\\
0.005	-2.60611009165854\\
0.006	-2.29490692836002\\
0.007	-2.31489410955881\\
0.008	-2.21840396788098\\
0.009	-2.43059025392427\\
0.01	-2.2008052665875\\
0.02	-2.07085361000935\\
0.03	-1.90146923747375\\
0.04	-1.90182157751497\\
0.05	-1.92927660918575\\
0.06	-1.89331228141987\\
0.07	-1.93720087018921\\
0.08	-1.82323180527079\\
0.09	-2.00031659978931\\
0.1	-1.88644636660246\\
0.2	-1.8715239180161\\
0.3	-1.85595792386618\\
0.4	-1.87826501345381\\
0.5	-1.84811790271148\\
0.6	-1.85167406052034\\
0.7	-1.85604052505765\\
0.8	-1.83840376729416\\
0.9	-1.86140181234165\\
1	-1.85239550351533\\
2	-1.85424838804933\\
3	-1.85252519676377\\
4	-1.85188061401099\\
5	-1.85077068551391\\
6	-1.85298392395197\\
7	-1.85193100542921\\
8	-1.85317271920897\\
9	-1.85264388450708\\
10	-1.85434878874541\\
20	-1.84358824568616\\
30	-1.75832158068441\\
40	-1.54934310227175\\
50	-1.34395704726312\\
60	-1.17415719304451\\
70	-1.02967144328329\\
80	-0.9184872590836\\
90	-0.834661080739887\\
};
\addplot [color=mycolor7,solid,line width=1pt,mark size=1pt,mark=*,mark options={solid},forget plot]
  table[row sep=crcr]{%
0.0001	-5.34240537072158\\
0.0002	-4.56253550335305\\
0.0003	-4.3763835516128\\
0.0004	-3.97886436433835\\
0.0005	-4.14114660225149\\
0.0006	-3.83647014140832\\
0.0007	-3.85342683486856\\
0.0008	-4.13894385658963\\
0.0009	-3.95169842339704\\
0.001	-4.12971507577503\\
0.002	-3.28366683411139\\
0.003	-3.28450494327191\\
0.004	-2.51977730030558\\
0.005	-2.86776757773489\\
0.006	-2.42234344562315\\
0.007	-2.64453889988182\\
0.008	-2.78712583224494\\
0.009	-2.2833493258604\\
0.01	-2.72874362635404\\
0.02	-2.48978587050267\\
0.03	-2.44611975984331\\
0.04	-2.45303458201064\\
0.05	-2.33826190035277\\
0.06	-2.37189079002521\\
0.07	-2.34881584759386\\
0.08	-2.33915317998092\\
0.09	-2.27907419979177\\
0.1	-2.30362984498629\\
0.2	-2.30101118651004\\
0.3	-2.29336390655928\\
0.4	-2.28580225030602\\
0.5	-2.27810119454721\\
0.6	-2.28052793475296\\
0.7	-2.28626248658365\\
0.8	-2.27488172057263\\
0.9	-2.27149511260904\\
1	-2.27760276718135\\
2	-2.27134213660606\\
3	-2.26997683295425\\
4	-2.27199287632424\\
5	-2.27063667662797\\
6	-2.27074187294737\\
7	-2.27154453449815\\
8	-2.26999324990849\\
9	-2.26997482937265\\
10	-2.26860669278106\\
20	-2.16903288791034\\
30	-1.83862895224765\\
40	-1.60653194206113\\
50	-1.40175450662727\\
60	-1.22331646410334\\
70	-1.08760091052236\\
80	-0.968778895973921\\
90	-0.881014628211398\\
};
\addplot [color=mycolor1,solid,line width=1pt,mark size=1pt,mark=*,mark options={solid},forget plot]
  table[row sep=crcr]{%
0.0001	-5.5787727537924\\
0.0002	-4.94581602320787\\
0.0003	-4.84122589832956\\
0.0004	-4.92201611204271\\
0.0005	-4.36531002133771\\
0.0006	-4.13850062047078\\
0.0007	-4.20050200101604\\
0.0008	-3.83820521727866\\
0.0009	-4.16143580932621\\
0.001	-3.63381039151916\\
0.002	-2.56179811002431\\
0.003	-2.91843596522514\\
0.004	-3.39404676152478\\
0.005	-2.72634359819405\\
0.006	-2.7127078757444\\
0.007	-2.71676470749407\\
0.008	-3.13335007474449\\
0.009	-2.81385351299773\\
0.01	-2.42777798280858\\
0.02	-2.6399781930436\\
0.03	-2.53542672186479\\
0.04	-2.3696741764426\\
0.05	-2.40191963463558\\
0.06	-2.42644861309256\\
0.07	-2.3942902562668\\
0.08	-2.42311600915561\\
0.09	-2.35101454969038\\
0.1	-2.45131893214905\\
0.2	-2.3794306906539\\
0.3	-2.39102282159942\\
0.4	-2.35445790678041\\
0.5	-2.3726450422534\\
0.6	-2.3726282937383\\
0.7	-2.38022293716278\\
0.8	-2.37298592369324\\
0.9	-2.37491540643138\\
1	-2.37151539734779\\
2	-2.37407513573152\\
3	-2.37265554228563\\
4	-2.37136596407646\\
5	-2.37149877828404\\
6	-2.37186905255698\\
7	-2.36985894901081\\
8	-2.37092438399477\\
9	-2.36886248168898\\
10	-2.37086968226344\\
20	-2.27109940632703\\
30	-1.90469184245855\\
40	-1.65833313405328\\
50	-1.44839160923094\\
60	-1.27325115285474\\
70	-1.13741759253004\\
80	-1.01306787626738\\
90	-0.925137451577137\\
};
\addplot [color=mycolor2,solid,line width=1pt,mark size=1pt,mark=*,mark options={solid},forget plot]
  table[row sep=crcr]{%
0.0001	-6.17103642923764\\
0.0002	-5.0416669462162\\
0.0003	-4.90613182633922\\
0.0004	-4.86162126295985\\
0.0005	-4.93013097336204\\
0.0006	-3.89275718668617\\
0.0007	-4.22769812963653\\
0.0008	-4.41938091583758\\
0.0009	-4.47189859507106\\
0.001	-4.17948334610046\\
0.002	-3.94941647211968\\
0.003	-3.68466753442195\\
0.004	-3.41927061584374\\
0.005	-3.53166017150104\\
0.006	-3.82138017485047\\
0.007	-3.68555280957763\\
0.008	-3.56582842309093\\
0.009	-3.55834926219281\\
0.01	-3.60318203846696\\
0.02	-3.46756407198908\\
0.03	-3.49795811817006\\
0.04	-3.53800491388701\\
0.05	-3.55337925996966\\
0.06	-3.51114943713585\\
0.07	-3.54058100538468\\
0.08	-3.53391539691317\\
0.09	-3.52168648766175\\
0.1	-3.51565021068343\\
0.2	-3.51128583787065\\
0.3	-3.50565184279513\\
0.4	-3.50853053676125\\
0.5	-3.50951900784407\\
0.6	-3.51180376474688\\
0.7	-3.50593854315119\\
0.8	-3.51148348454496\\
0.9	-3.51101057094141\\
1	-3.50977938806549\\
2	-3.51140611287314\\
3	-3.5103194746486\\
4	-3.51406910269584\\
5	-3.52911674659623\\
6	-3.50311927210861\\
7	-3.41556394892597\\
8	-3.29055699290078\\
9	-3.17162962272816\\
10	-3.06828791681676\\
20	-2.35631676457607\\
30	-1.97217196731505\\
40	-1.70234577125143\\
50	-1.49236667160035\\
60	-1.31694849310209\\
70	-1.17686265050892\\
80	-1.05155275695727\\
90	-0.963741276883302\\
};
\addplot [color=mycolor3,solid,line width=1pt,mark size=1pt,mark=*,mark options={solid},forget plot]
  table[row sep=crcr]{%
0.0001	-5.8781106136802\\
0.0002	-5.27949127769558\\
0.0003	-5.11150133889965\\
0.0004	-5.02750472932376\\
0.0005	-4.62190464113539\\
0.0006	-4.40112577347927\\
0.0007	-4.46010761694925\\
0.0008	-4.7418964304187\\
0.0009	-4.27390750163748\\
0.001	-4.15471679507821\\
0.002	-4.44808609554953\\
0.003	-4.12192323867988\\
0.004	-4.30944415438583\\
0.005	-4.31288137074903\\
0.006	-4.15106226280499\\
0.007	-4.20207568301944\\
0.008	-4.21167643686421\\
0.009	-4.22484753412538\\
0.01	-4.18069144801983\\
0.02	-4.20350201768046\\
0.03	-4.19065479287611\\
0.04	-4.17010408941873\\
0.05	-4.18034528907519\\
0.06	-4.18944379116746\\
0.07	-4.17096149184449\\
0.08	-4.1762932879802\\
0.09	-4.17979406768326\\
0.1	-4.17850213743169\\
0.2	-4.17714117144013\\
0.3	-4.17538399283034\\
0.4	-4.17719580638145\\
0.5	-4.17597066604236\\
0.6	-4.17817495606261\\
0.7	-4.17954997239677\\
0.8	-4.18179895095845\\
0.9	-4.18096215855921\\
1	-4.18263263380253\\
2	-4.18879077211281\\
3	-4.16034248999199\\
4	-3.92257125934466\\
5	-3.7380332394922\\
6	-3.59161820915258\\
7	-3.47236505397285\\
8	-3.34638876773383\\
9	-3.22918441923173\\
10	-3.12585869176913\\
20	-2.42248245608106\\
30	-2.02740913667961\\
40	-1.74391702563585\\
50	-1.53090872560103\\
60	-1.35373995132838\\
70	-1.21803342772182\\
80	-1.09570242585152\\
90	-0.996152126493199\\
};
\addplot [color=mycolor4,solid,line width=1pt,mark size=1pt,mark=*,mark options={solid},forget plot]
  table[row sep=crcr]{%
0.0001	-5.55635267010799\\
0.0002	-5.61618022682001\\
0.0003	-5.4533044770314\\
0.0004	-4.79442208720622\\
0.0005	-5.19322040555416\\
0.0006	-5.35006177149403\\
0.0007	-5.05978228945836\\
0.0008	-4.77286867910079\\
0.0009	-4.68029415448749\\
0.001	-4.8446661389748\\
0.002	-4.80366797926381\\
0.003	-4.78073506634732\\
0.004	-4.76578189883177\\
0.005	-4.67932305918063\\
0.006	-4.72884840323334\\
0.007	-4.67537627871458\\
0.008	-4.71085344342132\\
0.009	-4.70195950289621\\
0.01	-4.72355588030736\\
0.02	-4.69802487407658\\
0.03	-4.68591998740747\\
0.04	-4.69282155548998\\
0.05	-4.68513557513079\\
0.06	-4.68495951179488\\
0.07	-4.68460888063405\\
0.08	-4.68387946737604\\
0.09	-4.68481706539737\\
0.1	-4.68377892462264\\
0.2	-4.67899379529828\\
0.3	-4.6785422067948\\
0.4	-4.6807815738304\\
0.5	-4.6848148677895\\
0.6	-4.6874074733517\\
0.7	-4.69089505005825\\
0.8	-4.69400586921052\\
0.9	-4.69676020722729\\
1	-4.70063349513218\\
2	-4.56130661637666\\
3	-4.19428520967408\\
4	-3.95633127767018\\
5	-3.77545850503722\\
6	-3.63866965915441\\
7	-3.51281278284605\\
8	-3.39052459617064\\
9	-3.27849668158822\\
10	-3.1746348826774\\
20	-2.47571919762039\\
30	-2.07556257970732\\
40	-1.78462060383061\\
50	-1.56732033352166\\
60	-1.3926040491227\\
70	-1.2540148236272\\
80	-1.12897571941607\\
90	-1.03045137555019\\
};
\addplot [color=mycolor5,solid,line width=1pt,mark size=1pt,mark=*,mark options={solid},forget plot]
  table[row sep=crcr]{%
0.0001	-5.72138888536309\\
0.0002	-5.95992418865706\\
0.0003	-5.77214197741886\\
0.0004	-5.71357140919439\\
0.0005	-4.91245503055197\\
0.0006	-5.167079850649\\
0.0007	-4.94970543962915\\
0.0008	-5.29209207568951\\
0.0009	-5.06618164694027\\
0.001	-5.27051447825652\\
0.002	-5.09442442030654\\
0.003	-5.14497262889319\\
0.004	-5.14878172048487\\
0.005	-5.10537835581112\\
0.006	-5.16076373163451\\
0.007	-5.12825674563812\\
0.008	-5.08039290621417\\
0.009	-5.12733511424032\\
0.01	-5.13976855185547\\
0.02	-5.10722422480462\\
0.03	-5.11104869547591\\
0.04	-5.10290788961613\\
0.05	-5.11318008282206\\
0.06	-5.10912146381723\\
0.07	-5.1040234234964\\
0.08	-5.10585310094295\\
0.09	-5.10274974271477\\
0.1	-5.10516463179111\\
0.2	-5.10066891470719\\
0.3	-5.10069032332475\\
0.4	-5.10632714828591\\
0.5	-5.1126768879007\\
0.6	-5.12004955258559\\
0.7	-5.13051724253607\\
0.8	-5.14597419851439\\
0.9	-5.17747381740433\\
1	-5.25064958735049\\
2	-4.58621416417677\\
3	-4.22830176391578\\
4	-3.99537062049857\\
5	-3.81493266844724\\
6	-3.67619967978445\\
7	-3.55550989996398\\
8	-3.43343395199529\\
9	-3.32304155242168\\
10	-3.21841150017818\\
20	-2.52037896195242\\
30	-2.11020767200596\\
40	-1.8221051497029\\
50	-1.60316622236853\\
60	-1.42808735750305\\
70	-1.28468577788146\\
80	-1.15657026806178\\
90	-1.05784210924165\\
};
\addplot [color=mycolor6,solid,line width=1pt,mark size=1pt,mark=*,mark options={solid},forget plot]
  table[row sep=crcr]{%
0.0001	-6.1434226554934\\
0.0002	-5.85084148932507\\
0.0003	-5.08907337799516\\
0.0004	-5.74210069100907\\
0.0005	-5.57822801678367\\
0.0006	-5.9074328876882\\
0.0007	-5.56228895387517\\
0.0008	-5.30608324054344\\
0.0009	-5.60265265647317\\
0.001	-5.49053993621339\\
0.002	-5.48936701605587\\
0.003	-5.4842986826023\\
0.004	-5.47076772969203\\
0.005	-5.49168543902638\\
0.006	-5.47965131408205\\
0.007	-5.49868990354769\\
0.008	-5.4973599182923\\
0.009	-5.47669040130956\\
0.01	-5.48498226323994\\
0.02	-5.4675116435828\\
0.03	-5.46773798278703\\
0.04	-5.47140031745507\\
0.05	-5.46586654994733\\
0.06	-5.4712214094857\\
0.07	-5.47053593027652\\
0.08	-5.46906375201909\\
0.09	-5.47136617012147\\
0.1	-5.4673409899698\\
0.2	-5.46483381574538\\
0.3	-5.47134034858023\\
0.4	-5.48224313585394\\
0.5	-5.49769815120767\\
0.6	-5.52795496226847\\
0.7	-5.61948663091142\\
0.8	-5.47710126249661\\
0.9	-5.36465878831055\\
1	-5.28168702504797\\
2	-4.62083363699357\\
3	-4.26054795210606\\
4	-4.03094643087262\\
5	-3.85384789634984\\
6	-3.71538773444261\\
7	-3.59452239449248\\
8	-3.46759224258413\\
9	-3.35729235604787\\
10	-3.25246442932287\\
20	-2.55216157727048\\
30	-2.1432933738581\\
40	-1.84951040788904\\
50	-1.63150181353126\\
60	-1.45750395697555\\
70	-1.31375440488303\\
80	-1.18450324591152\\
90	-1.08825470163479\\
};
\addplot [color=mycolor7,solid,line width=1pt,mark size=1pt,mark=*,mark options={solid},forget plot]
  table[row sep=crcr]{%
0.0001	-6.23903388857025\\
0.0002	-5.95124383688762\\
0.0003	-5.78316095860514\\
0.0004	-5.93161600862636\\
0.0005	-6.00724852428312\\
0.0006	-5.65065442872188\\
0.0007	-6.06296074681821\\
0.0008	-5.84412290728064\\
0.0009	-5.80001933858098\\
0.001	-5.87168214002738\\
0.002	-5.84370028343131\\
0.003	-5.80501493885097\\
0.004	-5.83064069836979\\
0.005	-5.8126216016707\\
0.006	-5.81659012466926\\
0.007	-5.81873653878625\\
0.008	-5.81041383102349\\
0.009	-5.80389961513123\\
0.01	-5.79979370617959\\
0.02	-5.80347346002805\\
0.03	-5.79867138870985\\
0.04	-5.79815213095675\\
0.05	-5.79810432773999\\
0.06	-5.79321330847844\\
0.07	-5.79649989668247\\
0.08	-5.79671023551378\\
0.09	-5.79315118538668\\
0.1	-5.79323537679305\\
0.2	-5.79308491243678\\
0.3	-5.80570447639492\\
0.4	-5.8282088855372\\
0.5	-5.91945473637066\\
0.6	-5.7903699855386\\
0.7	-5.64614040664805\\
0.8	-5.49228772192086\\
0.9	-5.37848518089034\\
1	-5.31177043280866\\
2	-4.65125799332639\\
3	-4.29571111839785\\
4	-4.06909705630317\\
5	-3.89409615068031\\
6	-3.75330822092813\\
7	-3.62489634829682\\
8	-3.50043819469908\\
9	-3.39247791916144\\
10	-3.28801620966669\\
20	-2.58416918643549\\
30	-2.17189967220477\\
40	-1.87743735279085\\
50	-1.65956260737439\\
60	-1.48188266714305\\
70	-1.33987393309895\\
80	-1.21390130901791\\
90	-1.1096509432338\\
};
\addplot [color=mycolor1,solid,line width=1pt,mark size=1pt,mark=*,mark options={solid},forget plot]
  table[row sep=crcr]{%
0.0001	-6.6880580145181\\
0.0002	-6.19544109682881\\
0.0003	-6.6191662814556\\
0.0004	-5.73115547010957\\
0.0005	-6.17468640987682\\
0.0006	-6.18535529565969\\
0.0007	-6.00829850330481\\
0.0008	-6.15942272069088\\
0.0009	-6.09924479810841\\
0.001	-6.14967234465084\\
0.002	-6.14940169511937\\
0.003	-6.09862524486007\\
0.004	-6.1067392123284\\
0.005	-6.10126563388939\\
0.006	-6.10834469579467\\
0.007	-6.10167663781829\\
0.008	-6.0996809538865\\
0.009	-6.09875271639781\\
0.01	-6.0935355875501\\
0.02	-6.09155985508062\\
0.03	-6.08994629564551\\
0.04	-6.0900596104832\\
0.05	-6.09147779120095\\
0.06	-6.08722946826089\\
0.07	-6.08878581278012\\
0.08	-6.08719030369144\\
0.09	-6.08832627083568\\
0.1	-6.08772075787074\\
0.2	-6.09222652149726\\
0.3	-6.11387493108412\\
0.4	-6.15781985412454\\
0.5	-6.00531726760297\\
0.6	-5.82357302581753\\
0.7	-5.68057304624734\\
0.8	-5.53535555793112\\
0.9	-5.40810380392174\\
1	-5.33200643505958\\
2	-4.68536207077091\\
3	-4.33470488720565\\
4	-4.10746144980605\\
5	-3.93461705045837\\
6	-3.78838767486191\\
7	-3.65432643828153\\
8	-3.52985210577594\\
9	-3.41963954892913\\
10	-3.31451539262172\\
20	-2.60818722522215\\
30	-2.19661111086683\\
40	-1.90661000634795\\
50	-1.68555692904124\\
60	-1.50480977964244\\
70	-1.36207107286758\\
80	-1.23258551641423\\
90	-1.13220721880359\\
};
\addplot [color=mycolor2,solid,line width=1pt,mark size=1pt,mark=*,mark options={solid},forget plot]
  table[row sep=crcr]{%
0.0001	-6.90618096197073\\
0.0002	-6.5720460936682\\
0.0003	-6.52490464483754\\
0.0004	-6.42472366647141\\
0.0005	-6.33102501061555\\
0.0006	-6.41810257158601\\
0.0007	-6.44873650670132\\
0.0008	-6.46269451875248\\
0.0009	-6.35703989008608\\
0.001	-6.40639269074196\\
0.002	-6.32900147355661\\
0.003	-6.37757570195434\\
0.004	-6.37751952374931\\
0.005	-6.37824653780129\\
0.006	-6.3651464797299\\
0.007	-6.36749531147368\\
0.008	-6.3701811008675\\
0.009	-6.37490856737355\\
0.01	-6.36202260832439\\
0.02	-6.3603836941411\\
0.03	-6.35497326418096\\
0.04	-6.35475080041811\\
0.05	-6.35670886026503\\
0.06	-6.35598726511555\\
0.07	-6.35439869965728\\
0.08	-6.35716088066283\\
0.09	-6.35701012582573\\
0.1	-6.35644488772733\\
0.2	-6.36820096883223\\
0.3	-6.40002469649109\\
0.4	-6.24707284153474\\
0.5	-6.05501112192012\\
0.6	-5.84704220113857\\
0.7	-5.69644313375578\\
0.8	-5.56832739239029\\
0.9	-5.44098722986907\\
1	-5.35718766995429\\
2	-4.72170876725525\\
3	-4.37270990596922\\
4	-4.14386533916045\\
5	-3.96859485118166\\
6	-3.81571340633143\\
7	-3.6830589461998\\
8	-3.55601231100729\\
9	-3.44228224049288\\
10	-3.33853198576549\\
20	-2.63520128568602\\
30	-2.22252903304729\\
40	-1.9233117234905\\
50	-1.70653099404369\\
60	-1.53061244992765\\
70	-1.38455472096118\\
80	-1.25578585442853\\
90	-1.15119091984862\\
};
\addplot [color=mycolor3,solid,line width=1pt,mark size=1pt,mark=*,mark options={solid},forget plot]
  table[row sep=crcr]{%
0.0001	-6.77207381124983\\
0.0002	-6.45532417200693\\
0.0003	-6.46744504981477\\
0.0004	-6.62628380753138\\
0.0005	-6.67790425580502\\
0.0006	-6.73740414085901\\
0.0007	-6.61929009248044\\
0.0008	-6.61619542901181\\
0.0009	-6.65000322014655\\
0.001	-6.64541923806866\\
0.002	-6.63555528866551\\
0.003	-6.62452948575898\\
0.004	-6.62400802576178\\
0.005	-6.6253438075774\\
0.006	-6.61762187995724\\
0.007	-6.61529283838547\\
0.008	-6.61203566397625\\
0.009	-6.6160235761598\\
0.01	-6.61989027456901\\
0.02	-6.59973821934071\\
0.03	-6.59822621312623\\
0.04	-6.6002298632806\\
0.05	-6.60151398918115\\
0.06	-6.60020168258829\\
0.07	-6.60082674831507\\
0.08	-6.60296534650765\\
0.09	-6.60346926168368\\
0.1	-6.60432965055417\\
0.2	-6.61999346127034\\
0.3	-6.50682243445874\\
0.4	-6.31979235111776\\
0.5	-6.07324341683173\\
0.6	-5.87983082725783\\
0.7	-5.71793967515257\\
0.8	-5.59226052217606\\
0.9	-5.46441703547821\\
1	-5.38155855660029\\
2	-4.75581911465961\\
3	-4.40892194502319\\
4	-4.17755657976866\\
5	-3.99917975743615\\
6	-3.84604541317744\\
7	-3.70905849535715\\
8	-3.57951140014762\\
9	-3.46611140315991\\
10	-3.35963288978652\\
20	-2.65223846406706\\
30	-2.24412316126477\\
40	-1.94416624215747\\
50	-1.7262374439798\\
60	-1.5472415339967\\
70	-1.40323836817285\\
80	-1.27408737829817\\
90	-1.16935388660186\\
};
\addplot [color=mycolor4,solid,line width=1pt,mark size=1pt,mark=*,mark options={solid},forget plot]
  table[row sep=crcr]{%
0.0001	-6.63591077077199\\
0.0002	-6.94115265316406\\
0.0003	-6.73629017748683\\
0.0004	-6.87600345525497\\
0.0005	-6.87150718001975\\
0.0006	-6.72644517703946\\
0.0007	-6.89428860251914\\
0.0008	-6.87819386423995\\
0.0009	-6.88677567935092\\
0.001	-6.85949847446328\\
0.002	-6.86614862511843\\
0.003	-6.84766530376958\\
0.004	-6.84168742720019\\
0.005	-6.83305412214987\\
0.006	-6.85152721582589\\
0.007	-6.85378629808531\\
0.008	-6.83888259095123\\
0.009	-6.83846484724363\\
0.01	-6.83798306702818\\
0.02	-6.82693916349167\\
0.03	-6.82176085722756\\
0.04	-6.82350286481035\\
0.05	-6.82664108067011\\
0.06	-6.82687706029188\\
0.07	-6.8275994219847\\
0.08	-6.83161443483619\\
0.09	-6.83372502348158\\
0.1	-6.83519989275312\\
0.2	-6.81271921527293\\
0.3	-6.62027131099866\\
0.4	-6.33526021022417\\
0.5	-6.09292522555779\\
0.6	-5.90355124459751\\
0.7	-5.73896780688391\\
0.8	-5.60365743131151\\
0.9	-5.48299746001898\\
1	-5.40932780718966\\
2	-4.7974374602842\\
3	-4.44912945975859\\
4	-4.21227941695804\\
5	-4.03141287804763\\
6	-3.87189507995662\\
7	-3.73206028330736\\
8	-3.59938770693961\\
9	-3.48549504224099\\
10	-3.3819170463359\\
20	-2.6713205048732\\
30	-2.26074352950739\\
40	-1.96471351602238\\
50	-1.74493236234056\\
60	-1.56650713888627\\
70	-1.42144502093741\\
80	-1.28873219115279\\
90	-1.18842412472595\\
};
\addplot [color=mycolor5,solid,line width=1pt,mark size=1pt,mark=*,mark options={solid},forget plot]
  table[row sep=crcr]{%
0.0001	-7.1207796663259\\
0.0002	-7.24996123934905\\
0.0003	-7.09397692892394\\
0.0004	-7.03506930165605\\
0.0005	-7.08405073052947\\
0.0006	-7.07981346538906\\
0.0007	-7.05135101621559\\
0.0008	-7.09589707237199\\
0.0009	-7.07284558330477\\
0.001	-7.04331294639982\\
0.002	-7.06091623318534\\
0.003	-7.0638224316317\\
0.004	-7.06141660843623\\
0.005	-7.06381216556564\\
0.006	-7.05008372275698\\
0.007	-7.05158927236348\\
0.008	-7.05298556858924\\
0.009	-7.04932260631801\\
0.01	-7.05239700480849\\
0.02	-7.0347249839085\\
0.03	-7.02913663796424\\
0.04	-7.03287191120217\\
0.05	-7.03651440673088\\
0.06	-7.03741370744991\\
0.07	-7.03970294000406\\
0.08	-7.04571445921353\\
0.09	-7.04887449843622\\
0.1	-7.05107948328599\\
0.2	-6.91135225849738\\
0.3	-6.64077454883542\\
0.4	-6.35111316868815\\
0.5	-6.12075203384688\\
0.6	-5.91934328811252\\
0.7	-5.75730491476722\\
0.8	-5.63191118493626\\
0.9	-5.51281563377842\\
1	-5.44513124319326\\
2	-4.8328289361116\\
3	-4.48179774807518\\
4	-4.24204882423782\\
5	-4.05614381088968\\
6	-3.89509963657743\\
7	-3.75326014715796\\
8	-3.61805451700271\\
9	-3.50595203926251\\
10	-3.39875322730864\\
20	-2.68924848284147\\
30	-2.2778608566559\\
40	-1.98016943259466\\
50	-1.76124907504374\\
60	-1.58324721833595\\
70	-1.43954907290119\\
80	-1.30750443621459\\
90	-1.20491516743236\\
};
\addplot [color=mycolor6,solid,line width=1pt,mark size=1pt,mark=*,mark options={solid},forget plot]
  table[row sep=crcr]{%
0.0001	-7.05974981800151\\
0.0002	-7.3388600176482\\
0.0003	-7.40834512313918\\
0.0004	-7.29940189393594\\
0.0005	-7.30703643300046\\
0.0006	-7.25598943363134\\
0.0007	-7.29752808930728\\
0.0008	-7.23439518004338\\
0.0009	-7.29824317092304\\
0.001	-7.28845734033017\\
0.002	-7.26231529908267\\
0.003	-7.25922499173715\\
0.004	-7.25660313891177\\
0.005	-7.25875374293954\\
0.006	-7.25326953297293\\
0.007	-7.24755687563657\\
0.008	-7.25147245689922\\
0.009	-7.24520760200886\\
0.01	-7.2494719470361\\
0.02	-7.2272855987276\\
0.03	-7.22002294253108\\
0.04	-7.22512509783374\\
0.05	-7.23227621307118\\
0.06	-7.23441722404481\\
0.07	-7.23885799559214\\
0.08	-7.2465162035955\\
0.09	-7.25034272946156\\
0.1	-7.25080356704841\\
0.2	-7.03596557435546\\
0.3	-6.65103355389655\\
0.4	-6.37478885328413\\
0.5	-6.1370517965316\\
0.6	-5.934808122419\\
0.7	-5.78073853564142\\
0.8	-5.66156108854349\\
0.9	-5.57373193601068\\
1	-5.51115891972217\\
2	-4.87316675233092\\
3	-4.51762663357828\\
4	-4.26816241898814\\
5	-4.07899909376666\\
6	-3.91537117701171\\
7	-3.77029360468545\\
8	-3.63702659580392\\
9	-3.5213909700998\\
10	-3.4164664219925\\
20	-2.70643275930623\\
30	-2.29416743557041\\
40	-1.99729941716123\\
50	-1.77954735362919\\
60	-1.59881545637116\\
70	-1.45258721487833\\
80	-1.32105969629424\\
90	-1.21792434698807\\
};
\addplot [color=mycolor7,solid,line width=1pt,mark size=1pt,mark=*,mark options={solid},forget plot]
  table[row sep=crcr]{%
0.0001	-7.46343038328998\\
0.0002	-7.25304726983996\\
0.0003	-7.27798478083927\\
0.0004	-7.45370873182598\\
0.0005	-7.44493507139566\\
0.0006	-7.48183006608554\\
0.0007	-7.43989990173846\\
0.0008	-7.45714075371729\\
0.0009	-7.46546726298521\\
0.001	-7.44783027617271\\
0.002	-7.45575933846942\\
0.003	-7.44923633175768\\
0.004	-7.44253076881983\\
0.005	-7.43533661646121\\
0.006	-7.43396619434307\\
0.007	-7.43602317579001\\
0.008	-7.42773484203066\\
0.009	-7.42974200132055\\
0.01	-7.42036827180906\\
0.02	-7.40575238875127\\
0.03	-7.40021348848766\\
0.04	-7.4060666658338\\
0.05	-7.41622443330294\\
0.06	-7.42006503634293\\
0.07	-7.425922635519\\
0.08	-7.43355266316461\\
0.09	-7.43201828655608\\
0.1	-7.41909944009777\\
0.2	-7.0658508706325\\
0.3	-6.66948437028781\\
0.4	-6.38499891053247\\
0.5	-6.14714708071315\\
0.6	-5.95823626304611\\
0.7	-5.81937588359273\\
0.8	-5.73543594069678\\
0.9	-5.64453190957056\\
1	-5.53620245142757\\
2	-4.91005051285758\\
3	-4.54482367722764\\
4	-4.29387264997783\\
5	-4.10010449100039\\
6	-3.93117213808528\\
7	-3.78551369333688\\
8	-3.65042481378501\\
9	-3.53755592525923\\
10	-3.42814018362221\\
20	-2.71887357642498\\
30	-2.30966386155784\\
40	-2.00884723077185\\
50	-1.79283789192649\\
60	-1.61369239162065\\
70	-1.46901872007019\\
80	-1.33520913597127\\
90	-1.23112660022218\\
};
\addplot [color=mycolor1,solid,line width=1pt,mark size=1pt,mark=*,mark options={solid},forget plot]
  table[row sep=crcr]{%
0.0001	-7.76437991995423\\
0.0002	-7.64436056625692\\
0.0003	-7.72098526799195\\
0.0004	-7.65983891432764\\
0.0005	-7.60906321775611\\
0.0006	-7.63848985954377\\
0.0007	-7.66478335936526\\
0.0008	-7.63107917041984\\
0.0009	-7.60768038572695\\
0.001	-7.62308787761336\\
0.002	-7.61643928175061\\
0.003	-7.61561612702795\\
0.004	-7.61920312167509\\
0.005	-7.60866430638739\\
0.006	-7.61091065754312\\
0.007	-7.60513593467111\\
0.008	-7.60391666151525\\
0.009	-7.59806775627416\\
0.01	-7.59380884099553\\
0.02	-7.57312732407311\\
0.03	-7.5641319128163\\
0.04	-7.57614055472551\\
0.05	-7.58959821953793\\
0.06	-7.5946838254147\\
0.07	-7.598768855561\\
0.08	-7.59754393685435\\
0.09	-7.57235416452322\\
0.1	-7.52964339911421\\
0.2	-7.07194958254724\\
0.3	-6.68355077921972\\
0.4	-6.39380025021631\\
0.5	-6.16620523153583\\
0.6	-5.99464685832228\\
0.7	-5.90168100807165\\
0.8	-5.79335686572726\\
0.9	-5.67110253535638\\
1	-5.56099713870138\\
2	-4.94165848492611\\
3	-4.5736427204612\\
4	-4.31966752763543\\
5	-4.11891531607361\\
6	-3.94912242420067\\
7	-3.80252323254957\\
8	-3.66585778280654\\
9	-3.55145188022819\\
10	-3.44613660139515\\
20	-2.73362934194901\\
30	-2.31915135365785\\
40	-2.02443777648785\\
50	-1.80710789820799\\
60	-1.6267495324942\\
70	-1.4810035765482\\
80	-1.34833796645576\\
90	-1.2442863624827\\
};
\addplot [color=mycolor2,solid,line width=1pt,mark size=1pt,mark=*,mark options={solid},forget plot]
  table[row sep=crcr]{%
0.0001	-7.8168832203616\\
0.0002	-7.78350398403727\\
0.0003	-7.74096454584873\\
0.0004	-7.74155263097568\\
0.0005	-7.76119993446241\\
0.0006	-7.77795960815108\\
0.0007	-7.79202030406955\\
0.0008	-7.77694116184809\\
0.0009	-7.77298394066291\\
0.001	-7.78883616138794\\
0.002	-7.78797745693631\\
0.003	-7.78981790128497\\
0.004	-7.77180339082383\\
0.005	-7.77654309000665\\
0.006	-7.76438586423497\\
0.007	-7.76709024052946\\
0.008	-7.76015036629007\\
0.009	-7.76002431675767\\
0.01	-7.75414839976988\\
0.02	-7.72735696271694\\
0.03	-7.7207720472863\\
0.04	-7.7351170600025\\
0.05	-7.75161097765763\\
0.06	-7.75635185054244\\
0.07	-7.74777990988765\\
0.08	-7.71677319916766\\
0.09	-7.66215844449684\\
0.1	-7.61038896763777\\
0.2	-7.07782397817734\\
0.3	-6.69180082537502\\
0.4	-6.41337675621358\\
0.5	-6.19615194102872\\
0.6	-6.06464470801087\\
0.7	-5.95233726095872\\
0.8	-5.81717781248031\\
0.9	-5.69383215624734\\
1	-5.59288936479844\\
2	-4.97495487034281\\
3	-4.60351498255177\\
4	-4.33863517436206\\
5	-4.1353062872111\\
6	-3.96252034258071\\
7	-3.81455710021308\\
8	-3.67957741488186\\
9	-3.56391530607078\\
10	-3.4559449862262\\
20	-2.74556484569752\\
30	-2.33194820420864\\
40	-2.03881672327776\\
50	-1.8200827739259\\
60	-1.63974088657258\\
70	-1.49267203675121\\
80	-1.35995522874492\\
90	-1.25485811475711\\
};
\addplot [color=mycolor3,solid,line width=1pt,mark size=1pt,mark=*,mark options={solid},forget plot]
  table[row sep=crcr]{%
0.0001	-8.00998678774122\\
0.0002	-7.98017462203433\\
0.0003	-7.95916014259446\\
0.0004	-7.97433487002114\\
0.0005	-7.94502144030251\\
0.0006	-7.95101592730716\\
0.0007	-7.96371923550726\\
0.0008	-7.94180520820687\\
0.0009	-7.95168888335\\
0.001	-7.94791019997019\\
0.002	-7.93017914147021\\
0.003	-7.93286775469908\\
0.004	-7.92669442203329\\
0.005	-7.91918561745871\\
0.006	-7.91949711294338\\
0.007	-7.91268378543668\\
0.008	-7.91106818223066\\
0.009	-7.90326936862823\\
0.01	-7.90277259240351\\
0.02	-7.87246352697034\\
0.03	-7.86692662731237\\
0.04	-7.88562324274374\\
0.05	-7.90415172056024\\
0.06	-7.89617494137367\\
0.07	-7.8562494633425\\
0.08	-7.79590609582136\\
0.09	-7.74169206253489\\
0.1	-7.71753374105384\\
0.2	-7.08268169191206\\
0.3	-6.69875714835208\\
0.4	-6.43453269069607\\
0.5	-6.25994822537482\\
0.6	-6.13027137261372\\
0.7	-5.96878841947309\\
0.8	-5.83474466520098\\
0.9	-5.71869415419588\\
1	-5.62907958873467\\
2	-5.00432411962671\\
3	-4.62540444257798\\
4	-4.35679799412046\\
5	-4.15003527324205\\
6	-3.9755291069514\\
7	-3.82742608365074\\
8	-3.69087431750336\\
9	-3.5759857971803\\
10	-3.47051866263562\\
20	-2.75812776330019\\
30	-2.34539139858735\\
40	-2.04915025082535\\
50	-1.83108770543454\\
60	-1.65093089093729\\
70	-1.50359354016559\\
80	-1.37240308554381\\
90	-1.2666989507657\\
};
\end{axis}
\end{tikzpicture}%

%% file: ksoddperiodic_lyaps_L100varyT.tex
\definecolor{mycolor1}{rgb}{0.00000,0.44700,0.74100}%
\definecolor{mycolor2}{rgb}{0.85000,0.32500,0.09800}%
\definecolor{mycolor3}{rgb}{0.92900,0.69400,0.12500}%
\definecolor{mycolor4}{rgb}{0.49400,0.18400,0.55600}%
\definecolor{mycolor5}{rgb}{0.46600,0.67400,0.18800}%
\definecolor{mycolor6}{rgb}{0.30100,0.74500,0.93300}%
\definecolor{mycolor7}{rgb}{0.63500,0.07800,0.18400}%
\begin{tikzpicture}

\begin{axis}[%
width=4.25in,
height=1.75in,
at={(0in,0in)},
scale only axis,
xmode=log,
xmin=0.0001,
xmax=90,
xminorticks=true,
xlabel={$T$},
xmajorgrids,
xminorgrids,
ymin=-7.61492394399101,
ymax=3.3133037149797,
ytick={-7.60090270954199,-6.90775627898064,-5.29834236561059,0,5.29834236561059},
yticklabels={{$-10$},{$-5$},{$-1$},{$0$},{$1$}},
ylabel={$\lambda_i$},
ylabel shift=-4pt,
ymajorgrids,
axis background/.style={fill=white}
]
\addplot [color=mycolor1,solid,line width=1pt,mark size=1pt,mark=*,mark options={solid},forget plot]
  table[row sep=crcr]{%
0.0001	-7.30979115997746\\
0.0002	-6.87222896274444\\
0.0003	-6.40964027464237\\
0.0004	-6.07575398904167\\
0.0005	-6.11075182611611\\
0.0006	-5.74463019658206\\
0.0007	-5.51491661560522\\
0.0008	-5.18365293384459\\
0.0009	-5.33981498921897\\
0.001	-4.9032441586012\\
0.002	-2.83106401804558\\
0.003	-3.50924463000519\\
0.004	2.62677183670991\\
0.005	2.33382868340479\\
0.006	1.31055624142729\\
0.007	2.94991221056772\\
0.008	3.02137964662277\\
0.009	3.3133037149797\\
0.01	2.36284329847607\\
0.02	2.5652449860379\\
0.03	2.94464902684002\\
0.04	3.04031244269703\\
0.05	2.99766093548118\\
0.06	3.04694624876684\\
0.07	2.65423055680603\\
0.08	3.05644529915846\\
0.09	2.9957991940192\\
0.1	3.05746248748009\\
0.2	2.87239408733591\\
0.3	2.92746293393087\\
0.4	2.99363919276992\\
0.5	2.94355292536954\\
0.6	2.85589213666245\\
0.7	2.82720251353539\\
0.8	2.85446064085293\\
0.9	2.76096688719861\\
1	2.81924475225017\\
2	2.83422733058486\\
3	2.89901138657069\\
4	2.90587931783761\\
5	2.89363806161737\\
6	2.8462460524927\\
7	2.83051589844714\\
8	2.89825813114518\\
9	2.88345203383343\\
10	2.88175152842861\\
20	2.88009481056115\\
30	2.87531628527161\\
40	2.8759983112064\\
50	2.87107234759699\\
60	2.87185054422357\\
70	2.87456644294178\\
80	2.8554775035721\\
90	2.87204148551474\\
};
\addplot [color=mycolor2,solid,line width=1pt,mark size=1pt,mark=*,mark options={solid},forget plot]
  table[row sep=crcr]{%
0.0001	-7.5984066582442\\
0.0002	-6.86676452214157\\
0.0003	-6.24988657167415\\
0.0004	-6.01425413954339\\
0.0005	-5.85915541765201\\
0.0006	-5.77237374841377\\
0.0007	-5.6970817943026\\
0.0008	-5.1927552565472\\
0.0009	-5.42690096205062\\
0.001	-5.18965059097638\\
0.002	-4.29108970822893\\
0.003	-3.25203842658597\\
0.004	0.427251328828875\\
0.005	-0.736909021858986\\
0.006	2.31049435155605\\
0.007	1.85936042545604\\
0.008	3.02118748806481\\
0.009	1.55572928687976\\
0.01	2.40826258927207\\
0.02	2.23449941367572\\
0.03	2.37664272124451\\
0.04	2.36923158324402\\
0.05	2.26247220404956\\
0.06	2.51874475095974\\
0.07	2.93286612467203\\
0.08	2.76915073437843\\
0.09	3.01039731909425\\
0.1	3.02386787061335\\
0.2	2.76492706485019\\
0.3	2.77347815378093\\
0.4	2.65288850995115\\
0.5	2.73887848871993\\
0.6	2.77840297810915\\
0.7	2.75039502993907\\
0.8	2.76022324893393\\
0.9	2.68529356654394\\
1	2.65792077729336\\
2	2.70108545705631\\
3	2.72913427991433\\
4	2.73421954223303\\
5	2.67619873120094\\
6	2.72085790218902\\
7	2.69012833479995\\
8	2.7224841889028\\
9	2.69010586376143\\
10	2.69485772124115\\
20	2.70349255100679\\
30	2.68921229200982\\
40	2.72359719881818\\
50	2.7028108245134\\
60	2.69377456074393\\
70	2.69191234027191\\
80	2.70763582841084\\
90	2.70483832427901\\
};
\addplot [color=mycolor3,solid,line width=1pt,mark size=1pt,mark=*,mark options={solid},forget plot]
  table[row sep=crcr]{%
0.0001	-7.49685229323925\\
0.0002	-6.91116039160097\\
0.0003	-6.35738371060259\\
0.0004	-6.04513409570324\\
0.0005	-5.90698799210904\\
0.0006	-5.80254601763004\\
0.0007	-5.60944623843816\\
0.0008	-5.30871931946303\\
0.0009	-5.28452147666204\\
0.001	-5.25130970019869\\
0.002	-3.99975743579102\\
0.003	-3.60837922607618\\
0.004	-0.110788184771775\\
0.005	-1.85643577350856\\
0.006	2.44617600888881\\
0.007	-1.07272908780139\\
0.008	1.00656615035124\\
0.009	2.15808614719702\\
0.01	-0.531979038242177\\
0.02	2.00426732737502\\
0.03	2.644644831631\\
0.04	2.58200977195292\\
0.05	2.15498966251753\\
0.06	2.6548106596413\\
0.07	2.4292883586614\\
0.08	2.78437756763085\\
0.09	2.64564039199345\\
0.1	2.53282035617089\\
0.2	2.796102599312\\
0.3	2.53183160331858\\
0.4	2.54361028908274\\
0.5	2.65204505669471\\
0.6	2.63340466025786\\
0.7	2.48743656134718\\
0.8	2.49099113769959\\
0.9	2.56105735075935\\
1	2.56000478169289\\
2	2.47161404706404\\
3	2.54170888079167\\
4	2.53930164071301\\
5	2.5211455251254\\
6	2.54829343225882\\
7	2.54817998849241\\
8	2.53747701158808\\
9	2.52512811655941\\
10	2.55545548312019\\
20	2.52099729349205\\
30	2.52318811594994\\
40	2.53939092486755\\
50	2.53925225811195\\
60	2.5389218911969\\
70	2.5411926182135\\
80	2.5472103000342\\
90	2.56030543208955\\
};
\addplot [color=mycolor4,solid,line width=1pt,mark size=1pt,mark=*,mark options={solid},forget plot]
  table[row sep=crcr]{%
0.0001	-7.15201319554336\\
0.0002	-6.73263201388244\\
0.0003	-6.32329562666913\\
0.0004	-6.25112779842566\\
0.0005	-5.87788103933455\\
0.0006	-5.75327414628934\\
0.0007	-5.50093098545523\\
0.0008	-5.63523553885479\\
0.0009	-5.28162188019468\\
0.001	-5.26783624517443\\
0.002	-4.39152359358555\\
0.003	-3.49777568513852\\
0.004	-2.46106467859742\\
0.005	2.49516788958609\\
0.006	0.396929227591565\\
0.007	2.98169950239369\\
0.008	2.87860871694845\\
0.009	2.95545919874611\\
0.01	2.21692704418559\\
0.02	1.89871651645066\\
0.03	2.46843770385107\\
0.04	2.4596729660995\\
0.05	2.18259675163655\\
0.06	2.21251831144404\\
0.07	2.64800891758992\\
0.08	2.36498745407703\\
0.09	2.81360907108177\\
0.1	2.41091521091664\\
0.2	2.61945229701427\\
0.3	2.47448673381001\\
0.4	2.4810180003099\\
0.5	2.45958346953624\\
0.6	2.46808899376715\\
0.7	2.52740112599397\\
0.8	2.32652281112862\\
0.9	2.41035982848221\\
1	2.37143624153504\\
2	2.32149403393018\\
3	2.3778458892167\\
4	2.35236366450055\\
5	2.34627229097349\\
6	2.39794993690253\\
7	2.39734332791992\\
8	2.40039645089782\\
9	2.36436984365003\\
10	2.37467489385379\\
20	2.36166422164031\\
30	2.37449258072048\\
40	2.36158166088149\\
50	2.36880477428132\\
60	2.36743325143868\\
70	2.38400443170078\\
80	2.38997929496371\\
90	2.39974977541688\\
};
\addplot [color=mycolor5,solid,line width=1pt,mark size=1pt,mark=*,mark options={solid},forget plot]
  table[row sep=crcr]{%
0.0001	-7.57031387025008\\
0.0002	-7.0697235477921\\
0.0003	-6.43490575531354\\
0.0004	-6.05912043028606\\
0.0005	-5.90721636646481\\
0.0006	-5.69467944691848\\
0.0007	-5.46119548754209\\
0.0008	-5.4313971064647\\
0.0009	-5.18655751441843\\
0.001	-5.26148728738529\\
0.002	-4.6090533890191\\
0.003	-2.47448611507123\\
0.004	-1.96470974869044\\
0.005	-1.65366166120545\\
0.006	-2.92162683464388\\
0.007	1.36895696819824\\
0.008	-2.58722393376203\\
0.009	2.69848035041732\\
0.01	1.91405866152069\\
0.02	1.17656034030425\\
0.03	2.25562436627333\\
0.04	2.27852487342229\\
0.05	2.06053175000566\\
0.06	1.76304075848825\\
0.07	2.53379691491875\\
0.08	2.11704083858765\\
0.09	2.33820424890569\\
0.1	2.40936218469142\\
0.2	2.29602557709864\\
0.3	2.30011696146048\\
0.4	2.26014332173243\\
0.5	2.29322628303752\\
0.6	2.25362740489925\\
0.7	2.21276685055694\\
0.8	2.17726766971515\\
0.9	2.14110071172011\\
1	2.04789843929786\\
2	2.13434550302124\\
3	2.21404651840414\\
4	2.19097529672687\\
5	2.14791079280059\\
6	2.23734523817489\\
7	2.21157340332487\\
8	2.21842540522052\\
9	2.18328918777376\\
10	2.20700943862314\\
20	2.18689864580843\\
30	2.19640803023222\\
40	2.19163040659826\\
50	2.18763848322425\\
60	2.19554104189136\\
70	2.2195453709213\\
80	2.22437658166221\\
90	2.24382568124346\\
};
\addplot [color=mycolor6,solid,line width=1pt,mark size=1pt,mark=*,mark options={solid},forget plot]
  table[row sep=crcr]{%
0.0001	-7.55443486260963\\
0.0002	-6.76858386193183\\
0.0003	-6.51173912738674\\
0.0004	-6.37364858521686\\
0.0005	-6.03293647063624\\
0.0006	-5.86041034107577\\
0.0007	-5.80409055381228\\
0.0008	-5.65005505289558\\
0.0009	-5.24390515865956\\
0.001	-5.43775581848752\\
0.002	-4.40189178618887\\
0.003	-3.49913318937108\\
0.004	-2.47075005299343\\
0.005	1.88906616243877\\
0.006	-1.32635590968871\\
0.007	2.54392838687134\\
0.008	0.0280100482493259\\
0.009	0.0482938851779092\\
0.01	1.5078347768765\\
0.02	2.06247260716863\\
0.03	2.56680276492925\\
0.04	2.34793486161603\\
0.05	1.94159664337891\\
0.06	2.39930173063437\\
0.07	2.24030446846759\\
0.08	2.09880213531003\\
0.09	1.44505822919585\\
0.1	1.9830722752055\\
0.2	2.2701991243278\\
0.3	2.1318560268659\\
0.4	1.94185603896179\\
0.5	2.07524597107477\\
0.6	2.08326251859307\\
0.7	1.97718368445142\\
0.8	2.17253372810689\\
0.9	1.92924306564071\\
1	2.09148930809278\\
2	2.09850993289279\\
3	1.99216820709143\\
4	2.00408619639872\\
5	1.93309312233727\\
6	2.03538651177816\\
7	2.00856214225003\\
8	2.0293168372951\\
9	1.98463591035588\\
10	2.04346425580051\\
20	2.01059326782839\\
30	2.01811118402786\\
40	2.00775038623065\\
50	2.01318291246314\\
60	2.0249055883846\\
70	2.02792940047256\\
80	2.0492823311729\\
90	2.08688255422729\\
};
\addplot [color=mycolor7,solid,line width=1pt,mark size=1pt,mark=*,mark options={solid},forget plot]
  table[row sep=crcr]{%
0.0001	-7.4038981987984\\
0.0002	-6.73898146826291\\
0.0003	-6.50593453392368\\
0.0004	-6.22634462844594\\
0.0005	-6.08888751143439\\
0.0006	-5.74093294303106\\
0.0007	-5.61461307706679\\
0.0008	-5.28214584711926\\
0.0009	-5.44612344614859\\
0.001	-4.91591741811314\\
0.002	-4.81847159897299\\
0.003	-4.0394731678484\\
0.004	-2.96226123994459\\
0.005	1.00524909636864\\
0.006	-2.02046660595086\\
0.007	-2.80421213087505\\
0.008	-1.60708085042716\\
0.009	-1.53651669747584\\
0.01	2.18822475623111\\
0.02	1.59165324878094\\
0.03	1.99787121126716\\
0.04	2.24374005974881\\
0.05	0.506225743767803\\
0.06	1.80332640531445\\
0.07	1.39830864935632\\
0.08	1.71560150818688\\
0.09	2.14294406181587\\
0.1	1.58334496704447\\
0.2	1.56582041908805\\
0.3	1.84410338626762\\
0.4	1.86024547341198\\
0.5	1.83880089556551\\
0.6	1.91215647786185\\
0.7	1.8233050870113\\
0.8	1.88013351831615\\
0.9	1.69200889499757\\
1	1.70207747434143\\
2	1.81851153689738\\
3	1.78735908180177\\
4	1.7358999617272\\
5	1.82379639629382\\
6	1.89923608251505\\
7	1.80650541802288\\
8	1.86109865098213\\
9	1.81989325141276\\
10	1.73394276728782\\
20	1.80395984512884\\
30	1.82058612825748\\
40	1.82721840337815\\
50	1.82138129537561\\
60	1.82617404849406\\
70	1.84269751412261\\
80	1.86020839603339\\
90	1.91506991631637\\
};
\addplot [color=mycolor1,solid,line width=1pt,mark size=1pt,mark=*,mark options={solid},forget plot]
  table[row sep=crcr]{%
0.0001	-7.51056445396619\\
0.0002	-6.88366861475654\\
0.0003	-6.50193369814733\\
0.0004	-6.08257497093426\\
0.0005	-6.02530783937132\\
0.0006	-5.87234845653653\\
0.0007	-5.47789707669145\\
0.0008	-5.55876099372587\\
0.0009	-5.55702549487877\\
0.001	-5.16986962061029\\
0.002	-4.46650685171086\\
0.003	-3.55853937507039\\
0.004	-4.29609187370349\\
0.005	-2.94940760042276\\
0.006	-0.238263259508562\\
0.007	-3.02949649109405\\
0.008	-0.261758276887022\\
0.009	-3.29057632712564\\
0.01	0.805126967468149\\
0.02	1.71108531549135\\
0.03	2.16990511508832\\
0.04	1.69652801804989\\
0.05	1.6675493659074\\
0.06	1.70841420671401\\
0.07	0.954455613965296\\
0.08	1.7127554585244\\
0.09	1.66292949959793\\
0.1	1.69870249275361\\
0.2	1.90390135469881\\
0.3	1.50851012472205\\
0.4	1.52043750837185\\
0.5	1.74944419955148\\
0.6	1.81892221139387\\
0.7	1.67099053507961\\
0.8	1.50593212934459\\
0.9	1.62645126346037\\
1	1.33354455066228\\
2	1.35595602475644\\
3	1.6106729389091\\
4	1.48398341735774\\
5	1.50639814579519\\
6	1.62161657826712\\
7	1.59650223438002\\
8	1.55711551614261\\
9	1.52111735407223\\
10	1.5070057766677\\
20	1.52795952471514\\
30	1.58561301917789\\
40	1.58505841712361\\
50	1.57464735185546\\
60	1.61802805707821\\
70	1.62062788652467\\
80	1.64351384141504\\
90	1.72223658218148\\
};
\addplot [color=mycolor2,solid,line width=1pt,mark size=1pt,mark=*,mark options={solid},forget plot]
  table[row sep=crcr]{%
0.0001	-7.46717308820205\\
0.0002	-6.88436814314039\\
0.0003	-6.42939944847947\\
0.0004	-6.22360633674194\\
0.0005	-6.16224001091897\\
0.0006	-5.92619773045292\\
0.0007	-5.8213178461556\\
0.0008	-5.73303895684923\\
0.0009	-5.35596428842934\\
0.001	-5.25735605533924\\
0.002	-4.80543475924385\\
0.003	-3.98401038082981\\
0.004	-3.95924427335638\\
0.005	-4.19290219653201\\
0.006	-3.77080218371622\\
0.007	-1.37869696921865\\
0.008	-1.78229063119967\\
0.009	-2.87192500717486\\
0.01	-2.54527042455228\\
0.02	-1.55946455046628\\
0.03	-0.262802010881182\\
0.04	1.14756011354709\\
0.05	2.19804205613706\\
0.06	1.3718657432509\\
0.07	1.68985356223572\\
0.08	1.48060941089171\\
0.09	1.28815440610052\\
0.1	1.74858654313617\\
0.2	1.37965671345109\\
0.3	1.51031455004274\\
0.4	1.32079715114861\\
0.5	1.48884750755163\\
0.6	1.15453575492391\\
0.7	0.562074628575378\\
0.8	1.44229392434526\\
0.9	0.941235392544528\\
1	1.17622815682761\\
2	1.23961327518074\\
3	1.32245985737465\\
4	1.24031561980903\\
5	1.23268106223784\\
6	1.29458541554028\\
7	1.3180217029065\\
8	1.29696203651544\\
9	1.28055856451113\\
10	1.21204257881866\\
20	1.24767668699022\\
30	1.32451707815775\\
40	1.29902042828411\\
50	1.30157432356089\\
60	1.32298184352428\\
70	1.39097114157918\\
80	1.45201836309899\\
90	1.55139171893581\\
};
\addplot [color=mycolor3,solid,line width=1pt,mark size=1pt,mark=*,mark options={solid},forget plot]
  table[row sep=crcr]{%
0.0001	-7.50603787381121\\
0.0002	-7.05468533273805\\
0.0003	-6.71907434533379\\
0.0004	-6.22897955546562\\
0.0005	-5.69084148930989\\
0.0006	-5.95487051148183\\
0.0007	-5.70778817447663\\
0.0008	-5.57748452252611\\
0.0009	-5.57537890288294\\
0.001	-5.36886156479134\\
0.002	-4.37287484604665\\
0.003	-4.06547945740341\\
0.004	-3.16202952945283\\
0.005	-2.99603169029452\\
0.006	-3.31447763607601\\
0.007	-2.35788346875274\\
0.008	1.81113483558894\\
0.009	-3.09534136161266\\
0.01	-2.62731090606175\\
0.02	-1.9074667500357\\
0.03	-0.857185924561108\\
0.04	-0.17366359144407\\
0.05	-0.616141089316327\\
0.06	1.46959538588661\\
0.07	1.51654881670716\\
0.08	0.119020942750423\\
0.09	1.23714233190246\\
0.1	0.974125689615075\\
0.2	1.07369747435135\\
0.3	1.27385005959256\\
0.4	0.67286480251002\\
0.5	1.04046671414998\\
0.6	0.627957133616862\\
0.7	0.789103847569959\\
0.8	0.988097965084735\\
0.9	0.820163107267514\\
1	0.462041075739796\\
2	0.896497970942843\\
3	0.846134216488683\\
4	0.934561875210303\\
5	0.895494743563829\\
6	0.876571303217667\\
7	1.0249039997698\\
8	0.844183100341569\\
9	1.01386086203697\\
10	0.896025208438206\\
20	0.937890196223775\\
30	0.995469237485056\\
40	0.904649583159629\\
50	0.950996441082538\\
60	0.981074863161044\\
70	1.02789793667896\\
80	1.17921396765122\\
90	1.36113395128298\\
};
\addplot [color=mycolor4,solid,line width=1pt,mark size=1pt,mark=*,mark options={solid},forget plot]
  table[row sep=crcr]{%
0.0001	-7.60034030170155\\
0.0002	-6.8308604115793\\
0.0003	-6.50494630425806\\
0.0004	-6.2795207263831\\
0.0005	-6.20049430258877\\
0.0006	-5.96513280655932\\
0.0007	-5.95173560426437\\
0.0008	-5.22389212889252\\
0.0009	-5.4375258089585\\
0.001	-5.68976272359597\\
0.002	-4.94467474354872\\
0.003	-3.7254291213519\\
0.004	-4.09138274876976\\
0.005	-3.88057046037717\\
0.006	-3.32995393502671\\
0.007	-2.92771804590685\\
0.008	-2.99809707021016\\
0.009	-1.48803318683821\\
0.01	-3.08336911916747\\
0.02	-2.47847936108808\\
0.03	-0.380450389113486\\
0.04	-1.78465110045367\\
0.05	-1.65873423491518\\
0.06	-1.14588279559448\\
0.07	-1.37232043035531\\
0.08	1.64596998264106\\
0.09	-1.19891141671302\\
0.1	0.397712846938935\\
0.2	0.166367480185329\\
0.3	-0.0700622336554333\\
0.4	0.708238379151897\\
0.5	0.19960133145413\\
0.6	-0.0532046460769977\\
0.7	0.288812728437246\\
0.8	0.619977760778353\\
0.9	0.337911097694208\\
1	0.2810885265267\\
2	0.342276171211733\\
3	0.340418445689448\\
4	0.357717716044846\\
5	0.346851469812553\\
6	0.264581282081306\\
7	0.381617839392731\\
8	0.462292141483782\\
9	0.515402672628952\\
10	0.218137582734903\\
20	0.405750953929688\\
30	0.508087807623469\\
40	0.471939391296012\\
50	0.479235726993231\\
60	0.566268783602068\\
70	0.692874654282147\\
80	0.90455339621599\\
90	1.11281897741404\\
};
\addplot [color=mycolor5,solid,line width=1pt,mark size=1pt,mark=*,mark options={solid},forget plot]
  table[row sep=crcr]{%
0.0001	-7.55810614330082\\
0.0002	-6.94139447026731\\
0.0003	-6.57351478579857\\
0.0004	-6.231216702915\\
0.0005	-6.10011893785034\\
0.0006	-5.78745932700821\\
0.0007	-5.81943604246883\\
0.0008	-5.69586084739848\\
0.0009	-5.64947495561268\\
0.001	-5.50659090488894\\
0.002	-4.57250762525682\\
0.003	-4.46440480010385\\
0.004	-4.07121285338591\\
0.005	-3.7524569831952\\
0.006	-4.15077678743782\\
0.007	-2.66613934323101\\
0.008	-3.85298901269277\\
0.009	-2.91961411866867\\
0.01	-2.32016865339599\\
0.02	-1.80050312243104\\
0.03	-1.51960917027188\\
0.04	-2.39255020489606\\
0.05	-1.15833337473001\\
0.06	-1.42317238564731\\
0.07	-1.07735102722541\\
0.08	-1.48850758202937\\
0.09	-0.941081995711605\\
0.1	-0.332036601183503\\
0.2	-0.0171538987139965\\
0.3	-0.800955778801072\\
0.4	0.162991828157015\\
0.5	-0.410370863337986\\
0.6	-0.30585465136876\\
0.7	0.132523026547472\\
0.8	0.127673172403902\\
0.9	-0.0419291423111211\\
1	0.0237582888509919\\
2	-0.00438145708139328\\
3	-0.14945015151876\\
4	0.0654016444665359\\
5	0.126017668789965\\
6	-0.0252969798339756\\
7	-0.00161526509760691\\
8	0.134192459500527\\
9	0.0138287462408628\\
10	0.0277375491159408\\
20	0.00226050147485589\\
30	0.0155800046862626\\
40	0.0450998166563539\\
50	0.0787247007123751\\
60	0.164449231932112\\
70	0.274291599977409\\
80	0.509574102609075\\
90	0.851892488481817\\
};
\addplot [color=mycolor6,solid,line width=1pt,mark size=1pt,mark=*,mark options={solid},forget plot]
  table[row sep=crcr]{%
0.0001	-7.52422637251959\\
0.0002	-6.85112439472529\\
0.0003	-6.75346369053064\\
0.0004	-6.43894411941084\\
0.0005	-6.20852208054724\\
0.0006	-5.61391876727552\\
0.0007	-5.6139102339563\\
0.0008	-5.85154753961766\\
0.0009	-5.73213453015143\\
0.001	-5.36568332857228\\
0.002	-4.81788359501071\\
0.003	-4.14361748849128\\
0.004	-3.8750230281016\\
0.005	-3.0227395568171\\
0.006	-3.84058571340544\\
0.007	-3.70112368070617\\
0.008	-3.68517928380956\\
0.009	-3.2994450760018\\
0.01	-3.60761474560652\\
0.02	-2.22121365521622\\
0.03	-1.77215397977675\\
0.04	-1.83128125386309\\
0.05	-2.0668940897026\\
0.06	-1.64579097084043\\
0.07	-1.17139381984602\\
0.08	-1.42017171674748\\
0.09	-1.11526155764207\\
0.1	-1.99093259883582\\
0.2	-0.588195201935375\\
0.3	-0.465935060820963\\
0.4	-0.73144969186059\\
0.5	-0.39562106209417\\
0.6	-0.810808317629763\\
0.7	-0.706820505675791\\
0.8	-0.597269974198831\\
0.9	-0.78531810636744\\
1	-0.747742718645886\\
2	-0.507176648914928\\
3	-0.642839109284726\\
4	-0.500976357402223\\
5	-0.41698243304918\\
6	-0.333174886240184\\
7	-0.357779988420897\\
8	-0.351636592648094\\
9	-0.328263358101774\\
10	-0.328507619130204\\
20	-0.313973578937227\\
30	-0.382631032098714\\
40	-0.348123939670053\\
50	-0.373502560678884\\
60	-0.276675198954341\\
70	-0.164052377218772\\
80	0.122098776755515\\
90	0.539130974209235\\
};
\addplot [color=mycolor7,solid,line width=1pt,mark size=1pt,mark=*,mark options={solid},forget plot]
  table[row sep=crcr]{%
0.0001	-7.57922333170806\\
0.0002	-6.85946461898915\\
0.0003	-6.51724515220137\\
0.0004	-6.3992182724272\\
0.0005	-5.99364345795332\\
0.0006	-5.98324813622613\\
0.0007	-5.83384356425602\\
0.0008	-5.62853089342488\\
0.0009	-5.42881148767842\\
0.001	-5.51091576224728\\
0.002	-4.89219772610981\\
0.003	-4.65120583691101\\
0.004	-4.24770434390471\\
0.005	-4.27801881640448\\
0.006	-3.68173492395194\\
0.007	-3.46354583493273\\
0.008	-3.62007319154389\\
0.009	-3.82261926592123\\
0.01	-3.36870126704572\\
0.02	-2.80867159817355\\
0.03	-2.90840233855042\\
0.04	-1.86281144724347\\
0.05	-2.31767090192339\\
0.06	-2.16644771681415\\
0.07	-2.33067172600486\\
0.08	-2.18927281928695\\
0.09	-2.17239432743661\\
0.1	-1.98054990935892\\
0.2	-0.897387579461245\\
0.3	-1.44909956792608\\
0.4	-1.07052432790585\\
0.5	-1.15875129520733\\
0.6	-1.31713952885422\\
0.7	-1.30923927198305\\
0.8	-1.09082825065428\\
0.9	-1.27105605496158\\
1	-1.07909487114703\\
2	-0.705948873221121\\
3	-0.915062806944866\\
4	-0.982387198301549\\
5	-0.934287182834026\\
6	-1.01959792491138\\
7	-0.924239152796276\\
8	-0.941638476166108\\
9	-0.928568327450405\\
10	-1.00030818256007\\
20	-0.93444713184774\\
30	-0.903233510443107\\
40	-0.909135774537774\\
50	-0.903879623095779\\
60	-0.769655927632182\\
70	-0.611837348365667\\
80	-0.318063124272517\\
90	0.170033615624746\\
};
\addplot [color=mycolor1,solid,line width=1pt,mark size=1pt,mark=*,mark options={solid},forget plot]
  table[row sep=crcr]{%
0.0001	-7.53496498383303\\
0.0002	-6.84832157405521\\
0.0003	-6.63258563426195\\
0.0004	-6.38123224812114\\
0.0005	-6.08178001066606\\
0.0006	-5.69850159191734\\
0.0007	-5.85682545638227\\
0.0008	-5.76638872717765\\
0.0009	-5.37789949514449\\
0.001	-5.42859301803792\\
0.002	-4.89178073486899\\
0.003	-4.51479699526931\\
0.004	-4.27396205655637\\
0.005	-3.96439048909363\\
0.006	-3.76224337514069\\
0.007	-4.1652389979812\\
0.008	-3.60707836175424\\
0.009	-3.35222260662238\\
0.01	-3.19633805653915\\
0.02	-3.32857916072524\\
0.03	-2.80057487555905\\
0.04	-2.742600707623\\
0.05	-2.85779755653751\\
0.06	-2.28236907685113\\
0.07	-2.3211960804476\\
0.08	-2.00659670567997\\
0.09	-1.82253575240489\\
0.1	-2.25248831243408\\
0.2	-1.69687661280716\\
0.3	-1.68995574834876\\
0.4	-1.50462302660909\\
0.5	-1.60200784720278\\
0.6	-1.48493345178876\\
0.7	-1.60628232178168\\
0.8	-1.62782691527803\\
0.9	-1.74595563655044\\
1	-1.34671285302628\\
2	-1.47628221418754\\
3	-1.34927713104736\\
4	-1.4808033321215\\
5	-1.44149610674366\\
6	-1.45707980775805\\
7	-1.39345927501787\\
8	-1.38931850997361\\
9	-1.34510260466339\\
10	-1.38472145852298\\
20	-1.3960160876687\\
30	-1.40723078313133\\
40	-1.36956450569646\\
50	-1.35905675873987\\
60	-1.23011332133939\\
70	-1.03557052211286\\
80	-0.699318970046368\\
90	-0.195443811700215\\
};
\addplot [color=mycolor2,solid,line width=1pt,mark size=1pt,mark=*,mark options={solid},forget plot]
  table[row sep=crcr]{%
0.0001	-7.55997262633739\\
0.0002	-6.98168483047574\\
0.0003	-6.5925033337536\\
0.0004	-6.45331470712233\\
0.0005	-6.34200669517828\\
0.0006	-5.86997228657588\\
0.0007	-5.86201568247748\\
0.0008	-5.76397625420474\\
0.0009	-5.60098531605612\\
0.001	-5.29634087403412\\
0.002	-4.9216871688298\\
0.003	-4.50238848116603\\
0.004	-3.95060699702833\\
0.005	-3.92410614520403\\
0.006	-4.49107746155521\\
0.007	-4.14971529521479\\
0.008	-3.76510113100133\\
0.009	-3.71915955789414\\
0.01	-3.67555779527513\\
0.02	-3.29942473269076\\
0.03	-3.00627202373072\\
0.04	-3.0961634264438\\
0.05	-2.86747508879454\\
0.06	-2.63270287606747\\
0.07	-2.64885981463743\\
0.08	-2.2795000678947\\
0.09	-2.3543336915303\\
0.1	-2.37815835493134\\
0.2	-2.08084582862294\\
0.3	-1.97279448900717\\
0.4	-1.73545424819962\\
0.5	-1.91357110784618\\
0.6	-2.16925233384871\\
0.7	-1.99385105691579\\
0.8	-1.9162271833746\\
0.9	-1.78810930530799\\
1	-1.75798730342937\\
2	-1.81285911437415\\
3	-1.84714564423348\\
4	-1.83405843006717\\
5	-1.8608788991405\\
6	-1.85456046981021\\
7	-1.77158298066503\\
8	-1.91191235746258\\
9	-1.77673280266558\\
10	-1.90490367016372\\
20	-1.87836977297369\\
30	-1.79751229253218\\
40	-1.79914677985861\\
50	-1.73973232482434\\
60	-1.59003016661846\\
70	-1.37687706823094\\
80	-1.078162937725\\
90	-0.563031856563012\\
};
\addplot [color=mycolor3,solid,line width=1pt,mark size=1pt,mark=*,mark options={solid},forget plot]
  table[row sep=crcr]{%
0.0001	-7.51787695202007\\
0.0002	-7.19173165841521\\
0.0003	-6.62637801598776\\
0.0004	-6.4102390651552\\
0.0005	-6.36618573143441\\
0.0006	-6.13600203265602\\
0.0007	-5.85383071924346\\
0.0008	-5.78241409881316\\
0.0009	-5.47240277051639\\
0.001	-5.60795089040046\\
0.002	-4.66049033401383\\
0.003	-4.13945443130887\\
0.004	-4.23180917335367\\
0.005	-4.09725459457826\\
0.006	-3.87520726756371\\
0.007	-4.08537261560559\\
0.008	-4.19452305333495\\
0.009	-4.16415977412106\\
0.01	-3.84073964941514\\
0.02	-3.19348444123882\\
0.03	-3.11385052184709\\
0.04	-3.06841736549967\\
0.05	-3.27542660652754\\
0.06	-2.82034280895913\\
0.07	-2.77439197031961\\
0.08	-2.68566583576687\\
0.09	-2.72106935942492\\
0.1	-2.92961083745751\\
0.2	-2.52519330413364\\
0.3	-2.19614995324205\\
0.4	-2.14323544106839\\
0.5	-2.15296456665266\\
0.6	-2.12849855776502\\
0.7	-2.40702582743796\\
0.8	-2.27800046561469\\
0.9	-2.20333221577188\\
1	-2.15751267581458\\
2	-2.13010019573653\\
3	-2.2086966019003\\
4	-2.16740881185199\\
5	-2.18322076185149\\
6	-2.15409356457006\\
7	-2.21629054314472\\
8	-2.19744291505577\\
9	-2.21118796775248\\
10	-2.21843630648382\\
20	-2.18302999075584\\
30	-2.14764974095517\\
40	-2.16838061269518\\
50	-2.07329554557706\\
60	-1.92740605564831\\
70	-1.68133985009403\\
80	-1.35606379678611\\
90	-0.885159521967274\\
};
\addplot [color=mycolor4,solid,line width=1pt,mark size=1pt,mark=*,mark options={solid},forget plot]
  table[row sep=crcr]{%
0.0001	-7.60154030591542\\
0.0002	-6.83852460602629\\
0.0003	-6.78972776937888\\
0.0004	-6.54802950446943\\
0.0005	-6.15832785466482\\
0.0006	-5.9572523446261\\
0.0007	-5.84599570424075\\
0.0008	-5.85052923693531\\
0.0009	-5.9307752386949\\
0.001	-5.30504918916926\\
0.002	-4.70419313057316\\
0.003	-4.65494579234491\\
0.004	-4.47348357650576\\
0.005	-4.41751715589662\\
0.006	-4.04074746390194\\
0.007	-4.06983600259493\\
0.008	-4.07002767052101\\
0.009	-3.89392663377963\\
0.01	-3.85336544529847\\
0.02	-3.59427523137165\\
0.03	-3.35925924361185\\
0.04	-3.11722017497364\\
0.05	-2.96394908089531\\
0.06	-2.8849966587966\\
0.07	-2.81929023431105\\
0.08	-2.80230474950233\\
0.09	-2.97725517055747\\
0.1	-2.69994595322413\\
0.2	-2.48925687399563\\
0.3	-2.82379271497732\\
0.4	-2.73161890373001\\
0.5	-2.520310554815\\
0.6	-2.47611585113074\\
0.7	-2.61493583751599\\
0.8	-2.60794316160016\\
0.9	-2.62330092849552\\
1	-2.5308905555776\\
2	-2.51186223162068\\
3	-2.47185270771152\\
4	-2.58963154498356\\
5	-2.51992601435286\\
6	-2.57738426622256\\
7	-2.48146046939385\\
8	-2.54029300286552\\
9	-2.53919195410923\\
10	-2.54187623668218\\
20	-2.54529712179275\\
30	-2.49244960655429\\
40	-2.4782869715724\\
50	-2.32615552670666\\
60	-2.14500450388872\\
70	-1.9229142348829\\
80	-1.58963673749559\\
90	-1.14104431570964\\
};
\addplot [color=mycolor5,solid,line width=1pt,mark size=1pt,mark=*,mark options={solid},forget plot]
  table[row sep=crcr]{%
0.0001	-7.61492394399101\\
0.0002	-6.96110128576263\\
0.0003	-6.7529283459124\\
0.0004	-6.36408792502712\\
0.0005	-6.15339780590459\\
0.0006	-5.82136721969296\\
0.0007	-5.89315086154079\\
0.0008	-5.75882547700721\\
0.0009	-5.5754971434186\\
0.001	-5.70145711006898\\
0.002	-5.10836015614247\\
0.003	-4.95937146830654\\
0.004	-4.57044466758186\\
0.005	-4.69020729422613\\
0.006	-4.12779968473556\\
0.007	-4.42525753099467\\
0.008	-4.16660718763364\\
0.009	-4.00162466634837\\
0.01	-4.12079857462302\\
0.02	-3.46856656976466\\
0.03	-3.10682778253128\\
0.04	-3.35001964778242\\
0.05	-3.12680911016342\\
0.06	-3.09617579200824\\
0.07	-3.14910501498188\\
0.08	-3.20310028379567\\
0.09	-2.91679333092804\\
0.1	-2.86709336005096\\
0.2	-2.9583349016805\\
0.3	-2.81199862335243\\
0.4	-2.96416474866838\\
0.5	-2.73505075087643\\
0.6	-2.91582629483793\\
0.7	-3.00067243092038\\
0.8	-2.85635357686754\\
0.9	-2.85613195977486\\
1	-2.82985314729559\\
2	-2.81574723532561\\
3	-2.84419606097313\\
4	-2.81608771012551\\
5	-2.84164035483293\\
6	-2.86110710542045\\
7	-2.83587131357454\\
8	-2.85946675270908\\
9	-2.85390499241705\\
10	-2.8696253298842\\
20	-2.81638908682841\\
30	-2.78048996987192\\
40	-2.73228031576696\\
50	-2.56419829645107\\
60	-2.36319958771936\\
70	-2.10017341772531\\
80	-1.76473832883107\\
90	-1.33726310980508\\
};
\addplot [color=mycolor6,solid,line width=1pt,mark size=1pt,mark=*,mark options={solid},forget plot]
  table[row sep=crcr]{%
0.0001	-7.6062929533455\\
0.0002	-7.10152733112132\\
0.0003	-6.76832856546824\\
0.0004	-6.22966132814215\\
0.0005	-6.03225993342751\\
0.0006	-6.10929020880586\\
0.0007	-5.94256335066482\\
0.0008	-5.7855398089073\\
0.0009	-5.63347150678734\\
0.001	-5.564878695793\\
0.002	-4.81856740656262\\
0.003	-4.74723666770425\\
0.004	-4.55812356390735\\
0.005	-4.19663393882491\\
0.006	-4.32801488168413\\
0.007	-4.07012180075624\\
0.008	-4.37248906871928\\
0.009	-4.01991975741844\\
0.01	-3.91956060252893\\
0.02	-3.58747681942176\\
0.03	-3.59516793774264\\
0.04	-3.5706257957146\\
0.05	-3.14529130130956\\
0.06	-3.49013044943655\\
0.07	-3.43364395651559\\
0.08	-3.28461423950723\\
0.09	-3.2634246704061\\
0.1	-3.33287294461714\\
0.2	-3.09341123652569\\
0.3	-3.12909322791203\\
0.4	-3.15809822269278\\
0.5	-3.24798654058927\\
0.6	-3.180886182904\\
0.7	-3.13822196288509\\
0.8	-3.02791522535252\\
0.9	-3.1437474941467\\
1	-3.05255061520539\\
2	-3.12323856890331\\
3	-3.12289577405512\\
4	-3.09439927738566\\
5	-3.1145365244932\\
6	-3.11926962936579\\
7	-3.07671017626472\\
8	-3.1018158241425\\
9	-3.13209628637882\\
10	-3.09624086666636\\
20	-3.09788393470757\\
30	-3.03245069397934\\
40	-2.92143241384266\\
50	-2.74348428971241\\
60	-2.52264171949801\\
70	-2.23693234955811\\
80	-1.918621635566\\
90	-1.51411771350753\\
};
\addplot [color=mycolor7,solid,line width=1pt,mark size=1pt,mark=*,mark options={solid},forget plot]
  table[row sep=crcr]{%
0.0001	-7.61317271802468\\
0.0002	-7.07383644410888\\
0.0003	-6.80068030895105\\
0.0004	-6.31794269385943\\
0.0005	-6.17945208761242\\
0.0006	-5.98019581358\\
0.0007	-5.81803101490402\\
0.0008	-5.80412082564335\\
0.0009	-5.77677850093735\\
0.001	-5.50645480920296\\
0.002	-4.8867815514453\\
0.003	-4.83566834875473\\
0.004	-4.72698624350773\\
0.005	-4.35327401941703\\
0.006	-4.18122818247779\\
0.007	-4.31060844114678\\
0.008	-4.22214878935444\\
0.009	-4.37523891432264\\
0.01	-4.12783729438182\\
0.02	-3.62741165099757\\
0.03	-3.40608557466148\\
0.04	-3.67595932501852\\
0.05	-3.6996850951557\\
0.06	-3.52342463797412\\
0.07	-3.57474233218067\\
0.08	-3.37719764486624\\
0.09	-3.38215204046113\\
0.1	-3.55032028454821\\
0.2	-3.42166051497637\\
0.3	-3.46904867970897\\
0.4	-3.38846177933461\\
0.5	-3.42933810397325\\
0.6	-3.35707853088531\\
0.7	-3.30959697821657\\
0.8	-3.3457852494068\\
0.9	-3.36773812906474\\
1	-3.24063571218837\\
2	-3.35283980410964\\
3	-3.32342124122217\\
4	-3.30207987612189\\
5	-3.33850355897232\\
6	-3.3166607493066\\
7	-3.31357228766659\\
8	-3.31494063643217\\
9	-3.39598955713044\\
10	-3.30615503705657\\
20	-3.29406250803357\\
30	-3.21754251112515\\
40	-3.08675071353794\\
50	-2.86495581682821\\
60	-2.62204352940502\\
70	-2.34008193134981\\
80	-2.03369346679906\\
90	-1.63228962799666\\
};
\addplot [color=mycolor1,solid,line width=1pt,mark size=1pt,mark=*,mark options={solid},forget plot]
  table[row sep=crcr]{%
0.0001	-7.47254742271481\\
0.0002	-6.99558903064937\\
0.0003	-6.73059306246304\\
0.0004	-6.39130221093189\\
0.0005	-6.5057241104358\\
0.0006	-6.13395035898796\\
0.0007	-5.81441862230743\\
0.0008	-5.9854453558013\\
0.0009	-5.70498751270265\\
0.001	-5.43567109981474\\
0.002	-5.01153769116951\\
0.003	-4.93272399168982\\
0.004	-4.46266748913628\\
0.005	-4.45100073631727\\
0.006	-4.16099472789891\\
0.007	-4.1802097830199\\
0.008	-4.21745510147723\\
0.009	-3.83548554898153\\
0.01	-4.15527742578349\\
0.02	-3.77976568174328\\
0.03	-3.67938125333853\\
0.04	-3.63606789654038\\
0.05	-3.77530721186631\\
0.06	-3.65212412684969\\
0.07	-3.55908543581845\\
0.08	-3.62797060575941\\
0.09	-3.6623207448216\\
0.1	-3.62429863327328\\
0.2	-3.61633717368908\\
0.3	-3.42542001757354\\
0.4	-3.52860153891152\\
0.5	-3.55742037792448\\
0.6	-3.45444246255304\\
0.7	-3.50431360290751\\
0.8	-3.59754020764383\\
0.9	-3.54767042468692\\
1	-3.4368583555438\\
2	-3.51898247860692\\
3	-3.50326025885548\\
4	-3.50518235053309\\
5	-3.51139191449263\\
6	-3.502937716433\\
7	-3.48394093844441\\
8	-3.48407747707266\\
9	-3.54219785037449\\
10	-3.49353712916487\\
20	-3.47261982961318\\
30	-3.37572733096008\\
40	-3.2182119471072\\
50	-2.98510116257323\\
60	-2.72649868602104\\
70	-2.44365731216745\\
80	-2.1446213693609\\
90	-1.75524168077122\\
};
\addplot [color=mycolor2,solid,line width=1pt,mark size=1pt,mark=*,mark options={solid},forget plot]
  table[row sep=crcr]{%
0.0001	-7.58987918869338\\
0.0002	-7.26788298956683\\
0.0003	-6.7022616328176\\
0.0004	-6.5713388806925\\
0.0005	-6.45557162541139\\
0.0006	-5.95036091110812\\
0.0007	-6.0211171482435\\
0.0008	-5.78401481462424\\
0.0009	-5.6291466819332\\
0.001	-5.59474802330182\\
0.002	-5.20403301838298\\
0.003	-4.57395674627447\\
0.004	-4.62378027523935\\
0.005	-4.44526900307226\\
0.006	-4.61073611239484\\
0.007	-4.60566425452148\\
0.008	-4.2876056713566\\
0.009	-4.45516510595887\\
0.01	-4.17835430931629\\
0.02	-3.76661071535632\\
0.03	-3.59968140889596\\
0.04	-3.79763425222935\\
0.05	-3.72852930949542\\
0.06	-3.76111015597601\\
0.07	-3.71444414932644\\
0.08	-3.63460377975195\\
0.09	-3.73043245647619\\
0.1	-3.52615608543277\\
0.2	-3.60686714446546\\
0.3	-3.6053848220624\\
0.4	-3.6168692669969\\
0.5	-3.6279346015406\\
0.6	-3.58035337360119\\
0.7	-3.6369863898379\\
0.8	-3.70663979911102\\
0.9	-3.69096323053876\\
1	-3.61123063607531\\
2	-3.70684291380119\\
3	-3.66335086573333\\
4	-3.6743044286875\\
5	-3.65026653347446\\
6	-3.67258051543904\\
7	-3.62749198375443\\
8	-3.65271852180256\\
9	-3.68145763632687\\
10	-3.64670495453712\\
20	-3.62162212748649\\
30	-3.50910472218053\\
40	-3.31186542763224\\
50	-3.06994360125684\\
60	-2.79984396021585\\
70	-2.52012415389028\\
80	-2.22000200118555\\
90	-1.86032460796052\\
};
\addplot [color=mycolor3,solid,line width=1pt,mark size=1pt,mark=*,mark options={solid},forget plot]
  table[row sep=crcr]{%
0.0001	-7.55409110831041\\
0.0002	-7.1545229005477\\
0.0003	-6.61141834762088\\
0.0004	-6.54350672672667\\
0.0005	-6.38373175955667\\
0.0006	-6.16913139807518\\
0.0007	-5.87291295117281\\
0.0008	-5.88595184820571\\
0.0009	-5.73043111014908\\
0.001	-5.68555202068757\\
0.002	-5.20666568679622\\
0.003	-4.94636174924197\\
0.004	-4.71378338489224\\
0.005	-4.6192112331997\\
0.006	-4.69720807865852\\
0.007	-4.41849021406796\\
0.008	-4.23048527949615\\
0.009	-4.05115859865786\\
0.01	-4.07864818844149\\
0.02	-4.20590855500379\\
0.03	-4.0298658208987\\
0.04	-3.76430375368349\\
0.05	-3.90874134147619\\
0.06	-3.80220830515468\\
0.07	-3.83660789760768\\
0.08	-3.74209922853751\\
0.09	-3.72696908712438\\
0.1	-3.78949591690949\\
0.2	-3.83026416347067\\
0.3	-3.77662028481513\\
0.4	-3.80089474972552\\
0.5	-3.79166183109308\\
0.6	-3.72033470593703\\
0.7	-3.7142917271384\\
0.8	-3.80601931674674\\
0.9	-3.74488968685271\\
1	-3.74507370896615\\
2	-3.79677835691054\\
3	-3.79791589276094\\
4	-3.78645120943106\\
5	-3.81002783314484\\
6	-3.77915534624021\\
7	-3.75963582821037\\
8	-3.77868160162701\\
9	-3.81793904837723\\
10	-3.77911210516296\\
20	-3.74842624318806\\
30	-3.61700424601609\\
40	-3.39148835266854\\
50	-3.132141396301\\
60	-2.86658265082792\\
70	-2.59617182770586\\
80	-2.30880067315781\\
90	-1.95890685117519\\
};
\end{axis}
\end{tikzpicture}%

%% file: formLEindex_2.tex
\tikzsetfigurename{formLEindex}

\definecolor{mycolor1}{rgb}{0.00000,0.44700,0.74100}%
\definecolor{mycolor2}{rgb}{0.85000,0.32500,0.09800}%
\definecolor{mycolor3}{rgb}{0.92900,0.69400,0.12500}%
\definecolor{mycolor4}{rgb}{0.49400,0.18400,0.55600}%
\definecolor{mycolor5}{rgb}{0.46600,0.67400,0.18800}%
\definecolor{mycolor6}{rgb}{0.30100,0.74500,0.93300}%
\definecolor{mycolor7}{rgb}{0.63500,0.07800,0.18400}%
\begin{tikzpicture}

\begin{axis}[%
/tikz/line join=bevel,
width=4.25in,
height=2.75in,
xmin=1,
xmax=10,
xlabel={$i$},
xlabel shift={-2pt},
ymin=0,
ymax=0.1,
ylabel={$\lambda_i$: median and deviation},
yticklabel style={
  /pgf/number format/fixed,
  /pgf/number format/precision=4
},
scaled y ticks=false,
xmajorgrids,
ymajorgrids,
legend style={legend cell align=left, align=left, draw=white!15!black}
]
\addplot [color=mycolor1, line width=1pt,mark=*, mark size=1pt, mark options={solid, mycolor1}]
  table[row sep=crcr]{%
1	0.083449725\\
2	0.064299498\\
3	0.046616173\\
4	0.032270989\\
5	0.016496863\\
6	0.0030140255\\
};
\addlegendentry{$L=55$}
\addplot [forget plot,color=mycolor1, line width=1pt,mark=-, mark options={solid, mycolor1}]
  table[row sep=crcr]{%
1	0.086844816\\
1	0.080054634\\
};
\addplot [forget plot,color=mycolor1, line width=1pt,mark=-, mark options={solid, mycolor1}]
  table[row sep=crcr]{%
2	0.066678225\\
2	0.061920771\\
};
\addplot [forget plot,color=mycolor1, line width=1pt,mark=-, mark options={solid, mycolor1}]
  table[row sep=crcr]{%
3	0.049931591\\
3	0.043300755\\
};
\addplot [forget plot,color=mycolor1, line width=1pt,mark=-, mark options={solid, mycolor1}]
  table[row sep=crcr]{%
4	0.034270092\\
4	0.030271886\\
};
\addplot [forget plot,color=mycolor1, line width=1pt,mark=-, mark options={solid, mycolor1}]
  table[row sep=crcr]{%
5	0.017902986\\
5	0.01509074\\
};
\addplot [forget plot,color=mycolor1, line width=1pt,mark=-, mark options={solid, mycolor1}]
  table[row sep=crcr]{%
6	0.0035510198\\
6	0.0024770312\\
};

\addplot [color=mycolor2, line width=1pt,mark=*, mark size=1pt, mark options={solid, mycolor2}]
  table[row sep=crcr]{%
1	0.086232935\\
2	0.068435054\\
3	0.053913915\\
4	0.040668416\\
5	0.028227001\\
6	0.016864053\\
7	0.0031616403\\
};
\addlegendentry{$L=65$}
\addplot [forget plot,color=mycolor2, line width=1pt,mark=-, mark options={solid, mycolor2}]
  table[row sep=crcr]{%
1	0.090145313\\
1	0.082320557\\
};
\addplot [forget plot,color=mycolor2, line width=1pt,mark=-, mark options={solid, mycolor2}]
  table[row sep=crcr]{%
2	0.07174372\\
2	0.065126388\\
};
\addplot [forget plot,color=mycolor2, line width=1pt,mark=-, mark options={solid, mycolor2}]
  table[row sep=crcr]{%
3	0.056495387\\
3	0.051332443\\
};
\addplot [forget plot,color=mycolor2, line width=1pt,mark=-, mark options={solid, mycolor2}]
  table[row sep=crcr]{%
4	0.042526383\\
4	0.038810449\\
};
\addplot [forget plot,color=mycolor2, line width=1pt,mark=-, mark options={solid, mycolor2}]
  table[row sep=crcr]{%
5	0.029959911\\
5	0.026494091\\
};
\addplot [forget plot,color=mycolor2, line width=1pt,mark=-, mark options={solid, mycolor2}]
  table[row sep=crcr]{%
6	0.019219176\\
6	0.01450893\\
};
\addplot [forget plot,color=mycolor2, line width=1pt,mark=-, mark options={solid, mycolor2}]
  table[row sep=crcr]{%
7	0.0041873483\\
7	0.0021359323\\
};

\addplot [color=mycolor3, line width=1pt,mark=*, mark size=1pt, mark options={solid, mycolor3}]
  table[row sep=crcr]{%
1	0.088493924\\
2	0.072869015\\
3	0.058624817\\
4	0.048071848\\
5	0.037606329\\
6	0.026312758\\
7	0.017154006\\
8	0.0069666269\\
};
\addlegendentry{$L=75$}
\addplot [forget plot,color=mycolor3, line width=1pt,mark=-, mark options={solid, mycolor3}]
  table[row sep=crcr]{%
1	0.090673613\\
1	0.086314235\\
};
\addplot [forget plot,color=mycolor3, line width=1pt,mark=-, mark options={solid, mycolor3}]
  table[row sep=crcr]{%
2	0.074773555\\
2	0.070964475\\
};
\addplot [forget plot,color=mycolor3, line width=1pt,mark=-, mark options={solid, mycolor3}]
  table[row sep=crcr]{%
3	0.060001567\\
3	0.057248067\\
};
\addplot [forget plot,color=mycolor3, line width=1pt,mark=-, mark options={solid, mycolor3}]
  table[row sep=crcr]{%
4	0.049349352\\
4	0.046794344\\
};
\addplot [forget plot,color=mycolor3, line width=1pt,mark=-, mark options={solid, mycolor3}]
  table[row sep=crcr]{%
5	0.039731623\\
5	0.035481035\\
};
\addplot [forget plot,color=mycolor3, line width=1pt,mark=-, mark options={solid, mycolor3}]
  table[row sep=crcr]{%
6	0.02769239\\
6	0.024933126\\
};
\addplot [forget plot,color=mycolor3, line width=1pt,mark=-, mark options={solid, mycolor3}]
  table[row sep=crcr]{%
7	0.018378141\\
7	0.015929871\\
};
\addplot [forget plot,color=mycolor3, line width=1pt,mark=-, mark options={solid, mycolor3}]
  table[row sep=crcr]{%
8	0.0087238167\\
8	0.0052094371\\
};

\addplot [color=mycolor4, line width=1pt,mark=*, mark size=1pt, mark options={solid, mycolor4}]
  table[row sep=crcr]{%
1	0.087177231\\
2	0.07468853\\
3	0.062238518\\
4	0.053517134\\
5	0.042492309\\
6	0.034610777\\
7	0.026011249\\
8	0.015997667\\
9	0.0070956772\\
};
\addlegendentry{$L=85$}
\addplot [forget plot,color=mycolor4, line width=1pt,mark=-, mark options={solid, mycolor4}]
  table[row sep=crcr]{%
1	0.089872874\\
1	0.084481588\\
};
\addplot [forget plot,color=mycolor4, line width=1pt,mark=-, mark options={solid, mycolor4}]
  table[row sep=crcr]{%
2	0.076469134\\
2	0.072907926\\
};
\addplot [forget plot,color=mycolor4, line width=1pt,mark=-, mark options={solid, mycolor4}]
  table[row sep=crcr]{%
3	0.064172319\\
3	0.060304717\\
};
\addplot [forget plot,color=mycolor4, line width=1pt,mark=-, mark options={solid, mycolor4}]
  table[row sep=crcr]{%
4	0.055923793\\
4	0.051110475\\
};
\addplot [forget plot,color=mycolor4, line width=1pt,mark=-, mark options={solid, mycolor4}]
  table[row sep=crcr]{%
5	0.044715576\\
5	0.040269042\\
};
\addplot [forget plot,color=mycolor4, line width=1pt,mark=-, mark options={solid, mycolor4}]
  table[row sep=crcr]{%
6	0.036057251\\
6	0.033164303\\
};
\addplot [forget plot,color=mycolor4, line width=1pt,mark=-, mark options={solid, mycolor4}]
  table[row sep=crcr]{%
7	0.027421083\\
7	0.024601415\\
};
\addplot [forget plot,color=mycolor4, line width=1pt,mark=-, mark options={solid, mycolor4}]
  table[row sep=crcr]{%
8	0.017404398\\
8	0.014590936\\
};
\addplot [forget plot,color=mycolor4, line width=1pt,mark=-, mark options={solid, mycolor4}]
  table[row sep=crcr]{%
9	0.0090786587\\
9	0.0051126957\\
};

\addplot [color=mycolor5, line width=1pt,mark=*, mark size=1pt, mark options={solid, mycolor5}]
  table[row sep=crcr]{%
1	0.090751669\\
2	0.076792032\\
3	0.066455825\\
4	0.055832903\\
5	0.047310198\\
6	0.04043328\\
7	0.032446563\\
8	0.024618382\\
9	0.01659835\\
10	0.0081255925\\
};
\addlegendentry{$L=95$}
\addplot [forget plot,color=mycolor5, line width=1pt,mark=-, mark options={solid, mycolor5}]
  table[row sep=crcr]{%
1	0.095490893\\
1	0.086012445\\
};
\addplot [forget plot,color=mycolor5, line width=1pt,mark=-, mark options={solid, mycolor5}]
  table[row sep=crcr]{%
2	0.079454831\\
2	0.074129233\\
};
\addplot [forget plot,color=mycolor5, line width=1pt,mark=-, mark options={solid, mycolor5}]
  table[row sep=crcr]{%
3	0.068465498\\
3	0.064446152\\
};
\addplot [forget plot,color=mycolor5, line width=1pt,mark=-, mark options={solid, mycolor5}]
  table[row sep=crcr]{%
4	0.057024413\\
4	0.054641393\\
};
\addplot [forget plot,color=mycolor5, line width=1pt,mark=-, mark options={solid, mycolor5}]
  table[row sep=crcr]{%
5	0.048955588\\
5	0.045664808\\
};
\addplot [forget plot,color=mycolor5, line width=1pt,mark=-, mark options={solid, mycolor5}]
  table[row sep=crcr]{%
6	0.041748785\\
6	0.039117775\\
};
\addplot [forget plot,color=mycolor5, line width=1pt,mark=-, mark options={solid, mycolor5}]
  table[row sep=crcr]{%
7	0.033766654\\
7	0.031126472\\
};
\addplot [forget plot,color=mycolor5, line width=1pt,mark=-, mark options={solid, mycolor5}]
  table[row sep=crcr]{%
8	0.025651228\\
8	0.023585536\\
};
\addplot [forget plot,color=mycolor5, line width=1pt,mark=-, mark options={solid, mycolor5}]
  table[row sep=crcr]{%
9	0.018339491\\
9	0.014857209\\
};
\addplot [forget plot,color=mycolor5, line width=1pt,mark=-, mark options={solid, mycolor5}]
  table[row sep=crcr]{%
10	0.0087939597\\
10	0.0074572253\\
};

\end{axis}
\end{tikzpicture}

%% file: formLEerror_2.tex
\tikzsetfigurename{formLEerror}

\definecolor{mycolor1}{rgb}{0.00000,0.44700,0.74100}%
\definecolor{mycolor2}{rgb}{0.85000,0.32500,0.09800}%
\definecolor{mycolor3}{rgb}{0.92900,0.69400,0.12500}%
\definecolor{mycolor4}{rgb}{0.49400,0.18400,0.55600}%
\definecolor{mycolor5}{rgb}{0.46600,0.67400,0.18800}%
\definecolor{mycolor6}{rgb}{0.30100,0.74500,0.93300}%
\definecolor{mycolor7}{rgb}{0.63500,0.07800,0.18400}%
\begin{tikzpicture}

\begin{axis}[%
/tikz/line join=bevel,
width=3.5in,
height=2.5in,
xmin=0,
xmax=2,
xlabel={$p$},
xlabel shift={-2pt},
yminorticks=true,
ylabel={residual error},
ylabel shift={1pt},
yticklabel style={
  /pgf/number format/fixed,
  /pgf/number format/precision=4
},
scaled y ticks=false,
xmajorgrids,
ymajorgrids,
legend style={at={(0.02,0.96)},anchor=north west,legend cell align=left, align=left, draw=white!15!black}
]
\addplot [color=mycolor1, line width=1pt]
  table[row sep=crcr]{%
0.02	0.00367428006822522\\
0.04	0.00361070599448849\\
0.06	0.00354740242715804\\
0.08	0.00348439478799292\\
0.1	0.00342171000571346\\
0.12	0.00335937662711026\\
0.14	0.00329742493628783\\
0.16	0.0032358870823698\\
0.18	0.00317479721591942\\
0.2	0.0031141916342263\\
0.22	0.00305410893546856\\
0.24	0.00299459018157059\\
0.26	0.00293567906933\\
0.28	0.00287742210906902\\
0.3	0.00281986880966369\\
0.32	0.00276307186830081\\
0.34	0.00270708736269271\\
0.36	0.00265197494272447\\
0.38	0.00259779801760238\\
0.4	0.00254462393349998\\
0.42	0.00249252413545136\\
0.44	0.00244157430581851\\
0.46	0.00239185447007338\\
0.48	0.00234344905891674\\
0.5	0.00229644691396238\\
0.52	0.00225094122243812\\
0.54	0.00220702936472733\\
0.56	0.00216481265727873\\
0.58	0.00212439597268355\\
0.6	0.0020858872188442\\
0.62	0.00204939666046465\\
0.64	0.00201503606892422\\
0.66	0.00198291769127535\\
0.68	0.00195315303587517\\
0.7	0.00192585148111124\\
0.72	0.00190111872467431\\
0.74	0.00187905510341619\\
0.76	0.00185975382720444\\
0.78	0.00184329918318288\\
0.8	0.00182976477801497\\
0.82	0.00181921189343716\\
0.84	0.0018116880332853\\
0.86	0.00180722573695505\\
0.88	0.00180584172452788\\
0.9	0.00180753642287865\\
0.92	0.00181229390120215\\
0.94	0.00182008222053865\\
0.96	0.00183085417750685\\
0.98	0.00184454840011715\\
1	0.00186109073546284\\
1.02	0.00188039585683214\\
1.04	0.00190236901205371\\
1.06	0.00192690783551104\\
1.08	0.00195390415233648\\
1.1	0.00198324571343214\\
1.12	0.0020148178125699\\
1.14	0.00204850475035142\\
1.16	0.00208419112295746\\
1.18	0.0021217629254226\\
1.2	0.00216110846904452\\
1.22	0.00220211912022064\\
1.24	0.00224468987350876\\
1.26	0.00228871977523008\\
1.28	0.00233411221577565\\
1.3	0.00238077510928799\\
1.32	0.00242862097891265\\
1.34	0.00247756696465799\\
1.36	0.00252753476932713\\
1.38	0.0025784505561984\\
1.4	0.00263024481028447\\
1.42	0.00268285217320307\\
1.44	0.00273621126001407\\
1.46	0.00279026446485934\\
1.48	0.00284495776090188\\
1.5	0.00290024049890203\\
1.52	0.00295606520778469\\
1.54	0.00301238739972838\\
1.56	0.00306916538162838\\
1.58	0.00312636007423461\\
1.6	0.00318393483982154\\
1.62	0.00324185531889731\\
1.64	0.00330008927618545\\
1.66	0.00335860645590297\\
1.68	0.00341737844620086\\
1.7	0.00347637855251759\\
1.72	0.00353558167951391\\
1.74	0.00359496422120234\\
1.76	0.00365450395884994\\
1.78	0.00371417996621489\\
1.8	0.00377397252167122\\
1.82	0.00383386302677946\\
1.84	0.00389383393087067\\
1.86	0.00395386866122623\\
1.88	0.00401395155845386\\
1.9	0.00407406781668038\\
1.92	0.00413420342820302\\
1.94	0.00419434513226307\\
1.96	0.00425448036762715\\
1.98	0.00431459722868306\\
2	0.00437468442477765\\
};
\addlegendentry{\textsc{rms}}

\addplot [color=mycolor2, line width=1pt]
  table[row sep=crcr]{%
0.02	0.00218366659522587\\
0.04	0.00218395134785328\\
0.06	0.00218360687005291\\
0.08	0.00217894884386412\\
0.1	0.00218102978371311\\
0.12	0.00217879731288551\\
0.14	0.00217593588791419\\
0.16	0.0021540368733773\\
0.18	0.00204920831191292\\
0.2	0.001974612456981\\
0.22	0.00197096692686315\\
0.24	0.00196730880920122\\
0.26	0.00196954491076112\\
0.28	0.00194036851089908\\
0.3	0.00194465632056006\\
0.32	0.00189012070855259\\
0.34	0.00190518230057472\\
0.36	0.00189599362699783\\
0.38	0.0018487457124267\\
0.4	0.00188945076789201\\
0.42	0.00182757398460737\\
0.44	0.00178966010919848\\
0.46	0.00172186025171013\\
0.48	0.00168736938154552\\
0.5	0.00170730585553282\\
0.52	0.0016959649307378\\
0.54	0.00167646751293168\\
0.56	0.00163356718788892\\
0.58	0.00155903857079518\\
0.6	0.00155323399814319\\
0.62	0.00160749581041018\\
0.64	0.00157879044449274\\
0.66	0.00154172904103834\\
0.68	0.00151240286790362\\
0.7	0.0015216192980657\\
0.72	0.00149530848181252\\
0.74	0.00140010261207173\\
0.76	0.00141128546469468\\
0.78	0.00142252013671404\\
0.8	0.00143380627575387\\
0.82	0.00145301549951773\\
0.84	0.00144710267543623\\
0.86	0.00146817409862768\\
0.88	0.00146066143141765\\
0.9	0.00142837528696291\\
0.92	0.00140979783548825\\
0.94	0.00145475933282971\\
0.96	0.00140647727830481\\
0.98	0.00150401384326171\\
1	0.00145862358082444\\
1.02	0.00145801284081405\\
1.04	0.00149624291402686\\
1.06	0.00154522213893883\\
1.08	0.00157092246709545\\
1.1	0.00154316346347887\\
1.12	0.00160923313240914\\
1.14	0.0016199256112549\\
1.16	0.00162054486283964\\
1.18	0.00165753082755626\\
1.2	0.0016697677112655\\
1.22	0.00165253480183068\\
1.24	0.00166863027662219\\
1.26	0.00171970418760403\\
1.28	0.00176568425853066\\
1.3	0.00181153882532469\\
1.32	0.00186742712693779\\
1.34	0.00187004572461663\\
1.36	0.00189012665094208\\
1.38	0.00190960542906714\\
1.4	0.00192848318039299\\
1.42	0.0020209930304044\\
1.44	0.00206019461529233\\
1.46	0.00203499503170104\\
1.48	0.00208392292822683\\
1.5	0.00210189688647375\\
1.52	0.00217937209995923\\
1.54	0.00223651609646439\\
1.56	0.00222238882175789\\
1.58	0.00223189238901918\\
1.6	0.00228740057240358\\
1.62	0.00234286369007444\\
1.64	0.00239827526455506\\
1.66	0.00245362887814034\\
1.68	0.0024781671344232\\
1.7	0.00249316258754314\\
1.72	0.00249840583524174\\
1.74	0.0025217883161819\\
1.76	0.002544280811159\\
1.78	0.00260461112650604\\
1.8	0.00261794598097979\\
1.82	0.00263126283481282\\
1.84	0.00264456128719705\\
1.86	0.00265784093942257\\
1.88	0.00272516846252649\\
1.9	0.00283371570029863\\
1.92	0.00289169091313149\\
1.94	0.00289365351062713\\
1.96	0.00292577373422088\\
1.98	0.00291576987457289\\
2	0.00290510354877743\\
};
\addlegendentry{\textsc{mad}}

\end{axis}
\end{tikzpicture}

%% file: ksperiodic_dky.tex
\definecolor{mycolor1}{rgb}{0.00000,0.44700,0.74100}%
\definecolor{mycolor2}{rgb}{0.85000,0.32500,0.09800}%
\definecolor{mycolor3}{rgb}{0.92900,0.69400,0.12500}%
\definecolor{mycolor4}{rgb}{0.49400,0.18400,0.55600}%
\definecolor{mycolor5}{rgb}{0.46600,0.67400,0.18800}%
\definecolor{mycolor6}{rgb}{0.30100,0.74500,0.93300}%
\definecolor{mycolor7}{rgb}{0.63500,0.07800,0.18400}%
\begin{tikzpicture}

\begin{axis}[%
/tikz/line join=bevel,
width=4.25in,
height=1.75in,
at={(0in,0in)},
scale only axis,
xmin=0,
xmax=100,
xmajorgrids,
xlabel={$L$},
ymin=0,
ymax=23,
ymajorgrids,
axis background/.style={fill=white},
legend style={at={(0.02,0.96)},anchor=north west,legend cell align=left,align=left,draw=white!15!black}
]
\addplot [color=mycolor1,solid,line width=1pt]
  table[row sep=crcr]{%
0.1	0\\
0.2	0\\
0.3	0\\
0.4	0\\
0.5	0\\
0.6	0\\
0.7	0\\
0.8	0\\
0.9	0\\
1	0\\
1.1	0\\
1.2	0\\
1.3	0\\
1.4	0\\
1.5	0\\
1.6	0\\
1.7	0\\
1.8	0\\
1.9	0\\
2	0\\
2.1	0\\
2.2	0\\
2.3	0\\
2.4	0\\
2.5	0\\
2.6	0\\
2.7	0\\
2.8	0\\
2.9	0\\
3	0\\
3.1	0\\
3.2	0\\
3.3	0\\
3.4	0\\
3.5	0\\
3.6	0\\
3.7	0\\
3.8	0\\
3.9	0\\
4	0\\
4.1	0\\
4.2	0\\
4.3	0\\
4.4	0\\
4.5	0\\
4.6	0\\
4.7	0\\
4.8	0\\
4.9	0\\
5	0\\
5.1	0\\
5.2	0\\
5.3	0\\
5.4	0\\
5.5	0\\
5.6	0\\
5.7	0\\
5.8	0\\
5.9	0\\
6	0\\
6.1	0\\
6.2	0\\
6.3	1.65890743134896\\
6.4	1.70854296610162\\
6.5	1.31701815273784\\
6.6	1.60456537039204\\
6.7	1.60938182154853\\
6.8	1.76175876399696\\
6.9	1.66860474933502\\
7	1.55781947275818\\
7.1	1.73061443544364\\
7.2	1.64866218125129\\
7.3	1.70099013441241\\
7.4	1.65348939590849\\
7.5	1.78874182877377\\
7.6	1.67775882718026\\
7.7	1.71839690508091\\
7.8	1.70299809891133\\
7.9	1.55803276873283\\
8	1.76326454155456\\
8.1	1.75710329718997\\
8.2	1.44368621295835\\
8.3	1.53594001542072\\
8.4	0\\
8.5	1.65788052695905\\
8.6	1.69925850607657\\
8.7	1.83836324743396\\
8.8	1.56006944909418\\
8.9	1.76694540744873\\
9	1.68648407654058\\
9.1	1.80119495309645\\
9.2	1.55481825280356\\
9.3	1.77298634994757\\
9.4	1.71922257896962\\
9.5	1.36039200236479\\
9.6	1.38851306769884\\
9.7	1.36464428134872\\
9.8	1.6136430022694\\
9.9	1.68793628847208\\
10	1.84647216189227\\
10.1	1.68326713648262\\
10.2	1.48676538315215\\
10.3	1.85842470276686\\
10.4	1.63279857348358\\
10.5	1.7534953877526\\
10.6	1.59091525856225\\
10.7	1.67929391524392\\
10.8	1.57086656731437\\
10.9	1.75611961778182\\
11	1.75696650762599\\
11.1	1.87266908807651\\
11.2	1.72595727892272\\
11.3	2.06297004972936\\
11.4	1.93013749571587\\
11.5	1.74773813946751\\
11.6	1.96309426652498\\
11.7	1.66104549797252\\
11.8	1.74449378211831\\
11.9	1.51557535638458\\
12	1.66303840684752\\
12.1	1.65598014197933\\
12.2	1.6891889609562\\
12.3	1.64558452523826\\
12.4	1.50915716392013\\
12.5	1.75003506772742\\
12.6	1.79751012486061\\
12.7	1.86985074174142\\
12.8	1.81829285179003\\
12.9	1.65424963084174\\
13	1.70983217238284\\
13.1	3.25574029540908\\
13.2	3.19422271286204\\
13.3	3.26454695109791\\
13.4	3.29568499402471\\
13.5	3.25885091702577\\
13.6	3.29386460311263\\
13.7	3.24200020542178\\
13.8	3.16038570796265\\
13.9	3.17847957466534\\
14	3.12842375987197\\
14.1	3.10558326079378\\
14.2	3.12833216514757\\
14.3	3.1145985805\\
14.4	3.08245130548728\\
14.5	3.08248489069612\\
14.6	3.06180156334322\\
14.7	3.03897974126745\\
14.8	3.02366882016737\\
14.9	3.00297072746094\\
15	1.63537645098102\\
15.1	1.55017452003194\\
15.2	1.54471066439514\\
15.3	1.67991166021885\\
15.4	1.43814932216567\\
15.5	1.61633973964303\\
15.6	1.59654223039277\\
15.7	1.53082304064733\\
15.8	1.48023512908599\\
15.9	1.66372879457668\\
16	1.54132261625047\\
16.1	1.58813569288931\\
16.2	1.67059301705391\\
16.3	1.54661802687858\\
16.4	1.53083240764604\\
16.5	1.60086309034132\\
16.6	1.62428059364761\\
16.7	1.64689727832043\\
16.8	1.63551494423234\\
16.9	1.6500888744128\\
17	1.68886447523928\\
17.1	1.55466842258659\\
17.2	1.52227385954574\\
17.3	3.10455041176328\\
17.4	2.5990733311786\\
17.5	2.45292128827926\\
17.6	2.76308476025496\\
17.7	2.7407583819063\\
17.8	2.58303288160802\\
17.9	2.58886621220356\\
18	2.68012728808373\\
18.1	4.042283038276\\
18.2	4.05222590766294\\
18.3	4.0236097358342\\
18.4	4.0812234385157\\
18.5	4.05086023354762\\
18.6	4.05392757614693\\
18.7	4.0382071503488\\
18.8	4.05906199929041\\
18.9	4.12318331165175\\
19	4.12378362473544\\
19.1	4.14925889601516\\
19.2	4.15162715621855\\
19.3	4.06778212163359\\
19.4	4.19747662062977\\
19.5	4.2007638319004\\
19.6	4.22540405048079\\
19.7	4.22449996403994\\
19.8	4.40411568569524\\
19.9	4.20384054550144\\
20	4.28203613371491\\
20.1	4.27764125382036\\
20.2	4.30572626992845\\
20.3	4.50088084628086\\
20.4	4.45377057335091\\
20.5	4.50842345059612\\
20.6	4.57561694956972\\
20.7	4.86671338650429\\
20.8	4.57062011667795\\
20.9	3.221214387744\\
21	2.82866416614012\\
21.1	3.09100444618927\\
21.2	4.21910825579887\\
21.3	5.03018086957421\\
21.4	5.05112442057571\\
21.5	5.11276949951675\\
21.6	5.09975715387119\\
21.7	5.141950276498\\
21.8	5.14835464052541\\
21.9	5.15165452824013\\
22	5.19848295869525\\
22.1	5.23212941311931\\
22.2	5.30229077127081\\
22.3	1.51731179299057\\
22.4	1.73135490677231\\
22.5	1.04617419374648\\
22.6	1.68362821445552\\
22.7	1.60066243208172\\
22.8	1.42471763266238\\
22.9	2.76985059506498\\
23	2.64816415206677\\
23.1	1.56412266626931\\
23.2	4.68037457804476\\
23.3	5.00688472717451\\
23.4	2.8830568423958\\
23.5	5.35884147047169\\
23.6	1.47017257101261\\
23.7	1.63391772594827\\
23.8	2.69728961034629\\
23.9	1.62890273455636\\
24	1.55957119614506\\
24.1	1.27683673341501\\
24.2	1.20153345471251\\
24.3	1.48624954712404\\
24.4	1.631981325424\\
24.5	1.53512155196032\\
24.6	1.80893794659252\\
24.7	1.55844350858128\\
24.8	1.57444754891469\\
24.9	1.71411743677356\\
25	1.566871429417\\
25.1	1.57501031172651\\
25.2	1.54606154071346\\
25.3	1.49690711730702\\
25.4	1.65577796553226\\
25.5	1.49613401176203\\
25.6	2.16913503453965\\
25.7	2.85128106083308\\
25.8	4.34827964313723\\
25.9	5.24312899761035\\
26	5.7595262298694\\
26.1	5.67191632643735\\
26.2	5.87646361332884\\
26.3	2.31999567038043\\
26.4	1.86881521422834\\
26.5	1.41638056181888\\
26.6	1.99271543653571\\
26.7	1.63371323793153\\
26.8	1.23621164680327\\
26.9	1.81819097032102\\
27	1.48686592797713\\
27.1	1.84911831402564\\
27.2	1.56150008463437\\
27.3	1.7090683108927\\
27.4	1.57334089454784\\
27.5	1.77414466645176\\
27.6	1.57701322591207\\
27.7	1.49839681570909\\
27.8	1.74223420574317\\
27.9	1.63107149439484\\
28	1.57857485992868\\
28.1	1.64428650222291\\
28.2	1.54182405154797\\
28.3	1.72901724648802\\
28.4	1.62285125508259\\
28.5	1.79216393510806\\
28.6	2.35134478150086\\
28.7	3.0515109292005\\
28.8	2.80142240304434\\
28.9	2.45209682711407\\
29	2.96013248793243\\
29.1	4.32508729637468\\
29.2	2.7188595376295\\
29.3	2.49687399554849\\
29.4	2.58226895891396\\
29.5	2.67590555887456\\
29.6	2.79859210050705\\
29.7	3.0072860152461\\
29.8	2.91778926430975\\
29.9	5.74221752505062\\
30	6.5923707542921\\
30.1	6.65478975735127\\
30.2	6.72690454214833\\
30.3	6.78299394683989\\
30.4	3.21758901939482\\
30.5	1.56204543066519\\
30.6	1.46574673191219\\
30.7	1.56508737821766\\
30.8	1.39397978243599\\
30.9	1.48202221712606\\
31	1.58184488061915\\
31.1	1.59016014545784\\
31.2	1.51389384092659\\
31.3	1.39721092464604\\
31.4	1.43947540190668\\
31.5	1.49414643003683\\
31.6	1.54130236042887\\
31.7	1.4927942600416\\
31.8	1.52951226034263\\
31.9	1.44099094917656\\
32	1.54107490180837\\
32.1	1.52695576520466\\
32.2	1.68876776681893\\
32.3	1.49959858433978\\
32.4	1.59175839282208\\
32.5	1.47830841933149\\
32.6	1.58662772527894\\
32.7	1.52343230730515\\
32.8	1.61195736611253\\
32.9	1.2589693809149\\
33	1.53919926353151\\
33.1	1.55506854096048\\
33.2	1.40960652625391\\
33.3	1.55655193183268\\
33.4	1.49322489695556\\
33.5	1.62995844420785\\
33.6	1.82125262468334\\
33.7	3.06176388816351\\
33.8	6.0676886636833\\
33.9	6.32722077684739\\
34	7.31010072218253\\
34.1	7.20846309958142\\
34.2	7.42242803040432\\
34.3	7.53109604948317\\
34.4	7.47404427398944\\
34.5	7.55041559943563\\
34.6	7.69995428745186\\
34.7	7.629520865441\\
34.8	7.73143584044359\\
34.9	7.80136179332644\\
35	7.97739943010973\\
35.1	7.90142767832124\\
35.2	8.07104166483298\\
35.3	8.0499464564311\\
35.4	8.02262864526579\\
35.5	8.01527187326757\\
35.6	8.2032107424888\\
35.7	8.20281707859311\\
35.8	8.26673298245219\\
35.9	8.10165354265934\\
36	8.22892096867633\\
36.1	8.16733157403813\\
36.2	8.25115767347195\\
36.3	8.22033018676362\\
36.4	8.29389262751669\\
36.5	8.25355391796205\\
36.6	8.44501492544223\\
36.7	8.39162712711585\\
36.8	1.60194841877258\\
36.9	8.29825596009686\\
37	8.30539279017149\\
37.1	8.36446227397633\\
37.2	8.26461883074308\\
37.3	8.56145162216731\\
37.4	8.30234463112006\\
37.5	8.54821857820183\\
37.6	8.41673921967738\\
37.7	8.51957452369987\\
37.8	8.57373586067506\\
37.9	8.61988731258905\\
38	8.55247866306474\\
38.1	8.55196964285509\\
38.2	8.29056813809167\\
38.3	8.66094829021652\\
38.4	8.65932308735531\\
38.5	8.76737574352255\\
38.6	8.74524952691578\\
38.7	8.82837160652176\\
38.8	8.57514763713495\\
38.9	8.69219961726691\\
39	8.58299532770663\\
39.1	8.69693836958989\\
39.2	8.80593934674864\\
39.3	1.58831527843142\\
39.4	9.10382593016843\\
39.5	8.224754827216\\
39.6	5.39929819719918\\
39.7	1.62825978564125\\
39.8	1.5431759032609\\
39.9	1.42276318226578\\
40	1.51868898904215\\
40.1	6.28771655121968\\
40.2	1.49451420890574\\
40.3	1.58738095072622\\
40.4	1.63140016374136\\
40.5	1.45287391030742\\
40.6	1.34756177234508\\
40.7	1.47558735075838\\
40.8	1.33791987037394\\
40.9	1.22814834773152\\
41	1.57590101124374\\
41.1	1.4266201579345\\
41.2	1.51444825487879\\
41.3	1.52769738338535\\
41.4	1.40759724749361\\
41.5	1.67753854296813\\
41.6	1.82932595579136\\
41.7	2.6181350620688\\
41.8	9.4185853419595\\
41.9	9.12630369534169\\
42	7.16429405213044\\
42.1	9.15986156413167\\
42.2	9.35797977782834\\
42.3	9.30931484363458\\
42.4	9.30527962397076\\
42.5	9.41647967530471\\
42.6	9.46953530300228\\
42.7	9.32124801049247\\
42.8	9.64917977428815\\
42.9	9.48801195186814\\
43	9.61825259247828\\
43.1	9.565106123897\\
43.2	9.59191415332175\\
43.3	9.79600543610044\\
43.4	9.67870900837475\\
43.5	9.66141157961141\\
43.6	9.88636306574028\\
43.7	9.75658647493875\\
43.8	9.7924185494782\\
43.9	9.84086849082459\\
44	9.7906970435275\\
44.1	9.91448606152588\\
44.2	9.84151378264315\\
44.3	9.91141737353283\\
44.4	9.89650206924643\\
44.5	9.99711801830787\\
44.6	9.85215096020066\\
44.7	10.1872199056818\\
44.8	10.0480198544857\\
44.9	10.0619939255697\\
45	10.1673323082304\\
45.1	10.173833594238\\
45.2	10.1926408032364\\
45.3	10.2555239460745\\
45.4	10.1115840734244\\
45.5	10.2611867798548\\
45.6	10.3501382542156\\
45.7	10.287025700638\\
45.8	10.2463714145802\\
45.9	10.1307559009446\\
46	1.42015880675481\\
46.1	10.4619646353405\\
46.2	10.3706397848807\\
46.3	10.3569715377443\\
46.4	10.5164695914248\\
46.5	10.3834927969002\\
46.6	10.4525544277713\\
46.7	10.3663615275163\\
46.8	3.55616150958567\\
46.9	10.4443678879231\\
47	10.5954600171688\\
47.1	10.5288700206994\\
47.2	10.5843901253468\\
47.3	10.7355267831672\\
47.4	10.7324391873226\\
47.5	10.6928004403586\\
47.6	5.40164805376961\\
47.7	10.7181717008623\\
47.8	11.0041347691456\\
47.9	1.4634722342888\\
48	1.50975515525479\\
48.1	6.55558880682631\\
48.2	10.4121883190269\\
48.3	1.54677169284439\\
48.4	1.49647025307231\\
48.5	11.0758071300268\\
48.6	1.41098173493972\\
48.7	6.49051452104115\\
48.8	11.0893164190439\\
48.9	1.5690217032055\\
49	10.8576155409366\\
49.1	10.5913396510543\\
49.2	1.52293398372157\\
49.3	6.72462975387047\\
49.4	10.6318503667775\\
49.5	1.76530553029352\\
49.6	2.08228001347221\\
49.7	11.1531602679463\\
49.8	10.5315081567987\\
49.9	11.2083999192192\\
50	11.1161033006974\\
50.1	11.2380274928933\\
50.2	1.49497707417134\\
50.3	11.2379028057733\\
50.4	11.3607710298401\\
50.5	11.467305073745\\
50.6	11.1873262034767\\
50.7	11.3673742303439\\
50.8	11.6201613712746\\
50.9	11.4102552855156\\
51	11.3686060397559\\
51.1	11.4912371955366\\
51.2	11.3219328694222\\
51.3	11.3210556081385\\
51.4	11.5156219154791\\
51.5	11.6476942537725\\
51.6	2.966949976587\\
51.7	11.5792303698175\\
51.8	11.8984476897977\\
51.9	4.39472484918591\\
52	11.7053741654224\\
52.1	11.1308163080662\\
52.2	11.5709000933717\\
52.3	11.6518039435124\\
52.4	11.7670357539742\\
52.5	3.34343390692406\\
52.6	11.7263862613368\\
52.7	11.7858342759223\\
52.8	11.9551594835662\\
52.9	11.9045806386495\\
53	12.0136074770036\\
53.1	11.9123516931411\\
53.2	12.0906312432738\\
53.3	12.0859609461666\\
53.4	12.0962014067818\\
53.5	12.0558273659061\\
53.6	11.9979958662299\\
53.7	12.1355949122655\\
53.8	12.1388956687716\\
53.9	12.3743925764584\\
54	12.1398851039141\\
54.1	12.1838209508766\\
54.2	12.1473735185496\\
54.3	12.1679676715446\\
54.4	12.3137720283326\\
54.5	12.0351444148316\\
54.6	12.2543603450264\\
54.7	12.3318808338639\\
54.8	12.2277180477402\\
54.9	12.2741364230424\\
55	1.56807194575779\\
55.1	12.3602482545151\\
55.2	12.6787838244339\\
55.3	12.3927601396603\\
55.4	12.3034693559251\\
55.5	12.3966495842439\\
55.6	12.6319947813314\\
55.7	12.6203624323239\\
55.8	12.6569871913836\\
55.9	12.591167715339\\
56	12.5589080792146\\
56.1	1.53630288860567\\
56.2	12.7051081972956\\
56.3	12.6685083526751\\
56.4	1.47552763926732\\
56.5	12.729709561539\\
56.6	13.0053818094685\\
56.7	12.6644061473722\\
56.8	12.7304415095349\\
56.9	13.0045775397823\\
57	12.9216804500549\\
57.1	10.9191823884756\\
57.2	12.8891978446872\\
57.3	12.8461356412889\\
57.4	12.8836431226992\\
57.5	12.1404699977017\\
57.6	13.104376471596\\
57.7	13.188205468971\\
57.8	13.1814530270456\\
57.9	12.985935509669\\
58	13.0765129935653\\
58.1	13.1155216373933\\
58.2	13.2617924865959\\
58.3	12.671688131512\\
58.4	13.0290723956513\\
58.5	13.1401119639352\\
58.6	13.2312443886263\\
58.7	13.1511173305227\\
58.8	13.3750486845381\\
58.9	13.2316221266007\\
59	13.3149341572779\\
59.1	13.3047093287684\\
59.2	13.3813500605985\\
59.3	13.4909104121222\\
59.4	13.2528223482924\\
59.5	13.309173524145\\
59.6	13.3066788807022\\
59.7	13.5863760960589\\
59.8	13.3249495532745\\
59.9	13.3607159187147\\
60	13.5611882966392\\
60.1	13.4228387244719\\
60.2	13.480689521919\\
60.3	13.4105064647228\\
60.4	13.5788752345649\\
60.5	13.7283448710604\\
60.6	13.6132617994385\\
60.7	13.7448902033863\\
60.8	13.8757581269447\\
60.9	13.7465963990243\\
61	13.7883762055538\\
61.1	13.5436174687195\\
61.2	13.7797218009067\\
61.3	13.7482197334322\\
61.4	13.7474611353213\\
61.5	13.902864589268\\
61.6	14.1678382641732\\
61.7	13.8929548066265\\
61.8	13.7205312673771\\
61.9	13.8171855854533\\
62	13.9523126962601\\
62.1	14.0412894502692\\
62.2	14.0570532210981\\
62.3	14.0551403805594\\
62.4	13.9310201058213\\
62.5	14.0335172964118\\
62.6	14.0133072491696\\
62.7	14.0064789542734\\
62.8	14.1458473364053\\
62.9	14.0670557132319\\
63	14.2702382035326\\
63.1	14.1906250356636\\
63.2	14.1396494455022\\
63.3	14.3715772179819\\
63.4	14.2108328108432\\
63.5	14.2536071877957\\
63.6	14.2830055943817\\
63.7	14.5426356715173\\
63.8	14.258378722578\\
63.9	14.3624312625749\\
64	14.384752235451\\
64.1	1.52928996014676\\
64.2	14.3784724844834\\
64.3	14.6277121743896\\
64.4	14.3350054332175\\
64.5	14.4491012049498\\
64.6	14.3902938164777\\
64.7	14.4013806951279\\
64.8	14.3351409616249\\
64.9	14.651489362885\\
65	14.726939085133\\
65.1	14.5819074271973\\
65.2	14.790912915047\\
65.3	14.5452633222485\\
65.4	1.84553342048069\\
65.5	14.5055294552005\\
65.6	14.8052241310101\\
65.7	14.8333584143712\\
65.8	14.7728729366627\\
65.9	14.7354904457279\\
66	14.8295807390201\\
66.1	14.8797956440848\\
66.2	14.8312737657111\\
66.3	14.9130151634975\\
66.4	14.9461024497106\\
66.5	15.0883623880792\\
66.6	15.0226690649941\\
66.7	15.0034000307659\\
66.8	14.9942644632064\\
66.9	14.1643681779415\\
67	14.8922573746162\\
67.1	15.1634358446272\\
67.2	15.1263876986431\\
67.3	14.9899057168622\\
67.4	15.0879965395598\\
67.5	15.1198041014563\\
67.6	15.1577076836445\\
67.7	15.2402232173527\\
67.8	15.4192952497398\\
67.9	15.3188010513895\\
68	15.2374619373131\\
68.1	15.2689245759686\\
68.2	15.2683296054898\\
68.3	15.4428358707692\\
68.4	15.3137878555151\\
68.5	15.332825933779\\
68.6	15.441580841183\\
68.7	15.3699199456222\\
68.8	15.5044602603371\\
68.9	15.3156999470957\\
69	15.4887825053256\\
69.1	15.4943357908698\\
69.2	15.3539356773347\\
69.3	15.6573047343536\\
69.4	15.6069760751274\\
69.5	15.6901484163757\\
69.6	15.5837542725658\\
69.7	15.620021870847\\
69.8	15.7412853758839\\
69.9	15.6593544667771\\
70	15.8465036936557\\
70.1	15.8004724479133\\
70.2	15.7762969100318\\
70.3	15.7945618633669\\
70.4	15.8339848354728\\
70.5	15.8841964352967\\
70.6	15.8514654135292\\
70.7	15.8028234014969\\
70.8	15.8595509107918\\
70.9	15.79765205887\\
71	16.130976219895\\
71.1	16.0894873214678\\
71.2	15.8785652253788\\
71.3	15.9665803928471\\
71.4	16.0939499837973\\
71.5	16.0260012748699\\
71.6	16.0919602760183\\
71.7	15.997593128821\\
71.8	10.0603485633196\\
71.9	16.3736681811127\\
72	16.071804992493\\
72.1	16.178363971642\\
72.2	16.109049634811\\
72.3	16.3807201576735\\
72.4	16.2059703097711\\
72.5	16.2734467447335\\
72.6	16.202507665493\\
72.7	16.417899087198\\
72.8	16.2142414676903\\
72.9	16.4811734165617\\
73	16.4632644923359\\
73.1	16.3034260430605\\
73.2	16.6347802605313\\
73.3	16.209236585312\\
73.4	16.5423571805736\\
73.5	16.4373351870999\\
73.6	2.21955918937441\\
73.7	16.5816738917545\\
73.8	16.6182056787988\\
73.9	16.5925265781921\\
74	16.7235429555161\\
74.1	16.803469917031\\
74.2	16.5980504723062\\
74.3	16.5810419681834\\
74.4	16.8282843702382\\
74.5	16.8635877561893\\
74.6	16.6916935293102\\
74.7	16.7452073704702\\
74.8	16.8184729646438\\
74.9	16.904960280432\\
75	16.6936554908773\\
75.1	16.8039693213686\\
75.2	16.8699567784882\\
75.3	16.788477118414\\
75.4	16.8747936486268\\
75.5	17.0007322152988\\
75.6	17.0160674319276\\
75.7	17.1735965338666\\
75.8	16.9164363821408\\
75.9	17.1504731005947\\
76	17.0616034184739\\
76.1	17.0775703317725\\
76.2	17.1758046133693\\
76.3	17.1561883065136\\
76.4	17.0652641146144\\
76.5	17.2377885769213\\
76.6	17.2814301315931\\
76.7	17.3496689283029\\
76.8	17.4613432224193\\
76.9	17.2053927637629\\
77	17.1212972840355\\
77.1	17.3931872634171\\
77.2	17.2733663307447\\
77.3	17.4403021502378\\
77.4	17.3638237586814\\
77.5	17.273455315395\\
77.6	17.3766166828389\\
77.7	17.5069245407177\\
77.8	17.5115077390857\\
77.9	17.4502202289082\\
78	17.5957110039879\\
78.1	17.30475305216\\
78.2	17.6239592458398\\
78.3	17.4775843422507\\
78.4	17.8567580576209\\
78.5	17.842384670863\\
78.6	17.6356507071937\\
78.7	17.6860019013083\\
78.8	17.4616983576139\\
78.9	17.601186532567\\
79	17.7364590700102\\
79.1	17.9460162287299\\
79.2	17.6322132892081\\
79.3	18.0426930339714\\
79.4	17.6576826598916\\
79.5	18.0543376145594\\
79.6	17.8475694957497\\
79.7	17.9367144217451\\
79.8	18.0863693446268\\
79.9	17.8672321698615\\
80	17.8040872276797\\
80.1	17.8558442237547\\
80.2	17.9716196362268\\
80.3	17.9709649828322\\
80.4	18.1901815730626\\
80.5	18.1494132333479\\
80.6	18.0752173653275\\
80.7	18.0974382513112\\
80.8	18.1471050194669\\
80.9	18.0723652745756\\
81	18.1963432334433\\
81.1	18.2946366907672\\
81.2	18.0667481172268\\
81.3	18.1489110267248\\
81.4	18.2089773508365\\
81.5	18.3608216790096\\
81.6	18.2378286385726\\
81.7	18.790134524985\\
81.8	18.3594746747597\\
81.9	18.1738995114728\\
82	18.4293240394836\\
82.1	18.2911125081674\\
82.2	18.4758816037774\\
82.3	18.3640874799425\\
82.4	18.4339275395624\\
82.5	18.5214256521132\\
82.6	18.4823260305364\\
82.7	18.202997737157\\
82.8	18.4092851357345\\
82.9	18.344605288438\\
83	18.8098806929827\\
83.1	18.772857050079\\
83.2	18.7087032335421\\
83.3	18.7443829859789\\
83.4	18.6332011741099\\
83.5	18.9003477279239\\
83.6	18.6218963263039\\
83.7	18.5922105413575\\
83.8	19.0666585516558\\
83.9	18.5938256990666\\
84	18.7804674556901\\
84.1	18.71749799241\\
84.2	19.0865348650204\\
84.3	18.9644854904268\\
84.4	18.9804543945151\\
84.5	18.9672313922496\\
84.6	19.0046656315459\\
84.7	19.0504719364637\\
84.8	18.9865986129155\\
84.9	19.1117624954999\\
85	18.9876276336589\\
85.1	18.8959598527614\\
85.2	19.2129891565304\\
85.3	19.2339060201131\\
85.4	19.2730093136595\\
85.5	19.206697055299\\
85.6	19.0669107267579\\
85.7	19.1694354563026\\
85.8	19.2430992575242\\
85.9	19.2650565828458\\
86	19.4557881157225\\
86.1	19.0789489826262\\
86.2	19.4043846886185\\
86.3	19.4196307214216\\
86.4	19.5364110635889\\
86.5	19.5800569068633\\
86.6	19.6048687523349\\
86.7	19.4117942548658\\
86.8	19.5482030219377\\
86.9	19.3944468949136\\
87	19.4345418848227\\
87.1	19.3472528925696\\
87.2	19.7626241464629\\
87.3	19.643311960329\\
87.4	19.9373144188093\\
87.5	19.7201185867368\\
87.6	19.6517721089637\\
87.7	19.6453518762013\\
87.8	19.7399356337387\\
87.9	19.5696229686119\\
88	19.8996704447858\\
88.1	19.837667099353\\
88.2	19.5268733632336\\
88.3	19.8846262728527\\
88.4	19.8889323642903\\
88.5	19.7494630127164\\
88.6	19.9220559718402\\
88.7	20.0447169146466\\
88.8	20.1303139584168\\
88.9	20.3211638666083\\
89	20.0357857340143\\
89.1	20.1946373420825\\
89.2	19.9704543707776\\
89.3	20.1423760342389\\
89.4	19.9235146936377\\
89.5	19.9698716439976\\
89.6	20.0892675029422\\
89.7	20.0631544332054\\
89.8	20.1012058448157\\
89.9	20.2636204214814\\
90	20.1229053638908\\
90.1	20.1597042156756\\
90.2	20.324132578919\\
90.3	20.2673542917381\\
90.4	20.4174932132008\\
90.5	20.0518126096432\\
90.6	20.2422493703769\\
90.7	20.3466428210611\\
90.8	20.463709121785\\
90.9	20.2481082190932\\
91	20.4045521984357\\
91.1	20.3457934811475\\
91.2	20.8400822716374\\
91.3	20.782750032084\\
91.4	20.6142756132473\\
91.5	20.3979883455129\\
91.6	20.5617885646761\\
91.7	20.5524182761588\\
91.8	20.4701640776939\\
91.9	20.6154888775028\\
92	20.6542938123393\\
92.1	21.0310970530852\\
92.2	20.6972406879092\\
92.3	20.8637341345383\\
92.4	20.6768486061656\\
92.5	20.5017452562693\\
92.6	20.5594677155405\\
92.7	20.62746748772\\
92.8	20.7300522706165\\
92.9	20.880340901263\\
93	20.7457213747335\\
93.1	20.9079566921942\\
93.2	20.9411718339806\\
93.3	20.7900043321297\\
93.4	20.6384736532395\\
93.5	21.1299963561315\\
93.6	21.097533903969\\
93.7	20.9255480532589\\
93.8	20.8675918065667\\
93.9	21.1091626285181\\
94	21.086019918814\\
94.1	21.2571624706991\\
94.2	21.2109659448457\\
94.3	21.3789882986563\\
94.4	21.4769991350158\\
94.5	21.2408836360139\\
94.6	21.3207831917885\\
94.7	21.2313008740597\\
94.8	21.3523644523057\\
94.9	21.3002673535836\\
95	21.2168528618833\\
95.1	21.5218204439778\\
95.2	21.234605998824\\
95.3	21.5340784952893\\
95.4	21.4112985001799\\
95.5	21.269885032239\\
95.6	21.45356730638\\
95.7	21.3920509020535\\
95.8	21.4040469029043\\
95.9	21.3335555177185\\
96	21.6108677282668\\
96.1	21.5706912354308\\
96.2	21.5341781080267\\
96.3	21.6015847789013\\
96.4	21.6577520970836\\
96.5	21.8735168963057\\
96.6	21.7304419038627\\
96.7	21.9475057251209\\
96.8	21.7304784843724\\
96.9	21.6859905863257\\
97	21.9075023511828\\
97.1	21.5857705783428\\
97.2	22.1490531415719\\
97.3	21.7149364583976\\
97.4	21.8792121656568\\
97.5	21.9341330147674\\
97.6	21.8103868759285\\
97.7	22.0303303404532\\
97.8	21.8642621473005\\
97.9	21.9624317727434\\
98	22.0217137845544\\
98.1	22.3122843937143\\
98.2	22.1551122597954\\
98.3	22.0707551966822\\
98.4	21.9392532904665\\
98.5	22.171636198776\\
98.6	22.1620582345639\\
98.7	22.1553449443948\\
98.8	22.2479517163293\\
98.9	22.1676581932832\\
99	22.2617840489201\\
99.1	22.0979024131272\\
99.2	22.2856431846868\\
99.3	22.3727221533971\\
99.4	22.1883466072367\\
99.5	22.2629078234918\\
99.6	22.2829115953145\\
99.7	22.5234883286981\\
99.8	22.4924268763711\\
99.9	22.4125083509405\\
100	22.4431568055797\\
};
\addlegendentry{$D_{\text{KY}}$};

\addplot [color=mycolor2,dashed,line width=1pt]
  table[row sep=crcr]{%
0.1	-0.1374\\
100	22.44\\
};
\addlegendentry{$0.226L-0.160$};

\end{axis}
\end{tikzpicture}%

%% file: ksoddperiodic_dky.tex
\definecolor{mycolor1}{rgb}{0.00000,0.44700,0.74100}%
\definecolor{mycolor2}{rgb}{0.85000,0.32500,0.09800}%
\definecolor{mycolor3}{rgb}{0.92900,0.69400,0.12500}%
\definecolor{mycolor4}{rgb}{0.49400,0.18400,0.55600}%
\definecolor{mycolor5}{rgb}{0.46600,0.67400,0.18800}%
\definecolor{mycolor6}{rgb}{0.30100,0.74500,0.93300}%
\definecolor{mycolor7}{rgb}{0.63500,0.07800,0.18400}%
\begin{tikzpicture}

\begin{axis}[%
/tikz/line join=bevel,
width=4.25in,
height=1.75in,
at={(0in,0in)},
scale only axis,
xmin=0,
xmax=100,
xlabel={$L$},
xlabel shift={-1pt},
xmajorgrids,
ymin=0,
ymax=23,
ymajorgrids,
axis background/.style={fill=white},
legend style={at={(0.02,0.96)},anchor=north west,legend cell align=left,align=left,draw=white!15!black}
]
\addplot [color=mycolor1,solid,line width=1pt]
  table[row sep=crcr]{%
0.1	0\\
0.2	0\\
0.3	0\\
0.4	0\\
0.5	0\\
0.6	0\\
0.7	0\\
0.8	0\\
0.9	0\\
1	0\\
1.1	0\\
1.2	0\\
1.3	0\\
1.4	0\\
1.5	0\\
1.6	0\\
1.7	0\\
1.8	0\\
1.9	0\\
2	0\\
2.1	0\\
2.2	0\\
2.3	0\\
2.4	0\\
2.5	0\\
2.6	0\\
2.7	0\\
2.8	0\\
2.9	0\\
3	0\\
3.1	0\\
3.2	0\\
3.3	0\\
3.4	0\\
3.5	0\\
3.6	0\\
3.7	0\\
3.8	0\\
3.9	0\\
4	0\\
4.1	0\\
4.2	0\\
4.3	0\\
4.4	0\\
4.5	0\\
4.6	0\\
4.7	0\\
4.8	0\\
4.9	0\\
5	0\\
5.1	0\\
5.2	0\\
5.3	0\\
5.4	0\\
5.5	0\\
5.6	0\\
5.7	0\\
5.8	0\\
5.9	0\\
6	0\\
6.1	0\\
6.2	0\\
6.3	0\\
6.4	0\\
6.5	0\\
6.6	0\\
6.7	0\\
6.8	0\\
6.9	0\\
7	0\\
7.1	0\\
7.2	0\\
7.3	0\\
7.4	0\\
7.5	0\\
7.6	0\\
7.7	0\\
7.8	0\\
7.9	0\\
8	0\\
8.1	0\\
8.2	0\\
8.3	0\\
8.4	0\\
8.5	0\\
8.6	0\\
8.7	0\\
8.8	0\\
8.9	0\\
9	0\\
9.1	0\\
9.2	0\\
9.3	0\\
9.4	0\\
9.5	0\\
9.6	0\\
9.7	0\\
9.8	0\\
9.9	0\\
10	0\\
10.1	0\\
10.2	0\\
10.3	0\\
10.4	0\\
10.5	0\\
10.6	0\\
10.7	0\\
10.8	0\\
10.9	0\\
11	0\\
11.1	0\\
11.2	0\\
11.3	0\\
11.4	0\\
11.5	0\\
11.6	0\\
11.7	0\\
11.8	0\\
11.9	0\\
12	0\\
12.1	0\\
12.2	0\\
12.3	0\\
12.4	0\\
12.5	0\\
12.6	0\\
12.7	0\\
12.8	0\\
12.9	0\\
13	0\\
13.1	0\\
13.2	0\\
13.3	0\\
13.4	0\\
13.5	0\\
13.6	0\\
13.7	0\\
13.8	0\\
13.9	0\\
14	0\\
14.1	0\\
14.2	0\\
14.3	0\\
14.4	0\\
14.5	0\\
14.6	0\\
14.7	0\\
14.8	0\\
14.9	0\\
15	0\\
15.1	0\\
15.2	0\\
15.3	0\\
15.4	0\\
15.5	0\\
15.6	0\\
15.7	0\\
15.8	1.00609002989608\\
15.9	1.00901030365492\\
16	0\\
16.1	0\\
16.2	0\\
16.3	0\\
16.4	0\\
16.5	0\\
16.6	0\\
16.7	0\\
16.8	0\\
16.9	0\\
17	0\\
17.1	1.01545177552025\\
17.2	1.00084502418232\\
17.3	1.00805833546518\\
17.4	0\\
17.5	0\\
17.6	0\\
17.7	0\\
17.8	0\\
17.9	0\\
18	2.15078552880912\\
18.1	1.08143155140301\\
18.2	2.48176693278392\\
18.3	2.50202382137627\\
18.4	2.63437955677416\\
18.5	2.44552136810086\\
18.6	2.55886523793046\\
18.7	2.60582021238393\\
18.8	2.60455019773915\\
18.9	3.00358395086749\\
19	2.64393778087214\\
19.1	2.65122423351387\\
19.2	2.57552828057026\\
19.3	2.41269902240096\\
19.4	2.48831681057923\\
19.5	2.38980764544564\\
19.6	2.32167703389222\\
19.7	2.2646485201476\\
19.8	2.26466630343056\\
19.9	0\\
20	2.16961608140721\\
20.1	0\\
20.2	0\\
20.3	0\\
20.4	2.06914432561611\\
20.5	0\\
20.6	1.35026837934331\\
20.7	0\\
20.8	0\\
20.9	0\\
21	2.35428166336516\\
21.1	2.05214116127795\\
21.2	1.11156400865386\\
21.3	2.7212924769947\\
21.4	2.77192589457318\\
21.5	2.75385075946545\\
21.6	2.6953223224454\\
21.7	2.9153609743713\\
21.8	3.1765915543492\\
21.9	3.21290778758764\\
22	3.26384167595636\\
22.1	3.34545850058398\\
22.2	0\\
22.3	0\\
22.4	0\\
22.5	0\\
22.6	0\\
22.7	0\\
22.8	0\\
22.9	0\\
23	0\\
23.1	0\\
23.2	0\\
23.3	0\\
23.4	0\\
23.5	0\\
23.6	0\\
23.7	0\\
23.8	0\\
23.9	0\\
24	0\\
24.1	0\\
24.2	0\\
24.3	0\\
24.4	0\\
24.5	0\\
24.6	0\\
24.7	0\\
24.8	0\\
24.9	0\\
25	0\\
25.1	0\\
25.2	0\\
25.3	0\\
25.4	0\\
25.5	0\\
25.6	0\\
25.7	2.94108871410032\\
25.8	3.56951386954763\\
25.9	3.60088245224116\\
26	3.37591823049873\\
26.1	3.70047326793854\\
26.2	3.63376997224155\\
26.3	3.99984919262586\\
26.4	3.64836730192146\\
26.5	0\\
26.6	3.91458399409788\\
26.7	0\\
26.8	0\\
26.9	0\\
27	0\\
27.1	0\\
27.2	0\\
27.3	0\\
27.4	0\\
27.5	0\\
27.6	0\\
27.7	0\\
27.8	0\\
27.9	0\\
28	0\\
28.1	0\\
28.2	0\\
28.3	0\\
28.4	0\\
28.5	1.07164937381528\\
28.6	0\\
28.7	1.04379071314611\\
28.8	0\\
28.9	1.19866536992885\\
29	0\\
29.1	1.02758887747636\\
29.2	0\\
29.3	1.03854550250875\\
29.4	0\\
29.5	0\\
29.6	2.06062739669434\\
29.7	0\\
29.8	0\\
29.9	4.50755538586846\\
30	4.49225810849172\\
30.1	4.48429652472513\\
30.2	0\\
30.3	4.68763665650443\\
30.4	0\\
30.5	0\\
30.6	0\\
30.7	0\\
30.8	4.70557385284789\\
30.9	0\\
31	0\\
31.1	0\\
31.2	0\\
31.3	4.6962443452083\\
31.4	0\\
31.5	0\\
31.6	0\\
31.7	0\\
31.8	0\\
31.9	0\\
32	0\\
32.1	0\\
32.2	0\\
32.3	0\\
32.4	0\\
32.5	0\\
32.6	0\\
32.7	1.14061160217125\\
32.8	4.51166638432792\\
32.9	4.42647785495601\\
33	4.38829357183265\\
33.1	4.96069902926025\\
33.2	4.79729675017652\\
33.3	5.0566928236404\\
33.4	5.05970163907216\\
33.5	5.06389235367318\\
33.6	5.03943442178005\\
33.7	5.03925063959942\\
33.8	5.15654544379038\\
33.9	5.05443852275046\\
34	4.84629021138289\\
34.1	5.30639673872379\\
34.2	3.99438681617234\\
34.3	5.00805896760702\\
34.4	5.18534557381448\\
34.5	5.24102691871296\\
34.6	5.23897471314234\\
34.7	4.98481663260219\\
34.8	5.46292404616094\\
34.9	5.49094325848763\\
35	5.47747726979492\\
35.1	5.3975018281657\\
35.2	5.2005305385961\\
35.3	5.45215664539069\\
35.4	5.61566000460361\\
35.5	5.61477678990473\\
35.6	5.56232201164478\\
35.7	5.33259049250269\\
35.8	5.16497375820906\\
35.9	4.98658168580306\\
36	5.4489126291967\\
36.1	0\\
36.2	0\\
36.3	5.84110115041016\\
36.4	5.93230338118154\\
36.5	5.75697873624264\\
36.6	0\\
36.7	2.18252951926396\\
36.8	3.79001995299835\\
36.9	6.16589814578812\\
37	5.9444534900126\\
37.1	0\\
37.2	0\\
37.3	0\\
37.4	5.56451959699602\\
37.5	1.08692774301236\\
37.6	5.65590061264213\\
37.7	6.1648681795781\\
37.8	6.31515788075949\\
37.9	4.02212171035962\\
38	6.36102100202161\\
38.1	6.31283662806405\\
38.2	4.00080681353464\\
38.3	6.2188189328226\\
38.4	6.30427073954928\\
38.5	6.3728594004916\\
38.6	6.36096952158537\\
38.7	6.40444544896808\\
38.8	6.33935274708096\\
38.9	0\\
39	0\\
39.1	5.84801124513443\\
39.2	6.57888239834987\\
39.3	0\\
39.4	6.65333308149449\\
39.5	1.33337427398013\\
39.6	0\\
39.7	2.01608077000084\\
39.8	6.21919376323752\\
39.9	0\\
40	0\\
40.1	6.82350091613149\\
40.2	6.65706435050152\\
40.3	0\\
40.4	0\\
40.5	0\\
40.6	0\\
40.7	4.32001748021159\\
40.8	6.77267591264641\\
40.9	6.77338539155206\\
41	7.05578952958086\\
41.1	7.13046229386329\\
41.2	6.17009451716748\\
41.3	0\\
41.4	6.99658339226028\\
41.5	6.9414917624613\\
41.6	7.07560872470405\\
41.7	7.20176022496251\\
41.8	7.23463010387951\\
41.9	7.07441540482053\\
42	2.8400776975106\\
42.1	7.25149450460243\\
42.2	7.23801686990303\\
42.3	7.24057511154384\\
42.4	7.35236370799577\\
42.5	7.43743615975039\\
42.6	7.37385565937448\\
42.7	7.43294473889513\\
42.8	0\\
42.9	7.41809698274845\\
43	7.39315727673594\\
43.1	7.4057534838433\\
43.2	0\\
43.3	7.45290055243222\\
43.4	7.6005365889314\\
43.5	7.38787079074813\\
43.6	7.6494389702141\\
43.7	7.69105723385032\\
43.8	7.57316440759022\\
43.9	7.53683067818937\\
44	7.60324331669563\\
44.1	7.69597139724142\\
44.2	7.59032426105281\\
44.3	7.56805936853846\\
44.4	7.7125853567378\\
44.5	7.67373118423151\\
44.6	7.80668385201232\\
44.7	7.53628963912516\\
44.8	7.85797352511688\\
44.9	7.85752113080293\\
45	7.71198683329874\\
45.1	7.9246702449634\\
45.2	7.97795261424469\\
45.3	8.06130100725057\\
45.4	7.96646163285704\\
45.5	8.16395019045745\\
45.6	7.94374536241266\\
45.7	8.02995824771112\\
45.8	7.96333050812403\\
45.9	8.21743088112922\\
46	8.03518042455627\\
46.1	7.43165678219776\\
46.2	8.07851142241345\\
46.3	8.15555346868833\\
46.4	8.12573146361677\\
46.5	8.11596323576707\\
46.6	8.40419528768612\\
46.7	8.28901982163386\\
46.8	8.31128471423749\\
46.9	8.44572528293589\\
47	8.35373772590756\\
47.1	8.23104039102616\\
47.2	8.41526813744966\\
47.3	8.40665568464302\\
47.4	8.58030851620815\\
47.5	8.4633547046062\\
47.6	8.45678671227947\\
47.7	8.53805560150334\\
47.8	8.38308905019888\\
47.9	0\\
48	8.48517615013561\\
48.1	8.65295915944955\\
48.2	0\\
48.3	8.29128844203873\\
48.4	8.78449529047436\\
48.5	0\\
48.6	8.94410915543505\\
48.7	8.33460300685199\\
48.8	8.49237662772386\\
48.9	8.89077220484912\\
49	8.01730923346661\\
49.1	8.60968091019268\\
49.2	8.92755791718938\\
49.3	8.94833295064004\\
49.4	8.8448427923843\\
49.5	8.83818055505513\\
49.6	8.92875386474469\\
49.7	9.00169722636269\\
49.8	9.1005101997274\\
49.9	9.05441244375235\\
50	9.05030398274761\\
50.1	9.1145096964386\\
50.2	9.12909511364711\\
50.3	9.07638793196362\\
50.4	9.13244978891816\\
50.5	9.13648980372263\\
50.6	9.32571439879824\\
50.7	9.16689577109865\\
50.8	9.2787783754713\\
50.9	9.31353994181985\\
51	9.27910147467766\\
51.1	9.26465631890751\\
51.2	9.39054183998724\\
51.3	9.35062145274785\\
51.4	9.40159297177296\\
51.5	0\\
51.6	9.3482114384803\\
51.7	9.66072940540219\\
51.8	9.44757819503265\\
51.9	9.66986908895009\\
52	9.62376466133002\\
52.1	9.38696926655558\\
52.2	9.53922256641221\\
52.3	9.54672752307259\\
52.4	9.35654390629459\\
52.5	9.70565021297881\\
52.6	9.46125107527104\\
52.7	8.69354440637861\\
52.8	9.728007435245\\
52.9	9.82966343970192\\
53	9.40769485035785\\
53.1	9.74618084421686\\
53.2	9.84863457990154\\
53.3	9.71259534092147\\
53.4	10.0047433822259\\
53.5	10.0215067438497\\
53.6	9.65140471869532\\
53.7	9.83501244558952\\
53.8	10.1544105799018\\
53.9	10.0431248711089\\
54	10.0854859104687\\
54.1	9.93162985063359\\
54.2	9.97476718385665\\
54.3	10.0893214489705\\
54.4	10.0086794962364\\
54.5	10.1040942162976\\
54.6	10.1898817810722\\
54.7	10.1379668913549\\
54.8	10.059755332486\\
54.9	0\\
55	10.3329787752601\\
55.1	10.18516826533\\
55.2	10.4189236394464\\
55.3	10.2010166468494\\
55.4	0\\
55.5	10.3082332968093\\
55.6	10.3315156562215\\
55.7	10.3724564933809\\
55.8	10.3902310172736\\
55.9	10.4043681412134\\
56	10.4654365668206\\
56.1	10.2964398865974\\
56.2	10.4963100373998\\
56.3	10.3816761872165\\
56.4	10.5575903119541\\
56.5	10.6012771245736\\
56.6	10.6149809822942\\
56.7	10.7097919199436\\
56.8	1.09080134959954\\
56.9	10.618544124725\\
57	10.5546645729284\\
57.1	10.6644321303619\\
57.2	10.6624181411811\\
57.3	10.5984121357754\\
57.4	10.7525108847757\\
57.5	10.6576958845866\\
57.6	10.8359776164411\\
57.7	10.8400946492536\\
57.8	10.608638468865\\
57.9	11.1223697605557\\
58	11.1937429721325\\
58.1	10.7516675748457\\
58.2	10.8427616789832\\
58.3	10.997743110109\\
58.4	11.1620468549049\\
58.5	10.9949391689924\\
58.6	11.041028671544\\
58.7	10.8513325848085\\
58.8	11.0952967153191\\
58.9	11.2171168190968\\
59	11.1077125886064\\
59.1	11.1423740078575\\
59.2	11.2651945850247\\
59.3	11.2309332661687\\
59.4	11.2767387880643\\
59.5	11.4015045478422\\
59.6	11.1995011721853\\
59.7	11.2370052938091\\
59.8	11.3140805075255\\
59.9	11.278907192949\\
60	11.3465022433452\\
60.1	11.3857643914231\\
60.2	11.4431819407331\\
60.3	11.3684990221402\\
60.4	11.3814069489307\\
60.5	11.453098926282\\
60.6	11.2985073510803\\
60.7	11.4708023276645\\
60.8	11.4506591389123\\
60.9	11.6195696067799\\
61	11.6807353828235\\
61.1	11.7054992502012\\
61.2	11.5840173351964\\
61.3	11.6554445295895\\
61.4	11.5659322791682\\
61.5	11.6816402370173\\
61.6	11.7497395847462\\
61.7	11.7230079569239\\
61.8	11.6839866079669\\
61.9	11.7980065845405\\
62	11.8829042668941\\
62.1	11.9329311580179\\
62.2	11.9219220327213\\
62.3	11.8426955846748\\
62.4	11.8562088263529\\
62.5	11.8080637641363\\
62.6	11.8490483650973\\
62.7	11.9237986978096\\
62.8	12.1907354718673\\
62.9	12.098623260062\\
63	11.958914840216\\
63.1	12.0118942900971\\
63.2	12.028051932398\\
63.3	12.2465800540626\\
63.4	12.0990494749725\\
63.5	12.153946318024\\
63.6	12.1807411272954\\
63.7	12.247181570775\\
63.8	12.1339989518605\\
63.9	12.2084787798372\\
64	12.34652944802\\
64.1	12.3775020600999\\
64.2	12.3104382825098\\
64.3	12.5710746496726\\
64.4	12.4151932352696\\
64.5	12.3846549759157\\
64.6	12.3397708438446\\
64.7	12.4176881804867\\
64.8	12.2172100706686\\
64.9	12.5434724029539\\
65	12.5628550551539\\
65.1	12.4616085750663\\
65.2	12.3992965405931\\
65.3	12.4979652938257\\
65.4	12.7011149121496\\
65.5	12.414919402153\\
65.6	12.6792720765953\\
65.7	12.5954747405811\\
65.8	12.5646013599686\\
65.9	12.5624164759398\\
66	12.6533856293495\\
66.1	13.0733062014304\\
66.2	12.8373803242974\\
66.3	12.6861145434773\\
66.4	12.9965483981734\\
66.5	12.7308142478325\\
66.6	12.9366367399303\\
66.7	12.9365499047448\\
66.8	13.0950998078991\\
66.9	12.8356757434463\\
67	12.7405874207616\\
67.1	12.8559271438845\\
67.2	12.9801021804295\\
67.3	13.0190712152751\\
67.4	13.1130491915932\\
67.5	13.2254538069005\\
67.6	13.0047462416971\\
67.7	13.0161775298424\\
67.8	13.1397023893056\\
67.9	13.0837539287215\\
68	13.1201090132576\\
68.1	13.3159169364482\\
68.2	13.2589412473462\\
68.3	13.3009907662509\\
68.4	13.3393023084278\\
68.5	13.2447770883036\\
68.6	13.3592447026636\\
68.7	13.4613814897365\\
68.8	13.4694537933514\\
68.9	13.5309761345762\\
69	13.3112437699432\\
69.1	13.5840728368591\\
69.2	13.5181831603483\\
69.3	13.2994291819294\\
69.4	13.4389152514144\\
69.5	13.5268280193957\\
69.6	13.3947327939443\\
69.7	13.2161154159367\\
69.8	13.515456448916\\
69.9	13.5956361660012\\
70	13.5519271678234\\
70.1	13.5608776601818\\
70.2	13.7340878510066\\
70.3	13.6585428172668\\
70.4	13.6965055752096\\
70.5	13.8107026558822\\
70.6	13.8516993934045\\
70.7	13.6474494069951\\
70.8	13.8588795766845\\
70.9	13.7851153577535\\
71	14.0310741636469\\
71.1	13.9265009479757\\
71.2	13.9721304328388\\
71.3	13.5931981086215\\
71.4	13.7856098113332\\
71.5	13.9859073617633\\
71.6	14.1164117253706\\
71.7	14.0508305167487\\
71.8	14.1332409183773\\
71.9	14.0167969637752\\
72	14.1066582012752\\
72.1	14.1347672165758\\
72.2	14.432000971016\\
72.3	14.4523360261688\\
72.4	14.1907047876674\\
72.5	14.138863672984\\
72.6	14.1008566894147\\
72.7	14.206177639764\\
72.8	14.2335604238726\\
72.9	14.1665669206197\\
73	14.371528305669\\
73.1	14.3812122687758\\
73.2	14.3466906437622\\
73.3	14.3779974149227\\
73.4	14.3205369530069\\
73.5	14.723903925596\\
73.6	14.5177547697201\\
73.7	14.6139934018137\\
73.8	14.5264616252034\\
73.9	14.4318237834787\\
74	14.6669590227209\\
74.1	14.5458605684727\\
74.2	14.6303712657277\\
74.3	14.6179603190841\\
74.4	14.6641345959644\\
74.5	14.6017876186556\\
74.6	14.8001119158461\\
74.7	14.6790459401355\\
74.8	15.1134021822646\\
74.9	14.8500577802969\\
75	14.6112670273347\\
75.1	15.0319605778394\\
75.2	14.543155405137\\
75.3	14.7855164877733\\
75.4	15.0316088129545\\
75.5	14.7825109711478\\
75.6	14.8812470423307\\
75.7	14.8027566527637\\
75.8	14.9147080550688\\
75.9	14.89478130902\\
76	14.9685241544024\\
76.1	15.1199634668932\\
76.2	15.2748521086697\\
76.3	14.9560420940055\\
76.4	15.0117751384982\\
76.5	15.2115229778685\\
76.6	15.3117260509569\\
76.7	15.2078342923254\\
76.8	15.1155016025038\\
76.9	15.3004402653516\\
77	15.301868340192\\
77.1	15.4290886488138\\
77.2	15.4281546288476\\
77.3	15.3195951512399\\
77.4	15.3212721765965\\
77.5	15.2003206866056\\
77.6	15.0842814364172\\
77.7	15.4094572003334\\
77.8	15.4199624634087\\
77.9	15.5199594572371\\
78	15.5738086712004\\
78.1	15.4296819003749\\
78.2	15.5848615221463\\
78.3	15.2776549553053\\
78.4	15.4189787799088\\
78.5	15.694820731035\\
78.6	15.5886745636354\\
78.7	15.6045178789978\\
78.8	15.6341698596016\\
78.9	15.6309432165466\\
79	15.5402659535449\\
79.1	15.7869560066943\\
79.2	15.5618256288787\\
79.3	16.0127087463707\\
79.4	15.7524762211358\\
79.5	15.833043405339\\
79.6	15.9327321621424\\
79.7	15.8473389958814\\
79.8	16.1788590094583\\
79.9	15.9448628947939\\
80	15.8537841834995\\
80.1	0\\
80.2	15.9648336239567\\
80.3	15.9913037947155\\
80.4	16.2462771900239\\
80.5	15.9984989967461\\
80.6	16.1141585761862\\
80.7	16.1433968455235\\
80.8	15.9532702181908\\
80.9	16.1952110644248\\
81	16.4110664987715\\
81.1	16.222862646071\\
81.2	16.3009463717436\\
81.3	15.991538713573\\
81.4	16.2564757277526\\
81.5	16.2202432120913\\
81.6	16.4573031897551\\
81.7	16.17393821693\\
81.8	16.0726060417098\\
81.9	16.438354427705\\
82	16.3833468081425\\
82.1	16.3962365612704\\
82.2	16.4748267170476\\
82.3	16.4575084308946\\
82.4	16.5220833573845\\
82.5	16.5762078136467\\
82.6	16.4091970676235\\
82.7	16.6907201664662\\
82.8	16.4857333998647\\
82.9	16.7227428245481\\
83	16.6390868822743\\
83.1	16.822543855067\\
83.2	16.6220958203246\\
83.3	16.6566198241808\\
83.4	16.715079269143\\
83.5	16.99477715553\\
83.6	16.5721150592872\\
83.7	16.8900291272798\\
83.8	16.7858371505171\\
83.9	16.7395009691002\\
84	16.7402583121892\\
84.1	16.7770973493026\\
84.2	17.0224238421295\\
84.3	16.9252177805313\\
84.4	16.9991880417792\\
84.5	17.1223338179515\\
84.6	17.0569924026104\\
84.7	16.7914393999838\\
84.8	16.9035939515555\\
84.9	17.0775203296031\\
85	16.9611937843089\\
85.1	17.1160666382473\\
85.2	17.1864489172905\\
85.3	17.2154295322642\\
85.4	17.1780231205403\\
85.5	17.2206848926317\\
85.6	17.193980147114\\
85.7	17.2867211877261\\
85.8	17.3882364254094\\
85.9	17.1645363781468\\
86	17.4696267532188\\
86.1	17.1763602027727\\
86.2	17.1908229390992\\
86.3	17.5652089707089\\
86.4	17.4161494130499\\
86.5	17.1736836866351\\
86.6	17.5598766488398\\
86.7	17.4592024140289\\
86.8	17.3554558646576\\
86.9	17.352386361922\\
87	17.6187502126932\\
87.1	17.4670309121456\\
87.2	17.5895031661011\\
87.3	17.4921114750791\\
87.4	17.6513583426721\\
87.5	17.7404282956596\\
87.6	17.6670758514769\\
87.7	17.6233261706378\\
87.8	17.8310149149324\\
87.9	17.7488622692803\\
88	17.879794770259\\
88.1	17.6890662291362\\
88.2	17.7637277482366\\
88.3	17.7191220368847\\
88.4	17.9790886275307\\
88.5	17.8952526531515\\
88.6	17.6774641929997\\
88.7	17.8343412704123\\
88.8	17.8386765522432\\
88.9	17.9228132120163\\
89	18.0640671694887\\
89.1	18.3021562660381\\
89.2	17.8976895927042\\
89.3	18.1132943363539\\
89.4	18.0610203898957\\
89.5	18.145748337803\\
89.6	17.8375375597439\\
89.7	17.9985423234295\\
89.8	18.1699535085627\\
89.9	18.1754039342993\\
90	18.1005070733287\\
90.1	18.2554256265556\\
90.2	18.2163626727339\\
90.3	18.3831503768651\\
90.4	18.3399312315425\\
90.5	18.4677665269946\\
90.6	18.4740047711431\\
90.7	18.3497568111673\\
90.8	18.4605028505504\\
90.9	18.4607535125511\\
91	18.3568834683023\\
91.1	18.5179713606543\\
91.2	18.5169747431426\\
91.3	18.5033843820503\\
91.4	18.5639354165154\\
91.5	18.3436488849719\\
91.6	18.9313312047129\\
91.7	18.6406302730842\\
91.8	18.8316165694007\\
91.9	18.736703016093\\
92	18.3432813787846\\
92.1	18.8679410791928\\
92.2	18.7156116137999\\
92.3	18.6627881883123\\
92.4	18.7423205872015\\
92.5	18.8116981328314\\
92.6	18.6186840929025\\
92.7	18.9732464746527\\
92.8	18.9194025250001\\
92.9	18.9450761430211\\
93	18.8120462660582\\
93.1	19.0284412255986\\
93.2	18.9803156708327\\
93.3	19.1560464211672\\
93.4	18.7923961844238\\
93.5	19.013310116751\\
93.6	18.9822637265397\\
93.7	19.0446296038094\\
93.8	19.166613349276\\
93.9	19.1328618983165\\
94	19.0495071645482\\
94.1	19.1084545096176\\
94.2	19.0685586185415\\
94.3	19.2797478576158\\
94.4	19.2156631983606\\
94.5	19.0878987602593\\
94.6	19.2216042797006\\
94.7	19.21919848799\\
94.8	19.2769649009552\\
94.9	19.3274671234795\\
95	19.1400830450964\\
95.1	19.4971574595148\\
95.2	19.327226292793\\
95.3	19.279165826455\\
95.4	19.1552332515284\\
95.5	19.4774229923201\\
95.6	19.4924551893393\\
95.7	19.2902512610154\\
95.8	19.4118412240524\\
95.9	19.0833764260723\\
96	19.4340656732252\\
96.1	19.6178123281059\\
96.2	19.8226331388405\\
96.3	19.6179248513119\\
96.4	19.6236704268103\\
96.5	19.5865267638938\\
96.6	19.8473764739149\\
96.7	19.5472193430915\\
96.8	19.9411682730967\\
96.9	19.9384280676661\\
97	19.7351682401654\\
97.1	19.7348244137389\\
97.2	19.8385453366375\\
97.3	19.490562365512\\
97.4	19.8293565040603\\
97.5	19.6620620589358\\
97.6	20.0883028681827\\
97.7	20.0462955041151\\
97.8	20.0936178753522\\
97.9	19.9685744687508\\
98	19.9596319521676\\
98.1	20.1318881046667\\
98.2	19.8689029919308\\
98.3	20.1803890168002\\
98.4	19.9656189365087\\
98.5	20.3035331948659\\
98.6	20.3807588690177\\
98.7	20.2946963717796\\
98.8	20.2436477041572\\
98.9	20.276494929122\\
99	20.018037054308\\
99.1	20.1374572965341\\
99.2	20.4940660469471\\
99.3	19.9482183964311\\
99.4	20.4294716840382\\
99.5	20.4025954467486\\
99.6	20.6012123055866\\
99.7	20.3444215956578\\
99.8	20.3843447631953\\
99.9	20.3374674538082\\
100	20.7535952579485\\
};
\addlegendentry{$D_{\text{KY}}$};

\addplot [color=mycolor2,dashed,line width=1pt]
  table[row sep=crcr]{%
0.1	-2.0834\\
100	20.494\\
};
\addlegendentry{$0.226L-2.106$};

\end{axis}
\end{tikzpicture}%